\documentclass{amsart}
\usepackage{amsmath,amssymb,amsthm}
\usepackage
{graphicx}
\usepackage{epstopdf,pdfsync}
\usepackage{ifthen}

\newtheorem{thm}{Theorem}[section]
\newtheorem{lem}{Lemma}[section]
\newtheorem{cor}{Corollary}[section]

\newtheorem{defn}{Definition}[section]
\newtheorem{exmp}{Example}[section]
\newtheorem{rem}{Remark}[section]

\numberwithin{equation}{section}
\def\Z{\Bbb Z}

\def\R{\Bbb R}

\def\d{\partial}

\def\e{\epsilon}

\newcommand{\showcomments}{no}

\newcommand{\commentstyle}{\tiny}

\newcommand{\comment}[1]
{\ifthenelse{\equal{\showcomments}{yes}}
{\footnotemark\marginpar{\sffamily{\commentstyle
\addtocounter{footnote}{-1}\footnotemark#1 }\normalfont}}{}}

\def\flabel#1{\ifmmode #1\else$ #1$\fi}

\let\angle\undefined
\newsymbol\angle 115C
\def\degrees{\ifmmode^\circ\else$^\circ$\fi}
\def\a{\alpha}
\def\b{\beta}
\def\g{\gamma}

\def \pict #1 by #2 (#3) {\centerline{
\vbox to #2 {\hrule width #1 height 0pt depth 0pt
\vfill{\special{picture #3 }}}}}

\def \picture #1 by #2 (#3 scaled #4) #5{
\dimen0=#1 \dimen1=#2
\divide\dimen0 by 1000 \multiply\dimen0 by #4
\divide\dimen1 by 1000 \multiply\dimen1 by #4
\vbox{\pict \dimen0 by \dimen1 (#3 scaled #4)
\centerline { #5}}}

\catcode`\@=11

\title[]{Convexity of Morse Stratifications and Gradient Spines of $3$-manifolds}
\begin{document}
\author{Gabriel Katz}
\address{Department of Mathematics, William Paterson University, Wayne, NJ  
07470-2103}
\email{katzg@wpunj.edu}
\date{\today}
\maketitle
\qquad \qquad \qquad \emph{In memory of Jerry Levine,  friend and  mentor}

\begin{abstract} We notice that a generic nonsingular gradient field $v = \nabla f$  on a compact 3-fold $X$ with boundary canonically generates a simple  spine $K(f, v)$ of $X$. We study the transformations of  $K(f, v)$ that are induced by deformations of the data $(f, v)$. We link the Matveev complexity $c(X)$ of $X$ with counting the \emph{double-tangent} trajectories of the $v$-flow, i.e. the trajectories that are tangent to the boundary $\d X$ at a pair of distinct points. Let $gc(X)$ be the minimum number of such trajectories, minimum being taken over all nonsingular $v$'s.  We call $gc(X)$ the \emph{gradient complexity} of $X$.  Next, we prove that there are only finitely many $X$  of bounded gradient complexity, provided that $X$ is irreducible and boundary irreducible  with no essential annuli.  In particular, there exists only finitely many hyperbolic manifolds $X$ with bounded $gc(X)$. For such $X$, their normalized hyperbolic volume gives an upper bound  of $gc(X)$. If an irreducible and boundary irreducible $X$ with no essential annuli admits a nonsingular gradient flow with no double-tangent trajectories, then $X$ is a standard ball. All these and many other results of the paper rely on a careful study of the stratified geometry of $\d X$ relative to the $v$-flow. It is characterized by failure of $\d X$ to be \emph{convex} with respect to a generic flow $v$. It turns out, that convexity or its lack have profound influence on the topology of $X$. This interplay between intrinsic concavity of $\d X$ with respect to any gradient-like flow and complexity of $X$ is in the focus of the paper.
\end{abstract}

\bigskip

\section{Introduction}

Classical Morse theory links singularities of Morse functions  with   
topology of closed manifolds. Specifically, singularities of Morse functions $f: X \to \R$ cause  interruptions of the $f$-gradient flow, and the homology or even the topological type of a manifold  $X$ can  be expressed in terms of such interruptions (see [C]). These terms include descending and ascending disks, attaching maps, and spaces of  flow trajectories which connect the singularities. 

On manifolds with boundary, an additional  source of  the flow interruption occurs: it comes from a particular geometry of the boundary $\d X$, or rather  from the failure of the boundary to be \emph{convex} with respect to the flow (see Definition 4.1). In fact, on manifolds with boundary, one can \emph{trade} the $f$-singularities in the interior of $X$ for these boundary  effects. In our approach, the boundary effects take the central stage, while the singularities themselves remain in the  background. In the paper,  we apply this philosophy to  $3$-manifolds. Many of our results allow for straightforward multidimensional generalizations, the other are specifically three-dimensional. 

Some of our theorems are in the spirit  of  the pioneering work of I. Ishii  on, so-called, \emph{flow-spines}  [I], [I1] (see also a recent paper by Y. Koda [Ko] and an excellent monograph "Branched Standard Spines of 3-manifolds" by Benedetti and Petronio [BP], followed by [BP1]). In an earlier version of this paper [K], we managed to overlook all this line of research...  As we will introduce the relevant constructions, we will describe some technical differences between the flow-spines of [I] and the branched spines of [GR1], [BP], on one hand, and the gradient spines on the other. For now, it is sufficient to say that any generic gradient flow defines its gradient spine  in a \emph{canonical} way, and that we allow for fields $v$  that are not necessarily concave with respect  to the boundary $\d X$.

Our motivation comes from the desire to understand better the  interplay between the intrinsic concavity of $\d X$ with respect to generic gradient flows and the topology of the underlining 3-fold $X$. We conjecture that there exists a numerical topological invariant that measures the failure of convexity with respect to any nonsingular gradient flow---``\emph{some manifolds intrinsically are just more concave than others...}" In a sense, the gradient complexity $gc(X)$, introduced in this paper, can serve as a crude measure of intrinsic concavity of $X$.
In fact,  a 3-manifold $X$ with a connected boundary which admits a convex gradient-like field $v$ is a handlebody;  so a random manifold does not admit convex nonsingular gradient flows. For instance, $H_2(X; \Z) \neq 0$ constitutes an obstruction to the convexity for any nonsingular gradient flow. At the same time, any manifold with boundary admits a \emph{strictly concave} traversing (but not necessarily a gradient!) flow  [BP].
\smallskip

The combinatorial complexity theory of Matveev ([M]) helps us to uncover the behavior of generic nonsingular gradient-type flows on 3-folds or, rather,  the interactions of such flows with the boundary. Before describing these results in full generality, let us give to the reader their taste. For example, we prove that on a manifold $X$, obtained from the Poincar\`{e} homology sphere by removing an open disk, any nonsingular gradient-like flow has at least \emph{five} trajectories that are tangent to the boundary $\d X$, each one at a \emph{pair} of distinct points; moreover $X$ admits a gradient-like flow with not more than $6 \cdot 5 = 30$ such trajectories.  Another example is provided by the remarkable  hyperbolic manifold $M_1$ that has the  minimal (among hyperbolic manifolds) volume $V \approx 0.94272$. By removing an open disk from $M_1$ we get a manifold $X$ on which any nonsingular gradient-like flow has at least \emph{nine} trajectories, each one tangent to the sphere $\d X$ at a pair of distinct points; moreover $X$ admits a gradient-like flow with not more than $6\cdot 9 = 54$ double-tangent trajectories.  
\smallskip

\smallskip

A generic vector field $v$  on $X$ gives rise to a natural stratification 
\begin{eqnarray}
X \supset \d_1^+X \supset \d_2^+X \supset \d_3^+X 
\end{eqnarray}
by compact submanifolds,  where $dim(\d_j^+X) = 3 - j$. 
Here  $\d_1^+X$ is the part of the boundary $\d_1X := \d X$ where $v$ points inside $X$. $\d_2X$ is    a 1-dimensional locus where $v$ is \emph{tangent} to the boundary $\d X$. Its  portion $\d_2^+X \subset \d X$ consists of points  where $v$ points inside $\d_1^+ X$. Similarly, $\d_3X$ is a finite locus where $v$ is tangent to $\d_2 X$. Finally, $\d_3^+X \subset \d_3 X$ consists of points where  $v$ points inside $\d_2^+X$. \smallskip

In his groundbreaking 1929 paper [Mo], Morse discovered some beautiful connections of this stratification  to the index of the field $v$.\footnote{Actually, the results of [Mo] apply to compact manifolds  $X$ of any dimension.} \smallskip

Now, let us describe the content of our paper section by section.
\smallskip

\smallskip

{\bf Section 2} starts with a sketch  of  main results from [Mo]  (see Theorem 2.1 and Corollary 2.1). It also contains one remark  about the role that stratification (1.1) plays in the Gauss-Bonnet Theorem (see Theorem 2.2 and [G] for an interesting general discussion). More importantly, we notice that, at any point $x \in \d X$, the $v$-flow defines a projection of the boundary $\d X$ into a germ of the constant level surface $f^{-1}(f(x))$. At a generic point   $x \in \d_2 X$ this projection is a \emph{fold}, while at $x \in \d_3X$ it is a \emph{cusp}. Throughout the paper, these folds and cusps provide us with crucial measuring devices for probing the topology of $X$. A significant portion of the paper is preoccupied with role of the cusps. \smallskip

\smallskip

{\bf Section 3} As in [Mo], the stratification $\{\d_j^+X\}_j$ is in the focus of our investigation. Here we prove that the surface $\d_1^+X$ can be subjected to 1-sugery via a deformation of the gradient-like field $v$. This allows one to change the topology of the stratum $\d_1^+X$ almost at will (see Lemma 3.1 and Corollary 3.1).\smallskip

\smallskip

{\bf Section 4}  For given nonsingular Morse data $(f, v)$, we introduce the notion of $s$-\emph{convexity},  $s = 2$ or  $3$. The 2-convexity of $v$ is defined as the property $\d_2^+X = \emptyset$. It puts a severe restrictions on the topology of $X$ (see Theorem 4.2 and Corollary 4.5). In contract, the 3-convexity, $\d_3^+X = \emptyset$, by  itself has no topological significance: one can always deform $(f, v)$  to eliminate  $\d_3^+X$ together with all other cusps (Theorem 9.5).\footnote{Note that Theorem 4.1.9 in [BP]  implies $\d_3X = \emptyset$ for all, so called,  \emph{traversing} flows.} However, when we fix the topology of  $\d_1^+X$, some combinations of cusps from $\d_3X$ acquire topological  invariance (Corollary 9.2).

Although convexity or its lack are defined in terms of the gradient-like fields,  
we can arrange  for the 2-convexity if we know that singularities of $f|_{\d X}$ admit a 
particular \emph{ordering} induced  by $f$ (similar to the self-indexing property). Specifically, the singularities of $f|_{\d X}$ can be  divided into two groups: the positive $\Sigma_1^+$ where the gradient $v = \nabla f$ is directed inwards $X$, and the negative $\Sigma_1^-$ where $v$ is directed outwards (see Fig. 1). Theorem 4.1 claims that when $f(\Sigma_1^-)$ is \emph{above}  $f(\Sigma_1^+)$, then one can deform the riemannian metric on $X$ so that the convexity of the gradient flow will be guaranteed. Hence, it is impossible to find a nonsingular function $f$ with the property $f(\Sigma_1^-) > f(\Sigma_1^+)$ on 3-folds $X$ that are not handlebodies. In addition, Theorem 4.2 describes an interplay between the dynamics of the flow $v$ through the ``bulk" $X$ and of  the $v$-induced flow $v_1$ in $\d_1X$, on the one hand, and the convexity phenomenon, on the other. 

In Corollary 4.5, we prove that an acyclic $X$ is a 3-disk if and only if one of the two properties are satisfied: (1) $X$ admits nonsingular 2-convex Morse data $(f, v)$, (2) $X$ admits nonsingular 3-convex Morse data $(f, v)$ with a connected 
$\d_1^+X$.
\smallskip

\smallskip

{\bf Section 5} is devoted to properties of \emph{gradient spines}, a construction central to our investigations. In spirit, but not technically, it represents a special class of flow-spines [I]. The difference between the two classes reflects the difference between the spaces of nonsingular vector  and  gradient fields on a given manifold. 

Recall that a \emph{spine} $K \subset X$ is a compact cellular two-dimensional subcomplex $K$ of the 3-fold $X$, such that  $X \setminus K$ is homeomorphic to the product $$[\d X \setminus (\d X \cap K)]\times [0, 1).$$ 

The relation between general spines and ambient 3-folds is subtle: a manifold $X$ has many non-homeomorphic spines $K$, and there are topologically distinct $X$ that share the same spine. In order to make the reconstruction of $X$ from $K$ possible, $K$ has to be rather special (cf. [M], [BP]).

In fact, a generic nonsingular  gradient-like field $v$ canonically gives rise to a spine that we call \emph{gradient} (see Fig. 10, 12).  A gradient spine $K$ is a union of $\d_1^+X$ with the descending  $v$-trajectories that pass through 
$\d_2^+X$. Like branched spines (Definition 6.6), gradient spines inherit \emph{orientations} from the boundary $\d X$ and have a \emph{preferred side} in the ambient $X$. Exactly these properties of  a gradient spine $K$ allow for its resolution into a surface $S$ homeomorphic to $\d_1^+X$ and, eventually, for a reconstruction of $X$ from $K$ (Theorem 6.1). By modifying the field $v$, we can arrange for $\d_1^+X$, and thus for $S$, to be homeomorphic to a  disk $D^2$. As a result, we get our Origami Theorem 5.2: any  3-manifold $X$ with a connected boundary has a gradient spine $K$ obtained from a disk $D^2$ by identifying certain arcs in $\d D^2$ with the appropriate arcs in the interior of $D^2$ (see Fig. 14, 15). This statement is similar in flavor to Theorem 1.2 , [I], where  disk-shaped sections of  generic (non-gradient) flows on closed manifolds are employed for the same goal.\smallskip

\smallskip

{\bf Section 6} deals with combinatorial structures that generalize the notion of a gradient spine $K$ (see Definitions 6.1---6.4). We start with a 2-complex $K$ whose local geometry is modeled  after gradient spines (see Fig. 20).  Adding a system of, so called, $TN$-\emph{markers} to $K$ along its singularity set $s(K)$ produces an object which captures the  topology of the ambient $X$ and admits a canonic resolution into an oriented surface. Such a polyhedron $K$ with markers is called an \emph{abstract gradient spine}. Unlike generic 2-complexes, each abstract gradient spine $K$ is a spine of some manifold (see [BF] and [K], Appendix, Theorem 10.1).\footnote{In a sense, the category of abstract gradient spines is equivalent to the category of compact 3-manifolds with  non-empty boundary.}

The notion of a $\vec Y$-\emph{spine} (see Definition 6.5) is still another generalization of gradient spines. It is a very close relative of branched spines. In fact, for an oriented  $X$, the notions of  an oriented branched spine $K \subset X$ and of a $\vec Y$-spine are equivalent, provided $K^\circ$ being orientable (Lemma 6.3). Unlike  abstract gradient or branched spines, the $\vec Y$-spines $K$ are defined \emph{extrinsically}, that is, in terms of an  embedding in $X$ of the vicinity of the singular set $s(K) \subset K$.  By Lemma 6.2, any gradient spine is a $\vec Y$-spine. Moreover, according to Theorem 8.1,  $\vec Y$-spines admit a ``nice" approximation by the gradient spines of the same complexity. 
\smallskip

\smallskip

{\bf Section 7}\; We apply ideas and results of [M], which revolve around  Matveev's  notion of combinatorial complexity of simple 2-complexes and compact 3-folds, to the gradient and $\vec Y$-spines. We introduce the \emph{gradient complexity} $gc(X)$ of a 3-fold $X$ with boundary as the minimal number of \emph{double-tangent}  trajectories that  a nonsingular gradient-like field on $X$ can have. A double-tangent trajectory is tangent to the boundary $\d_1X$ at a pair of distinct points. In general, $gc(X) \geq c(X)$, where $c(X)$, the Matveev combinatorial complexity, is defined to be the minimal number of special isolated singularities\footnote{called, butterflies in [M] and $Q$-singularities in this paper} that a simple spine $K \subset X$ can have.  One can restrict the scope of this definition only to $\vec Y$-spines (equivalently, to oriented branched spines) in order to get the notion of  $\vec Y$-complexity $c_{\vec Y}(X)$. We prove that $gc(X) \geq c_{\vec Y}(X) \geq c(X)$. In fact, Theorem 8.2 claims that  $gc(X) = c_{\vec Y}(X)$, and, for the geometrical pieces $X$ of the SJS decomposition, $gc(X) \leq 6 \cdot c(X)$.

The inequality $gc(X)  \geq c(X)$ helps us to restate many results from $[M]$ in the language of  double-tangent trajectories. For instance, by Theorem 7.3, for any natural $c$, there is no more than finitely many irreducible and boundary irreducible with no essential annuli 3-folds $X$ that admit nonsingular gradient-like flows with $c$  double-tangent trajectories. The number $N(c)$ of such 3-folds has a crude upper bound $\Gamma_4(c) \cdot 12^c$, where $\Gamma_4(c)$ stands for the number of topological types of regular four-valent graphs with $c$ vertices at most. In particular, there is no more than $\Gamma_4(c) \cdot 12^c$ hyperbolic manifolds with $c$ double-tangent trajectories. 

Let $X$ being obtained from a closed hyperbolic 3-fold $Y$ by removing a number of open balls. By Theorem 7.5, any non-singular gradient-like flow (as  well as any convex traversing flow) on $X$ has at least $V(Y)/V_0$ double-tangent trajectories. Here $V(Y)$ stands for the hyperbolic volume of $Y$ and $V_0$ for the volume of the perfect ideal tetrahedron. 

Fortunately, all orientable irreducible and closed 3-manifolds of complexity at most six (there are 74 members in this family) have been classified and their minimal spines have been listed [M]\footnote{By definition, a spine of a closed manifold $Y$ is a spine of the punctured $Y$, that is, of $Y \setminus D^3$.}. Some partial results are available for the 1155 closed irreducible manifolds of complexity at most nine. This has been accomplished by an algorithmic computation coupled with  ``hands on" analysis of spines that look different, but share the same values of the Turaev-Viro invariants [TV]. The bottom line is that all $X$ with $c(X) \leq 6$ are distinguished  by their Turaev-Viro invariants! Thus, for each manifold $Y$ in the Matveev list and any generic nonsingular gradient flow on $X = Y \setminus D^3$, we get a lower bound on the number of double-tangent (to the boundary $\d_1X  \approx S^2$) trajectories. 

Consider any irreducible and orientable  3-manifold $X$ produced from a closed manifold $Y$ by removing a ball. In Corollary 7.1, we prove that if $X$  admits a nonsingular gradient-like flow with no double-tangent trajectories, then $X$ is the standard disk. In view of Perelman's work [P1], [P2], this is not  an exciting fact, but it shows how far we can get with our Morse-theoretic techniques. 

Section 7 contains  a few more results about upper and lower estimates of $gc(X)$ for manifolds obtained from closed manifolds $Y$ by removing a number of 3-balls. Theorem 7.6 provides a lower bound for $gc(X)$ in terms of the presentational complexity of the fundamental group $\pi_1(X)$. At the same time, any self-indexing Morse function $h$ on $Y$ gives rise to an upper estimate of $gc(X)$ (given in terms of the attaching maps for the unstable 2-disks of index two $h$-critical points).

In Theorem 7.4, we notice that $gc(X)$ can increase only as a result of 2-surgery on $X$. 
\smallskip

\smallskip

{\bf Section 8} \; Here we are addressing a natural question: Which spines are of the gradient type? The main result of the section, Theorem 8.1, claims that any $\vec Y$-spine $K \subset X$ can be approximated by  a gradient spine $\tilde K$; moreover,  $c(\tilde K) \leq c(K)$. Furthermore, one can get $K$ from $\tilde K$ by controlled elementary collapses of certain 2-cells. Theorem 8.1 depends on some results from Section 9 about possible cancellations  of cusps from $\d_3X$. Theorem  8.2 establishes the equality $c_{\vec Y}(X) = gc(X)$ and the crucial inequality $c(X) \leq gc(X) \leq 6\cdot c(X)$.

\smallskip

\smallskip

{\bf Section 9} (together with Section 7) contains our main results. Here we analyze the effect of deforming Morse data $(f, v)$ on the gradient spine they generates. Theorem 9.3 claims that, in the process, the gradient spine goes through a number of elementary expansions and collapses of two-cells mingled  with so called $\a$- and $\b$-\emph{moves} (see Fig. 29, 30). These are analogs of the second and third Reidemeister moves for link diagrams. Theorem 9.1 describes possible cancelations of cusps from $\d_3X$ that accompany generic deformations of $v$.  One of our main results, Theorem 9.4, is a combination of Theorem 9.3 with a special case of Phillips' Theorem [Ph]. We prove that when two nonsingular functions $f_0$ and $f_1$ on $X$ produce the same invariants $h(f_0), h(f_1) \in H^2(X; \Z)$---the same $Spin^c$-structures in the sense of Turaev [T]---, then their gradient spines are linked by a sequence of elementary $2$-expansions, $2$-collapses, and $\a$- and $\b$-moves (see the proof of Corollary 8.1 for the definition of the invariant $h(f) \in H^2(X; \Z)$).

Deformations of $(f, v)$ that cause jumps in the value of $h(f)$ (in the  $v$-induced $Spin^c$-structure) manifest themselves as a ``disk-supported  surgery on the preferred spine orientation".  We call them \emph{mushroom flips} (see Fig. 35). 

In Theorem 9.5, we prove that, given generic Morse data $(f, v)$,  it is possible to deform them so that  all the cusps from $\d_3X$ will be eliminated, but the number of double-tangent trajectories $gc(f, v)$ will be preserved.
\smallskip

\smallskip


Finally, it should be said that the Morse theory on stratified spaces, in general,  and  on manifolds with boundary, in particular, has been an area of an active advanced and interesting research. For a variety of perspectives on this topic see  [Mo], [GM],  [F], [C], [Ha]. Our intension is to bring the stratified Morse theory and the complexity theory of 3-folds under a single roof.  

\bigskip

{\it Acknowledgments} \quad This paper is  shaped by numerous and valuable  discussions I had with Kiyoshi Igusa.  My deep gratitude goes to him. I am grateful to Yakov Eliashberg for pointing  that some propositions below (Theorem 9.6,   Lemma 3.1 and Corollary 3.1) are intimately related to and similar in spirit  with his general theory of folding maps [E1], [E2]\footnote{They also bare resemblance to some results of Harold Levine [L].}. I am also very grateful to the referee of [K] who informed me about  the existing results of [I], [I1], [BP], [BP1]  and others.


\section{The Morse Stratification on Manifolds with Boundary}

Let $X$ be a compact $3$-manifold with boundary $\d X$. Let
$f: X \rightarrow \R$ be a generic smooth function. Then $f$ has  
non-degenerate critical points in the interior of $X$ and the  
restriction of $f$ to the boundary $\d X$ is also a Morse function. Let  
$v$ be a gradient-like vector field for $f$, that is, $df(v) > 0$ away from the $f$-critical points. 
Instead of working with  such pairs $(f, v)$, we can pick a  Riemannian metric on $X$ and choose $v = \nabla f$, the gradient field. Both points of view are equivalent, but we prefer  the first. \smallskip

The singularities of $f|_{\d X}$ come in two flavors: \emph{positive} and \emph{negative}.  
At a positive singularity, the field $v$ is directed inward $X$, and at a negative singularity, ---
outward. This distinction between positive and negative critical points  of $f|_{\d X}$ depends on
$f$, not on $v$. At a positive singularity and
in an appropriate coordinate system $\{x_1, x_2, x_3\}$ with
$\{x_1 = 0\}$ defining  $\d X$ and $x_1 > 0$ --- the interior of $X$,
$$
f(x) = c + x_1 +  a_2 x_2^2 + a_3 x^2_3,
$$
where $c$ and $a_i \neq 0$ being constants.  At a negative singularity, one  
has
$$
f(x) = c - x_1 +  a_2 x_2^2 + a_3 x^2_3.
$$

Let $\Sigma^{\pm}_1$ be the set of positive/negative singularities of  
$f|_{\d X}$ and  let $\Sigma_0$---the set of singularities of $f$ in the interior of $X$.
Denote by $X_{\leq c}$ the set $\{x\in X |\; f(x) \leq c\}$.\smallskip

Crossing the critical value $c_\star$ of a positive singularity causes  
the topological type
of $X_{\leq c}$ to change, while crossing $c_\star$ of a negative  
singularity has no effect
on the topology of  $X_{\leq c}$ as illustrated in Figure 1.


\begin{figure}[ht]
\centerline{\includegraphics[height=1.5in,width=3in]{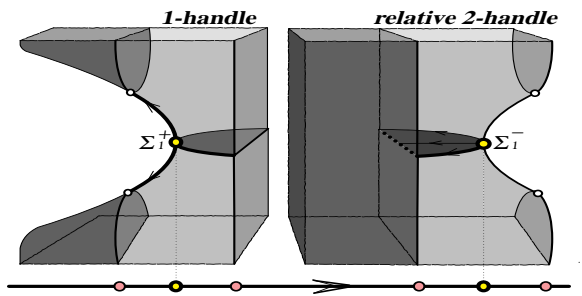}}
\bigskip
\caption{\small{A positive singularity of index 1 and a negative singularity  
of index 2 on the boundary of a solid.
The gradient-like field $v$ is horizontal.}}
\end{figure}

For a generic field $v$,  the locus $L$ where the field is tangent to $\d X$ is a \hfill\break 1-dimensional submanifold of the boundary; $L$ divides $\d X$ into two domains:  $\d^+X$ where $v$ is directed inwards of $X$, and  in  $\d^-X$, where it  is directed outwards.

Morse noticed that, for a generic vector
field $v$, the tangent locus $L$ inherits a structure in relation to $\d^+X$  
analogous to that of  $\d X$ in relation to $X$ \cite{Mo}. To explain
this point we need to revise our notations in a way which will be  
amenable to recursive definitions.\smallskip

Let $\d_0X := X$,  and $\d_1X := \d X$.   Denote by $\d_2X \subset \d_1X$ the locus where $v$ is tangent to $\d_1X$. For a generic $v$,  $\d_2X$ divides $\d_1X$ into a domain $\d_1^+X$ where $v$ is directed inwards $X$ and a domain $\d_1^-X$ where $v$ is outwards inwards $X$. Evidently, $\d_1^\pm X \supset \Sigma_1^\pm$. Consider the set $\d_3X$ where $v$ is tangent to $\d_2X$. The set $\d_3X$ divides $\d_2X$ into a set $\d_2^+X$ where $v$ is directed inwards $\d_1^+X$ and a set $\d_2^-X$ where $v$ is directed outwards $\d_1^+X$. Finally, $\d_3X = \d_3^+X \coprod \d_3^-X$, where $v$ is directed inwards 
$\d_2^+X$ at the points of  $\d_3^+X$. 
\smallskip

From now and on, we call $(f, v)$ \emph{generic} if 1) all the strata ${\d_j X}_{1 \leq l \leq 3}$ are regularly embedded smooth manifolds and 2) all the restrictions $f|_{\d_j X}$ are Morse functions. Most of the time, the second property will be  irrelevant, but when we need it, we do not want to modify our definition.  At some point, the word ``generic" will mean an additional  general position requirement imposed on the field $v$ (see Definition 5.2). When we say that a Riemannian metric is generic, we imply that $(f, \nabla f)$ is generic.\smallskip

We introduce critical sets $\Sigma_j^\pm \subset \d_j^\pm X$ of $f|_{\d_jX}$ in a way similar to the one we used to define  $\Sigma_1^\pm$. With some generic metric in place, let $v_j$ be the orthogonal projection of $v$ onto $\d_jX$, and let $n_j$ denote the  normal field to $\d_jX$ inside $\d_{j-1}X$ that points inside $\d_{j-1}^+X$.  
Note that, away from the singularities from $\Sigma_j$, $df (v_j) > 0$. 
\bigskip

\begin{figure}[ht]
\centerline{\includegraphics[height=1.8in,width=4in]{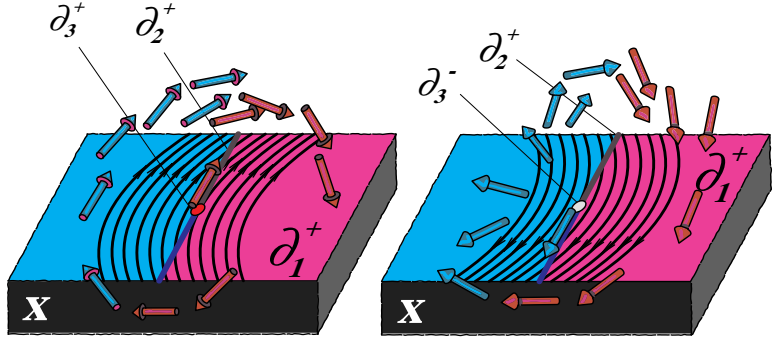}}
\bigskip
\caption{\small{The patterns of  fields $v$ (the $3D$-arrows) and $v_1$ (the  
parabolic flow) in vicinity of a point from $\d_3^+X$ (on the left) and a 
point from $\d_3^-X$ (on the right).}}
\end{figure}

For a vector field $v_k$ as above on $X_k$ with isolated singularities $\{x_\star  
\in \Sigma_k \subset Int(X_k)\}$, denote
by $Ind_{x_\star}(v_k)$ its index at $x_\star$, and by $Ind^+(v_k)$ --- the  
sum $\sum_{x_\star \in \Sigma^+_k} \; Ind_{x_\star}(v_k)$.
Then, according to [Mo], one has two sets of equivalent relations:

\begin{thm} {\bf (Morse Law of Vector Fields)}.  For any   
generic metric and $0 \leq k \leq 3$,
\begin{itemize}
\item $\chi(\d_k^+X) = Ind^+(v_k) + Ind^+(v_{k+1})$ \footnote{By  
definition,
$Ind^+(v_3) = \#(\Sigma_3^+)$, and $Ind^+(v_4) = 0$.}

\item $Ind^+(v_k) = \sum_{j = k}^3 \; (-1)^j \chi(\d_j^+X).$
\end{itemize}
\end{thm}

\begin{figure}[ht]
\centerline{\includegraphics[height=1.3in,width=1.8in]{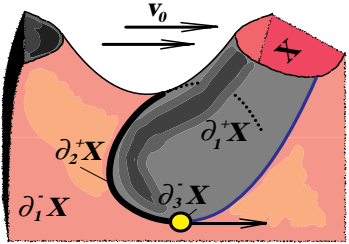}}
\bigskip
\caption{\small{A more realistic picture of the boundary $\d_1X$ in vicinity  
of $\d_3^-X$ in relation to the
horizontal gradient field $v$.}}
\end{figure}

\begin{cor}
For generic vector field $v$ and  a metric on $X$,
$$Ind(v) = \sum_{k = 0}^3 (-1)^k \chi(\d_k^+X). \qed$$
\end{cor}

For an engaging discussion of the Morse Theorem 2.1 see  the paper of Gottlieb
[G].\footnote{That nice paper attracted my attention to the topic 
of Morse theory on manifolds with boundary.} In particular, it describes a link between the Morse 
stratification $\{\d_j^+X\}_j$ and the geometry (normal curvature $K$) of  $\d_1X$:

\begin{thm}
Let $\Phi: X \to \R^3$ be a smooth map with a nonzero Jacobian on the boundary $\d X$ and 
$p: \R^3 \to \R$ a generic linear function, so that the function $f := p\circ\Phi$ has only isolated 
singularities in $Int(X)$.
Then the degree of the Gauss map $g: \d X \to S^2$ can be calculated either by integrating the normal 
curvature $K$ of $\Phi(\d X) \subset \R^3$ (Gauss-Bonnet Theorem), or in terms of the  $v$-induced stratification $\d_3^+X \subset \d_2^+X \subset \d_1^+X \subset X$:
\begin{eqnarray} 
deg(g) = \frac{1}{4\pi} \int_{\d X} K d\mu = \chi(X) - Ind(v) \\ \nonumber
= \chi(\d_1^+X) - \chi(\d_2^+X) + \chi(\d_3^+X)
\end{eqnarray}
\end{thm}

We notice that formulas from Theorem 2.1 and Corollary 2.1 admit \emph{equivariant generalizations}.  For any compact Lie group $G$, a $G$-manifold $X$, an equivariant function $f: X \to \R$,  and a generic eqivariant metric on $X$ 
(generating $G$-equivariant gradient-like fields $(v, v_1, v_2)$),  the invariants $\{\chi(\d_k^+ X)\}$, as well as the degree $deg(g)$, can be interpreted as taking values in the Burnside ring $\mathcal B(G)$ of $G$ (see [TD] for the definitions). \smallskip

There is another  degree-type  invariant of $(X, f)$ linked to generic Morse data $(f, v)$.  The set $\d_2X \subset \d_1^+X$ carries two non-zero vector fields: the normal field $n_2$ that points inside $\d_1^+X$ and trivializes  the oriented tangent bundle of $\d_1X$ along $\d_2X$, and the field $v = v_1$. Therefore, $v$ defines a map  $h: \d_2X \to S^1$. We view $h$ as an element in the one-dimensional oriented bordism group $\Omega_1(S^1)$ of the circle. This group splits as $\Omega_1(pt) \oplus \Omega_0(pt)  \approx \Omega_0(pt) $ (see [CF]), i.e., an element  $h: M^1 \to S^1$   in $\Omega_1(S^1)$ is determined by the degree class $deg(h) = [h^{-1}(pt)] \in \Z$. 
 
Any deformation of  $v$ preserves the class of $h: \d_2X \to S^1$ in $\Omega_1(S^1)$ and thus the degree $[h^{-1}(pt)]$. Deformations of $f$ that change the singularity set $\Sigma_1$ do change the degree class. This degree can be easily computed in terms of the cusp sets $\d_3^+X$ and $\d_3^-X$.
 
 \begin{lem} For a fixed $f$, the number $\#(\d_3^+X) - \#(\d_3^-X)$ equals to twice the degree of the map 
 $h: \d_2X \to S^1$ and is independent of the fields $v, n_2$.
 \end{lem}
 
{\it Proof}\; Each loop $\gamma$ from $\d_2X$ either entirely belongs to one of the two sets $\d_2^+X$ and to  $\d_2^-X$, or the arcs of $\gamma$ belonging to $\d_2^+X$ and to  $\d_2^-X$ alternate. In the first case, the 
contribution of $\gamma$ to $deg(h)$ is zero. In the second case, the contribution of each arc with the ends of opposite polarity is also zero. Each arc with two positive ends contributes a rotation of $v$ by $+\pi$, while 
each arc with two negative ends contributes a rotation by $-\pi$  (see Fig. 7). Hence the total rotation along $\gamma$ is $\pi[\#(\d_3^+X) - \#(\d_3^-X)]$. \qed
\smallskip

By Corollary 9.2,  a more refined count of the cusps from $\d_3X$ will produce a very different formula for the degree of $h: \d_2X \to S^1$.
\bigskip

For a given nonsingular $f: X \to \R$, each choice of a gradient-like field $v$ \emph{locally} gives rise to a map 
$p: X \to \R^2$.  Let us outline the construction of $p$. Add an external collar $W$ to $X$ and extend the Morse data $(f, v)$ into $Y : = X \cup W$ without adding new singularities. At each point $x \in Int(W)$ the $(-v)$-flow defines a surjection $p_x$ of a neighborhood $U_x \subset Y$ onto a  neighborhood $V_x$ of $x$ in $f^{-1}(f(x))$. Consider the  restriction  
$p_x : U_x \cap \d_1 X \to V_x$, $x\in \d_1 X$, to the boundary $\d_1 X$. According to Whitney [W], generic smooth maps $\R^2 \to \R^2$ have only folds and cusps as their stable singularities. Therefore, for  generic Morse data $(f, v)$ and $x \in \d_1X \setminus \d_2X$, $p_x: U_x \cap \d_1 X \to V_x$ is a surjection,  for $x \in \d_2X \setminus \d_3X$, $p_x$ is a folding along an arc of $\d_2X$, and at $x \in \d_3X$, $p_x$ is a cusp map with $p_x(\d_2 X)$ being the cuspidal  curve. Note that along $\d_2^+X$, $p_x: X \to V_x$ is locally onto, while along $\d_2^-X$, it is not. 

It is especially easy to visualize the stratification $\{\d_jX\}$ when $X$ is embedded or immersed in $\R^3$ and $f$ is induced from a generic  linear function $l$ on $\R^3$. In such a case, a global surjection $p: X \to \R^2$ is available.  Its fibers are parallel to the gradient vector $v = \nabla l$. Now, $\d_2X$ can be identified with the folds of the map $p: \d_1X \to \R^2$ and  $\d_3X$ with its cusps. 
\smallskip

As we deform a nonsingular field $v$ within  generic one-parameter families, the local structure of the projections $p_x$ can be described in terms of a few canonical forms.
One of them, the \emph{cusp}, 
\begin{eqnarray}
F(x, y)= (x^3+xy,\, y)
\end{eqnarray}
is a stable singularity of a map from $\R^2$ to itself.\footnote{It comes from the universal unfolding of the $A_3$ singurality $f(x)=x^3$. }

The \emph{dove tail} $t$-parameter family
\begin{eqnarray}
F_t(x, y)=(x^4+x^2t+xy,\, y)
\end{eqnarray}
describes a cancellation of two cusps that will play a significant role in Section 9.\footnote{It
is the universal unfolding of the codimension 1 singularity of a mapping $\R^2$ to $\R^2$ coming from the universal (two-parameter) unfolding of $A_4$  singularity.}


\section{Surgery on the Morse Stratification}

Let $X$ be a  compact 3-manifold $X$ with boundary $\d_1X$.
Given a smooth function $f: X \rightarrow \mathbb R$ with isolated singularities, we can construct a  
new function with \emph{no singularities inside} $X$: just cut from $X$ 
a number of tunnels. Each tunnel starts at the 
boundary  $\d_1X$ and has a dead end which engulfs a singularity. Denote by $T$ the interior of the tunnels.  Then $f$, being restricted to $X \setminus T \approx X$, is nonsingular, and its perturbation can be assumed to be of  the Morse type on $\d(X \setminus T)$.
\bigskip

\begin{lem} Let $X$ be a compact 3-manifold with boundary $\d_1X$.
Let $f: X \rightarrow \R$ be a smooth function with no  
singularities in a regular neighborhood $N$ of $\d_1X$. Denote by $v$ be its gradient-like field.
Let $\g \subset \d_1^\pm X$ be a simple path that connects two points from $\d_2X$ and 
has an empty intersection with the critical set $\Sigma_1^\pm$.

Then one can deform $v$ in $N$ to a new $f$-gradient-like vector field  
$\tilde v$ for which the  new set $\d_1^\mp X$ will be obtained from 
the original one by the one-surgery along $\g$. 
Outside of $N$,  $v = \tilde v$.
\smallskip

A similar statement holds for any field $v$\footnote{not necessarily of the gradient type} which is nonsingular along $\d_1X$ and in general position to it.
\end{lem}
\begin{figure}[ht]
\centerline{\includegraphics[height=1.2in,width=2.4in]{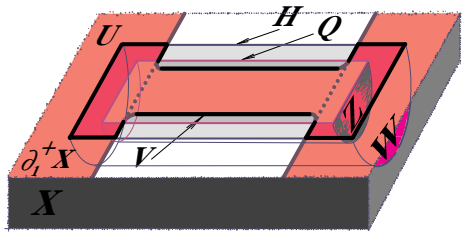}}
\bigskip
\caption{\small{Performing $1$-surgery on $\d_1^+X$.}}
\end{figure}
{\it Proof.}\; Let $v_1$ be the orthogonal projection of 
$v$ onto $\d_1X$ in the metric on $X$ in which $v = \nabla f$. The idea is to perform surgery on $\d_1^+X$ 
by a homotopy of the field $v - v_1$, while keeping  $f$ and $v_1$ fixed. 
We start with  ``1-surgery" on the fields along a  
band $H$ which connects two arcs, say $A \subset \d_2X$ and $B \subset \d_2X$.
The band, with the exception of small neighborhoods of its two ends,  
resides in $\d_1^+X$ (alternatively, in $\d_1^-X$)
as shown in Fig. 4. The band avoids the singularities
of the function $f|_{\d_1X}$, so that $v_1 \neq 0$ everywere in the  
band. Denote by $Q$ a smller band which is contained in $H$ (see
Fig. 4).

Let $n$ denote the interior normal to $\d_1X$. We decompose the field  
$v$ as $v_1 + h\cdot n$,
where the function $h$ is positive in the open domain $U$ --- the shaded  
area without the handle (it is bounded on the left and right by the two dotted segments)--- and  is
negative in the interior of the complement to $U$.
In fact, we can assume that $0$ is a regular value of $h$.

At each point $x \in X$, the differential $df$ picks a particular open  
half-space $T^+_{f,x}$ in the
tangent space $T_x$, and $v \in T^+_f$.  Along the boundary $\d_1X$,  
another family of half-spaces is available: let $T^+_x$
denote the set of tangent vectors at $x \in \d X$ which point inside of  
$X$. Note that, away from the singularities of $f|_{\d_1 X}$, the cone 
$T^+_{f, x} \cap T^+_x$  is open.
\begin{figure}[ht]
\centerline{\includegraphics[height=1.1in,width=1.8in]{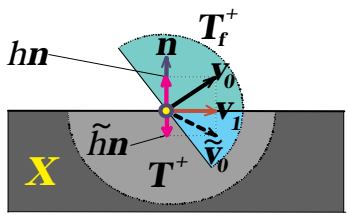}}
\bigskip
\caption{\small{Changing the field $v = v_1 + hn$ at a point  $x  
\in \d_1^-X$ into a field $\tilde v = v_1 + \tilde hn$ for which $x \in  
\d_1^+X$.}}
\end{figure}

Consider a smooth function $\tilde h : H \rightarrow \mathbb R$ which  
satisfies the following properties: 1) $\tilde h^{-1}([0, +\infty))  
= V$, 2) zero is a regular value of $\tilde h$ and $\tilde h^{-1}(0) = \d V$,  
3) the field $v_1 + \tilde h\cdot n \in T^+_{f}$. The last property can be achieved  
by starting with any $\tilde h$ subject to 1) and 2) and rescaling it by a variable 
factor $a > 0$, so that  
$v_1 + (a\tilde h)n \in T^+_f \cap T^+$ (see Fig. 5). At each point,  the existence 
of an appropriate $a$ follows from the fact that $v_1 \in T^+_f$. Because the 
cone $T^+_f \cap T^+$ is open and the domain $H$ is compact,
the global existence of such an $a$, by a partition-of-unity argument, follows from 
its existance at each point  of $H$.

We extend the field $v_1 + \tilde h\cdot n$ inside $X$ to get a smooth an  
$f$-gradient-like field $w$ in a small regular neighborhhod $W$ of $V$. 
We use a smooth partition of  unity $1 = \alpha + \beta$
subordinate to the cover $W$, $X \setminus Z$. The function $\alpha$  
vanishses in $Z$ and $\beta$ in
$X \setminus W$. Now consider the field $\tilde v :=\alpha v +  
\beta w$. Since $T^+_x$ is convex,
$\tilde v \in T^+_f$. Moreover, in $V$,  $\tilde v$ points inside $X$ and, in 
$H \setminus V$,  outside  $X$. Also, outside $W$,  $\tilde v = v$.
\smallskip

Finally, for a fixed $f$, the set of all $f$-gradient-like fields is open and convex.  
Hence, any modification of a $f$-gradient-like field can be obtained by its deformation.

The arguments for generic (non-gradient) fields $v$ are similar and simpler.
\qed
\smallskip

\begin{cor} 
Under hypotheses and notations of Lemma 3.1, the following claims are valid. There is a  
deformation of a given gradient-like field  in the neighborhood $N$ of $\d_1X$ so that, for 
the new gradient-like field, both portions $\d_1^\pm X_j$ of $\d_1^\pm X$ residing in each 
connected component  $\d_1X_j$ of $\d_1X$ are nonempty, and $\d_1^+X_j$ is 
homeomorphic to any given domain in $\d_1X_j$ with a nonempty complement.

In particular, for a given generic $f$ and all $j$'s, there exists a gradient-like field $v$ such that anyone of the two properties is satisfied:
\begin{itemize}
\item  $\d_1^+X_j$ is homeomorphic to a disk.
\item  $\d_1^+X_j$ and $\d_1^-X_j$ are homeomorphic surfaces.
\end{itemize}
A similar statement is valid in a category of generic nonsingular vector fields.
\end{cor}

{\it Proof.}\;
Let  $\d_1X_j$ be a component of $\d_1X$. 
When $\d_1X_j = \d_1^+X_j$ (or $\d_1X_j = \d_1^-X_j$), we can 
pick a point  $x \in \d_1 X$ where $v_1  \neq 0$. Then, employing an argument 
depicted in Fig. 5, we can deform the  
field $v$ in the vicinity of $x$ so that, with respect to the  
modified gradient-like field,   $x \in \d_1^-X$ ($x \in \d_1^+X$, correspondingly).  Thus we can assume that  
$\d_1^-X_j, \d_1^+X_j \neq \emptyset$. 

Now, by one-surgery on both $\d_1^+X$ and $\d_1^-X$, we can 
change  the topology of $\d_1^+X$ any way we like, as long as we keep 
keep the sets of both polarities nonempty. No matter how we change the two 
sets, we must keep $\Sigma^+_{1, j}$ inside $\d_1^+X$ and $\Sigma^-_{1, j}$ inside $\d_1^-X$.
In particular, we can deform the field so that $\d_1^+X_j$ is a 2-disk or, say, to 
insure that  $\d_1^+X_j$ is homeomorphic to  $\d_1^-X_j$.

Note that surgery on $\d_1^+X$  typically will change  the sets $\d_3^{\pm}X$.
\qed

   
\section{Morse Strata and Convexity}

\begin{defn} Given generic Morse data $(f, v)$ on a  manifold $X$ with  boundary  
we say that $v$  
is $s$-\emph{convex} (\emph{concave}), if $\d_s^+X = \emptyset$ ($\d_s^-X = \emptyset$, correspondingly). In particular, if $\d_2^+X  = \emptyset$ ($\d_2^-X  = \emptyset$), we say that $v$  is 
symply \emph{convex} (\emph{concave}).
\end{defn}

An existence of convex Morse data has strong topological implications. 
Let $\Sigma$ be a surface with boundary. We denote by $\mathcal  
L(\Sigma)$ a smooth
$3$-manifold with boundary obtained from the product $\Sigma\times[-1,  
1]$ by rounding
its corners $\d\Sigma\times\{\pm 1\}$ and by replacing a narrow  
cylindrical band
$\Sigma\times[-\epsilon, \epsilon]$ with a "curved parabolic" one as  
shown in Fig. 6.
The projection $\mathcal L(\Sigma) \rightarrow [-1, 1]$ defines a  
nonsingular function $f$.
The  vertical field $v$ in $\Sigma\times[-1, 1]$ is of the  
$f$-gradient type.
With respect to it, $\d_2^-\mathcal L(\Sigma) = \d\Sigma$.
We call the triple $(\mathcal L(\Sigma), f, v)$  a \emph{lense} based  
on $\Sigma$.
\begin{figure}[ht]
\centerline{\includegraphics[height=1.1in,width=2.2in]{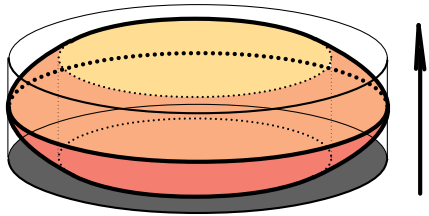}}
\bigskip
\caption{\small{A lense $\mathcal L(\Sigma)$ is convex with respect to the  
vertical field:
$\d_2^+\mathcal L(\Sigma) = \emptyset$. The set $\d_2^-\mathcal  
L(\Sigma)$ is  the equator of the lense.}}
\end{figure}

\begin{lem} 
A connected $3$-manifold $X$ 
admits convex nonsingular  data $(f, v)$  if and only if $X$ is diffeomorphic to a handle body   
$\mathcal L(\d_1^+X)$.

In particular, if an acyclic $3$-manifold $X$   
admits convex nonsingular data, it is a  $3$-disk.
\end{lem}

{\it Proof.}\; 
The first claim  is straightforward (see [BP], Proposition 4.2.2). Note that the convexity on a connected $X$ implies that $\d_1X$ is connected.

When $X$ is acyclic, a homological argument  
implies that $\d_1X \approx S^2$. Thus $\partial_1^+X$ must be   
a contractible domain in $S^2$, that is, a $2$-disk.
Therefore, the manifold $X$ must be shaped as a lens, one face of which  
is that disk. 
\qed
\bigskip

Fig. 7 shows a typical behavior of a vector field $v_1$ in a neighborhood  
of $\d_2X$.  The  arcs of $\d_2^+X$  come in tree flavors: $A$  is bounded by a pair of points from $\d_3^+X$, $B$ 
 is bounded by a pair of points from $\d_3^-X$ and $C$ is bounded by a pair of mixed polarity. \smallskip

This time, we play our convexity game in the dimension $2$, 
not $3$. At points of  $\d_3^-X$ the field $v_1$ in $\d_1^+X$
is \emph{convex}, at points of $\d_3^+X$ it is \emph{concave}.
\smallskip

\begin{figure}[ht]
\centerline{\includegraphics[height=1in,width=4in]{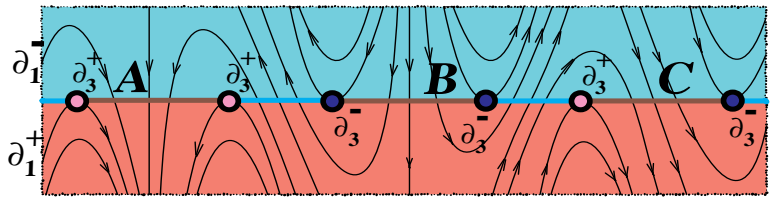}}
\bigskip
\caption{\small{The arcs from $\d_2^+X$ of the types A, B, and C.}}
\end{figure}

According to Theorem 9.5, we can deform $v$ (in the space of nonsingular gradient-like fields) so that $\d_3 X = \emptyset$, in other words, \emph{there are no topological obstructions to the 3-convexity and 3-concavity of (gradient) fields}!  On the way to establishing this fact, we need  to  perform $0$-surgery on $\d_3X \subset \d_2^+X$. 
\begin{lem} Let $v$ be a gradient-like field for $f: X \to \R$, $v \neq 0$ along $\d_1X$.
Let $C \subset \d_2^+ X$ be an arc with one of its ends  
$a\in \d_3^+X$,  the other end $b\in \d_3^-X$,
and  no other points of $\d_3X$ in its interior. Assume also that $f|_C$ has no critical 
points.\footnote{See Fig. 7, arc $C$. Note that the absence of critical points of $f|_C$ implies that 
$f(a) < f(b).$}Then we can deform $v$ in the vicinity of $C  
\subset X$ in such a way that: 
\begin{enumerate}
\item with respect
to the new $f$-gradient-like field $\tilde v$ the strata  
$\d_1^+X$ and $\d_2X$ remain the same, 
\item arc $C$ changes its polarity (from being in $\d_2^+X$ to  
being in $\d_2^-X$), and the points $a$,  $b$ are eliminated from the set $\d_3X$.
\end{enumerate}
\end{lem}

{\it Proof.}\; Let $v_2$ be an orthogonal projection of $v$ on $\d_2X$. The argument is analogous to the one in Lemma 3.1. However, this time, we will keep both the direction of $n_1$-component $v - v_1$ of the field $v$ and the field $v_2 \neq 0$ fixed in the vicinity of $C$, while deforming the field $v_1$. Because the  direction of the normal component $v - v_1$ remains the unchaged, the stratum $\d_1^+X$ and its boundary $\d_2X$  will be preserved, but the arc $C$ will change its polarity. \smallskip


\qed

\begin{cor}
For any  generic function $f: X \to \R$ with no critical points in $\d_1 X$, there exists a gradient-like field 
$v$ with the following property. Each arc $C \subset \d_2X$, which connects a minimum $x$ of $f|_{\d_2X}$ with a consecutive maximum $y$,  has:  
\begin{enumerate}
\item  a single point from $\d_3^-X$, provided $x \in \Sigma_2^+$ and $y \in \Sigma_2^-$, 
\item a single point from $\d_3^+X$, provided $x \in \Sigma_2^-$ and $y \in \Sigma_2^+$, and 
\item no points from $\d_3X$ when $x, y \in \Sigma_2^{\pm}$. In that case, the polarity of  $C$ is the same as the polarity of $x$ and $y$ in $\Sigma_2$.
\end{enumerate}
\end{cor}

{\it Proof.}\; We can assume that  $f|_{\d_2X}$ is Morse and its maxima and minima  alternate.
By Lemma 4.2, one can change the polarity of arcs $C \subset \d_2X$ between consecutive points $a, b$ from $\d_3X$, provided $v_2|_C \neq 0$. Note that the polarity of $a$ and $b$ must be opposite. 
\qed
\smallskip
 
The Morse formula for the vector fields (Theorem 2.1) helps to link the  topology of $\d_1^+X$ with the distribution of arcs from $\d_2^+X$ and  points from $\d_3^{\pm}X$ along its boundary $\d_2X$.

The next lemma is similar in spirit to Theorem 4.8  from [E1].  That theorem is a very special case of the Eliashberg general  theory of folding maps surgery (see [E1], [E2]). However, we cannot  apply Eliashberg's  results directly:  our $v$-generated foldings $p_x: \d_1X \to f^{-1}(f(x))$ have a nice target space only \emph{locally}; a natural target space in our setting is a 2-complex, typically with singularities.

\begin{lem} Let $X$ be a  compact $3$-manifold with a generic nonsingular vector field $v$.
Then  the $v$-generated stratification  $\{\d_k^+X\}_{0 \leq k \leq 3}$ of $X$ has the following properties:
$$\chi(\d_1^+X) - \chi(\d_1^-X) = \#(\d_3^+X) - \#(\d_3^-X) = 2[\#(B\; \mathrm{arcs}) - \#(A\; \mathrm{arcs})]$$
and
$$\chi(\d_1^+X) = \chi(X) + [\#(\d_3^-X) - \#(\d_3^+X)]/2.$$
\end{lem}

{\it Proof.}\; Since $v \neq 0$  in $X$, the index $I(v) = 0$. 
By the Morse formula,
$\chi(X) - \chi(\d_1^+X) + \chi(\d_2^+X) - \chi(\d_3^+X) = 0$. Thus,  
$\chi(\d_1^+X) = \chi(X) + \chi(\d_2^+X) - \chi(\d_3^+X)$.
Loops  in $\d_2^+X$ do not contribute to $\chi(\d_2^+X)$, so   
$\chi(\d_2^+X) = \#(\mathrm{arcs\; in\;} \d_2^+X)$.
Hence, $\chi(\d_1^+X) = \chi(X) + \#(\mathrm{arcs\; in\;} \d_2^+X) -  
\chi(\d_3^+X)$. Note, that the $C$-arcs and their ends
do not contribute to the difference $\#(\mathrm{arcs\; in\;} \d_2^+X) -  
\chi(\d_3^+X)$: such arcs have a single end in
$\d_3^+X$. On the other hand, the $B$-arcs do not contribute to  
$\d_3^+X$. Therefore the difference $\chi(\d_2^+X) - \chi(\d_3^+X)$ is equal
$[\#(A\; \mathrm{arcs}) + \#(B\; \mathrm{arcs})] - 2 \#(A\;  
\mathrm{arcs}) =
\#(B\; \mathrm{arcs}) - \#(A\; \mathrm{arcs}) = [\#(\d_3^-X) -  
\#(\d_3^+X)]/2$.

Recall that for any 3-manifold $X$, $\chi(X) = \frac{1}{2}\chi(\d_1X) = \frac{1}{2}[\chi(\d_1^+X) + \chi(\d_1^-X)]$. Replacing $\chi(X)$ with$ \frac{1}{2}[\chi(\d_1^+X) + \chi(\d_1^-X)]$ in the formulas above,
leads to the relation $\chi(\d_1^+X) - \chi(\d_1^-X) = \chi(\d_3^+X) - \chi(\d_3^-X)$. \qed 
\smallskip

Combining  Lemma 4.3 with Corollary 3.1, we get

\begin{cor} 
$\#(\d_3^+X) = \#(\d_3^-X)$ if and only if $\chi(\d_1^+X) = \chi(\d_1^-X) = \chi(X).$
\smallskip

When $\d_1X$ is connected, by deforming $v$, we can arrange for $\d_1^+X$ to be a 2-disk. For any such choice of Morse data $(f, v)$, 
$$[\#(\d_3^+X) - \#(\d_3^-X)]/2 =  1 - \chi(X).$$
\end{cor}

\begin{cor} 
If a 3-manifold $X$ with $\chi(X) > 0$ admits a nonsingular function $f$ with a $3$-convex  
Morse data, 
 then the  restriction $f|_{\partial_1^+X}$ must have at least $\chi(X)$ extrema.
\end{cor}

{\it Proof.}\; The hypotheses $\partial_3^+X = \emptyset$ implies that only $B$-arcs could 
be present in $\d_2^+X$. The positive contribution to $\chi(\partial_1^+X)$ comes from the 
components of $\partial_1^+X$ shaped as disks. We divide disks into two types: 1) disks
with no $B$-arcs in their boundary (which entirely belongs to $\d_2^+X$ or to $\d_2^-X$) 
and 2) the rest of the disks. Any disk of the first type must contain at least 
one local extremum of $f|_{\partial_1^+X}$. Any disk of the second type contains at least 
one $B$-arc. Since $\partial_3^+X = \emptyset$,  we get  
$\chi(\partial_1^+X) = \chi(X) + \#(B-\mathrm{arcs})$. Now the lemma follows from writing 
down $\chi(\partial_1^+X)$ as the Euler class of all disks of the first type plus the Euler 
class of the rest of $\partial_1^+X$. \qed

\begin{cor} 
For any connected  $3$-manifold $X$ with a connected boundary and Euler number $\chi$, there 
exit nonsingular Morse data $(f, v)$  so that  $\d_1^+X$ is a disk $D^2$.  
For such data, we get $\#(\d_3^+X) \geq   2\chi -  2$ and $\#(\d_3^-X) \geq  2 - 2\chi$. 
As a result, when $\chi > 1$, the disk cannot be convex with respect to the field $v_1$; 
as $\chi$ grows, the disk $\d_1^+X$ becomes more ``wavy". Similarly, when $\chi < 1$, 
the disk cannot be concave with respect to  $v_1$, that is, $\d_3^-X \neq \emptyset$.
\end{cor}
 
{\it Proof.}\; By Lemma 3.1, appropriate deformations of $v$ will  
produce $\d_1^+X \approx D^2$. 
In  view of Lemma 4.3, the claim follows.
\qed

\bigskip

The following theorem shows that convexity of Morse data is equivalent  to the 
possibility of special \emph{ordering} of the $f|_{\d_1X}$-critical points by their 
critical values, and thus, in general, fails.  However, if we formally attach index 
$i + 1$ to each critical points $x \in \Sigma_1^-$ of classical index $i$, the new 
self-indexing of $f|_{\d_1X}$ becomes possible. 

\begin{thm} Let $(f, v)$ be Morse data whose  
restriction $(f_1, v_1)$ to the boundary $\d_1X$ is also of
the Morse type. If $\d_2^+X = \phi$, then there is no ascending  
trajectory $\gamma(t) \subset \d _1X$
of the vector field $v_1$, such that $[lim_{t \rightarrow +\infty}  
\gamma(t)] \in \Sigma_1^+$ and
$[lim_{t \rightarrow -\infty} \gamma(t)] \in \Sigma_1^-$.

Conversely, if no such $\gamma(t)$ exists, one can deform the  
gradient-like vector fields $\{v, v_1\}$
(equivalently, the metric $g$ in which $v = \nabla f$) to a new gradient-like pair $\{\tilde  
v, \tilde v_1\}$ (to a new metric $\tilde g$),
in such a way that, with respect to the new fields, $\d_2^+X = \phi$.   
In particular, if
$f(\Sigma_1^+) < f(\Sigma_1^-)$, then $f$ admits convex Morse data  
(convex metric $\tilde g$).

In contrast, no nonsingular Morse data $(f, v)$ are $2$-concave: $\d_2^+X \neq \emptyset$\footnote{Note that any $X$ admits a field $v \neq 0$ with respect to which $\d_2^+X = \emptyset$ (cf. [BP])}.
\end{thm}

{\it Proof.}\; For a generic metric $g$, consider the vector field $v =  
\nabla f$  and its orthogonal projection $v_1$ on $\d_1X$.
The function $h_1: \d_1 X \rightarrow \R$, defined via the formula $v  
= v_1 + h_1\cdot n_1$, where $n_1$ is the
inward normal field,
has  zero as  a regular value.  Then $\d^+_1X = h_1^{-1}([0,  
+\infty))$,  $\d^-_1X = h_1^{-1}((-\infty, 0])$ and $\d_2X =  
h_1^{-1}(0)$.
Now, if the ascending trajectory $\gamma(t)$ which links  $\Sigma_1^-$  
with $\Sigma_1^+$ exists, it must cross somewhere the
boundary $\d_2X$ of $\d^-_1X$. By definition, such crossing belongs to   
$\d_2^+X$ which must non-empty.

On the other hand, if no such $\gamma(t)$ exists, then we claim the  
existence of codimension 1 closed
submanifold $N \subset \d_1X$, which separates $\d_1X$ in two domains 
$A \supset \Sigma_1^+$ and $B \supset \Sigma_1^-$
($\d A = N = \d B$) and, in addition, has the following  property. The  
vector field $v_1$ is \emph{transversal} to $N$ and points
\emph{outward} of $A$.  Indeed, one can take a small regular  
neighborhood (in $\d_1X$) of the union of descending trajectories of  
all critical points  from $\Sigma_1^+$ for the role of $A$. Here we are 
employing  the fact that no descending trajectory originating at $\Sigma_1^+$
reaches $\Sigma_1^-$.

Since, away from $\Sigma_1^+ \cup \Sigma_1^-$, $v_1 \neq 0$, in the  
tangent space
$T_x$ of $X$ ($x \in \d_1X$), there is an open cone $T^+_{f, x}$  
containing $v_1$ and comprised of gradient-like vectors.

With such a separator $N$ in place, consider a smooth function $\tilde  
h_1: \d_1X \rightarrow \R$ with the properties: 1) zero is a regular
value of $\tilde h_1$ and  $\tilde h_1^{-1}(0) = N$;\; 2) $\tilde  
h_1^{-1}((-\infty, 0) = A$, $\tilde h_1^{-1}([0, +\infty) = B$;\;  3)
$\tilde h_1 = h_1$ in the vicinity of $\Sigma_1^+ \cup \Sigma_1^-$;  
and 4) $v_1 + \tilde h_1\cdot n \in T^+_f$.
Note that the field $\tilde v := v_1 + \tilde h_1\cdot n$ points  
inside $X$ along $A$ and outside along $B$.

Now we can find a metric $\tilde g$, in which $\tilde v$ is the  
gradient of $f$. In $\tilde g$,
$\tilde v$ is orthogonal to the plane $K_x := Ker(df)$. Denote by  
$\tilde v_1$ the $\tilde g$-orthogonal projection
of $\tilde v$ on $T_x$. Since $\tilde h_1$ vanishes on $N$,  $\tilde  
v_1 = v_1$ along $N$, and therefore, is
transversal to $N$ and points outward $A$. As a result, with respect to  
  $\tilde g$, $\d_2^+X = \phi$.

We can deform the original metric $g$ into $\tilde g$, thus deforming  
the gradient fields $v, v_1$ into the gradient fields $\tilde v, \tilde v_1$. 
\smallskip

The last claim follows from the observation that since $v \neq 0$ the absolute maximum (minimum) of $f$ on $X$ must be realized at a point from $\Sigma_1^-$ (from $\Sigma_1^+$). \qed
\bigskip

In view of Theorem 4.1 and Lemma 4.1, we get

\begin{thm} A connected $3$-manifold $X$ with a connected boundary is a  
handlebody if and only if one of the
following properties is valid:
\begin{itemize}
\item $X$ admits a smooth nonsingular function $f$  
whose restriction on the boundary $\d_1X$ is
Morse;  moreover,  no \emph{ascending} trajectory of a gradient-like  
field $v_1$ links in $\d_1X$ a singularity
from $\Sigma_1^-$ to a singularity from $\Sigma_1^+$.\smallskip

\item $X$ admits a smooth nonsingular function $f$   
whose  boundary $\d_1X$ is
convex with respect to a gradient-like field $v$, that is, $\d_2^+X =  
\emptyset$.
\end{itemize}
\end{thm}

Given the remarkable proof of the Geometrization Conjecture [P1], [P2], 
the  proposition below must be viewed just as an illustration.  It shows what advances 
towards the Poincar\`{e} Conjecture are possible by modest means of the Morse Theory alone. 
We get the following criteria for  recognizing standard 3-disks in terms of Morse data: 

\begin{cor}
An  acyclic $3$-manifold $X$  is a $3$-disk if and only if one of  
the following properties is valid:
\begin{enumerate}
\item $X$ admits a smooth function $f$ with no critical points in $X$  
whose restriction on the boundary $\d_1X \approx S^2$ is Morse;  moreover,  no  
\emph{ascending} trajectory of a gradient-like field $v_1$ links in $\d_1X$ a
singularity  from $\Sigma_1^-$ to a singularity from  
$\Sigma_1^+$.\smallskip

\item $X$ admits a smooth function $f$ with no critical points in $X$  
whose boundary
$\d_1X$ is convex with respect to a gradient-like field $v$.\smallskip

\item $X$ admits a smooth function $f$ with no critical points in $X$,  
so that
$\d_1^+X \subset S^2$ is connected and $\d_3^+X = \emptyset$ (i.e. 
the Morse data are 3-convex.)
\end{enumerate}
\end{cor}

{\bf Proof \quad} Evidently, a standard 3-disk admits Morse data with the  
properties described in $(1)$ --- $(3)$.\smallskip

By a homological argument based on the Poincar\`{e} duality,  an acyclic  
$3$-manifold $X$ has a  spherical boundary. \smallskip

By Theorems 4.1, 4.2, properties (1) and (2) are equivalent,  and by Lemma 4.1, (2)  
implies that $X$ is a $3$-disk.
To prove (3) we use the last formula from Lemma 4.3. 
Since $\chi(X) = 1$ and $\d_1^+X$ is connected, it follows that
$\#(\d_3^+X) - \#(\d_3^-X)$ is twice the number of holes in $\d_1^+X$.  
Hence, $\#(\d_3^+X) = 0$ implies that $\#(\d_3^-X) = 0$ and, therefore, $\d_1^+X$ must be a $2$-disk with no points from $\d_3X$ along its boundary. 
However, the boundary of the disk $\d_1^+X$ cannot belong entirely  to  
$\d_2^+X$:  $f$ must attend its absolute maximum and minimum in  
$S^2$. We have seen already that this implies that
$\d_2^-X \neq \emptyset$. Hence, $\d_3X = \emptyset$ which implies that
$\d_2^+X = \emptyset$. Thus, under the hypotheses, the  
$3$-convexity implies the $2$-convexity. In turn, the $2$-convexity  
implies that  the manifold is a $3$-disk.
By Corollary 3.1, we can find Morse data so that $\d_1^+X$ is a disk. \qed 
\bigskip

\begin{exmp}
\end{exmp}
Let $X$ be the Poincar\'{e} homological sphere  
from which a $3$-disk being has been removed. It follows from the 
theorems above that,  for any smooth function $f: X \rightarrow \R$ with no singularities in  $X$, there is an ascending $v_1$-trajectory which links a singularity from  $\Sigma_1^-$ to a singularity from
$\Sigma_1^+$. Also, no nonsingular Morse data $(f, v)$ can insure 
both the connectivity of $\d_1^+X$ and the 3-convexity. Therefore, a ``mild" connectivity restriction on $\d_1^+X$ turns the obstructions to 3-convexity into a topological phenomenon.  


\section{Cascades, $2$-Spines  and Concavity}

\begin{defn} Let $K$ be a finite  two-dimensional polyhedron (cellular complex) imbedded as a subcomplex in a compact 3-manifold $X$ with boundary $\d_1X$. 
We say that $K \subset X$ is a \emph{spine}, if $X \setminus K$ is 
homeomorphic (diffeomorphic) to the product $(\d_1X \setminus (K \cap \d_1X))\times [0, 1)$.
\end{defn}

It follows that $X$ is collapsible onto $K$, in particular, $K$ is a strong deformation retract of $X$. Furthermore, $X$ is homeomorphic to a cylinder of a map $g:  \d_1X \to K$ which is an identity on
$K \cap \d_1X$ (cf. [M], Theorem 1.1.7).
\smallskip

Our next goal is to use nonsingular Morse data $(f, v)$ in order to construct rather special spines that we call \emph{gradient}. 
First, we  focus on the complications arising from the concave locus  
$\d_2^+ X$. It is comprised of a finite
number of disjoint arcs or loops $\{E_j\}$. For each connected curve $E_j  
\subset \d_2^+ X$, consider the set of points in $X$
which can be reached from $E_j$ moving down along the trajectories of  
$-v$. Denote by $W_j$ the closure in $X$ of
this set. We call such a set $W_j$ a \emph{waterfall}. The union  
$\cup_j \; W_j$ of all waterfalls is called a
\emph{cascade} and is denoted by $\mathcal C(\d_2^+ X)$. We denote by 
$\mathcal C(\d_3^-X)$ the finite union of the
downward trajectories through the points of $\d_3^-X$. These trajectories are called \emph{free}.
\begin{figure}[ht]
\centerline{\includegraphics[height=1.5in,width=2.5in]{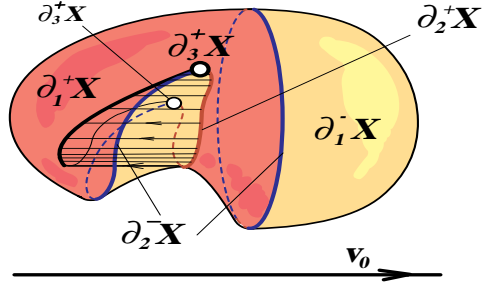}}
\bigskip
\caption{\small{Killing a $1$-cycle in $\d_1^+X$ by attaching the $2$-cell  
$\mathcal C(\d_2^+ X)$.}}
\end{figure}
\begin{figure}[ht]
\centerline{\includegraphics[height=1.4in,width=2.8in]{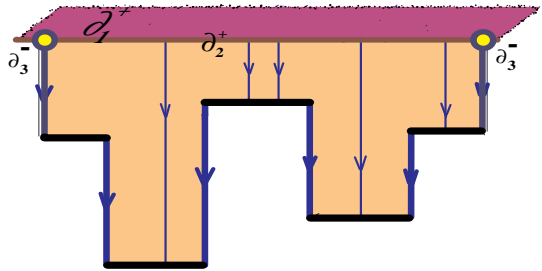}}
\bigskip
\caption{\small{A waterfall of an $A$-arc in $\d_2^+X$.}}
\end{figure}
\begin{figure}[ht]
\centerline{\includegraphics[height=1.4in,width=2.8in]{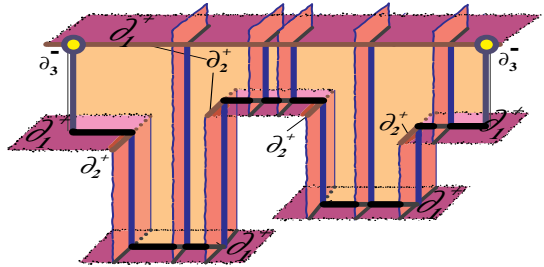}}
\bigskip
\caption{\small{The same waterfall in a cascade.}}
\end{figure}
\begin{figure}[ht]
\centerline{\includegraphics[height=2in,width=3in]{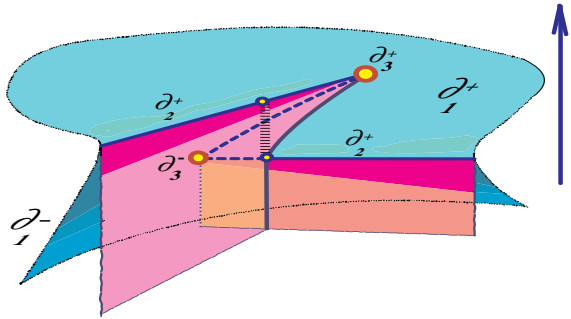}}
\bigskip
\caption{\small{A cascade of a dove tail singularity with two cusps from 
$\d_3^+X$ and $\d_3^-X$. Note the trajectory shared by the two waterfalls.}}
\end{figure}

\begin{defn} We say that the gradient-like field $v$ is  
$\d_2^+$-\emph{generic} if
\begin{itemize}
\item for each $x \in  \d_2^+X$ there is an open interval $V_x$ centered on $x$ such that  the surface $W_x$, formed by the (downward) trajectories through $V_x$, away from $V_x$, has only transversal  
intersections with $\d_2^+X$ ;
\item the downward trajectories of points from the finite set $\d_3 X$  
are all distinct and each trajectory
belongs to a single waterfall $W_j$.
\end{itemize}
\end{defn}

Figures 8, 10, 11 show mechanisms by which $\d_2^+$-generic waterfalls  
are created, as well as their typical shapes.\smallskip

\begin{lem} A small pertubation of a gradient-like field $v$ turns it  
into a  $\d_2^+$-generic field.
For such a field, the downward trajectories of points from $\d_3^-X$  
terminate in the interior of the surface $\d_1^+ X$.
\end{lem}
{\it Proof.}\; For $v$ being  $\d_2^+$-generic is an open dense property established  by standard transversality arguments. \qed
\smallskip

The following proposition describes how any generic $v$ generates a unique gradient spine.
\begin{thm} Let $X$ be a compact $3$-manifold.   We assume that a smooth function 
$f: X \to \R$ has no singularities and  its  
gradient-like field $v$ is $\d_2^+$-generic. Then the $2$-complex
$K = \d_1^+ X \cup \mathcal C(\d_2^+ X)$ is a spine of  
$X$.
\end{thm}

\begin{figure}[ht]
\centerline{\includegraphics[height=2in,width=2.8in]{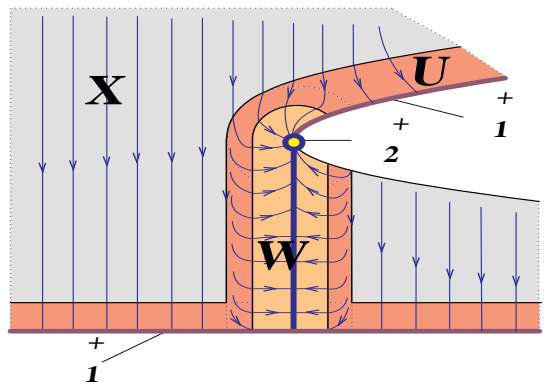}}
\bigskip
\caption{\small{A  section of the retraction $X \to K$ by a plane transversal to
 $\d_2^+X$.}}
\end{figure}

{\it Proof.}\;  We can modify the gradient-like field $v$ by a nonnegative 
conformal factor  $\phi: X \to \R_+$ which is zero on $\d_1^+X$ and 
positive everywhere else.  Abusing notations, denote this modification by $v$ as well. Let $\mathcal U$ be a regular neighborhood of the  
complex $K$ in $X$ which is invariant under
the flow produced by the modified $-v$. The existence of $U$ is implied by the  
invariancy of the set $K$. For example, one can
use a metric in $X$, invariant under the flow, in order to choose  
$\mathcal U$ as an $\epsilon$-neighborhood of $K$.
Then we pick a regular $\epsilon/2$-neighborhood $\mathcal W \subset X$ of the  
cascade $\mathcal C(\d_2^+ X)$. It is properly contained
in $\mathcal U$. Let $w$ be a field  that defines a retraction of $\mathcal U$  
onto $K$. It  vanishes on $K$, is tangent to $\mathcal U \cup \d_1^-X$, as well to $\mathcal W \cup \d_1X$, 
and inward transversal to the portion of $\d \mathcal W$ that is not in $\d_1X$. Moreover, as Fig. 
12 testifies, $w$ can be chosen so that, away from $K$, it is never positively proportional to $v$. 
Consider a smooth partition of unity $1 = \alpha + \beta$, where  
$\alpha$ is supported in
$X \setminus (\mathcal W \cup \d_1^+X)$ and $\beta$ in $\mathcal U$.  
Form the vector field $v = -\alpha v + \beta w$. It is
defined globally and vanishes only on $K$. 
By the construction of $\mathcal U$, $v = -v$ on the boundary of  
$\mathcal U$.  Therefore, $v$ or is tangent to $\d\mathcal U$ or
points inside  $\mathcal U$ (see Fig. 12). As a result,  positive  
$v$-trajectories of  points $x \in \mathcal U$ must reach either the
cascade, or the set $\d_1^+X$. Hence, the $v$-flow governs the  
retraction of $X$ on $K$, and $K$ is a strong deformation retract of $X$.
Moreover, the $(-v)$-flow gives a product structure $\d_1^-X\times[0, 1)$ 
to $X \setminus K$. \qed
\begin{rem}
\end{rem}
For a given field $v$, one can introduce an equivalence relation 
$\sim_{v}$ among points of $X$: two points are defined to be equivalent if they both belong to the closure of the same $v$-trajectory. The quotient space $X/\sim_{v}$ with the quotient topology is called the \emph{orbit-space} of $v$. For a generic $v \neq 0$,  $X/\sim_{v}$ can be given the structure of a two-dimensional CW-complex which is homotopy equivalent to $X$. In fact,  for a generic $v \neq 0$, the obvious maps $P: X \to X/\sim_{v}$ and 
$p: \d_1^+ X \cup  \mathcal C(\d_2^+ X) \to X/\sim_{v}$ are  Serre fibrations. Thus, for such $v$, the gradient spine can be also regarded as a homotopy substitute for the space of $v$-trajectories.  \qed
\begin{rem}
\end{rem}
The 2-cells  in which $K$ is  subdivided admit a preferred  
orientation induced by the preferred orientation of $\d X$, 
so that they form an integral $2$-chain $[K]$. Its  boundary $\d[K]$ consists of the $1$-chain $\d_2^- X \cup \mathcal C(\d_3^-X)$ together with the singularity locus  
$s(K)$ of $K$. This locus is comprised of curves from the intersection 
$\d_1^+ X \cap  \mathcal C(\d_2^+ X)$ and the orbits shared by pairs of waterfalls. In short, each edge from the support of $\d[K]$  contributes to the cycle $\d [K]$ with multiplicity $\pm 1$ (and not $\pm 3$)---an important property which the gradient spines share with the branched spines (see Corollary 3.1.7, [BP]) and which originally was studied  in [GR] and [GR1] .
\begin{rem}
\end{rem}
 Note that changing $f$ to $-f$ and $v$ to $-v$ exchanges the strata $\d_1^+X \Leftrightarrow \d_1^-X$, $\d_3^+X \Leftrightarrow \d_3^-X$ and keeps the strata $\d_2^+X, \d_2^-X$ fixed. Therefore, the gradient spines $K(f, v)$ and $K(-f, -v)$ ``complement" each other in $X$: they share $\d_2X$ and their cascades $\mathcal C(f, v)$ and  $\mathcal C(-f, -v)$ complement each other in the set spanned by all the trajectories through $\d_2^+X$. Evidently, the topologies of  $K(f, v)$ and $K(-f, -v)$ could be radically different. However,
their  complexities (see Definition 7.1) are equal.
\bigskip
\begin{figure}[ht]
\centerline{\includegraphics[height=2.7in,width=3.5in]{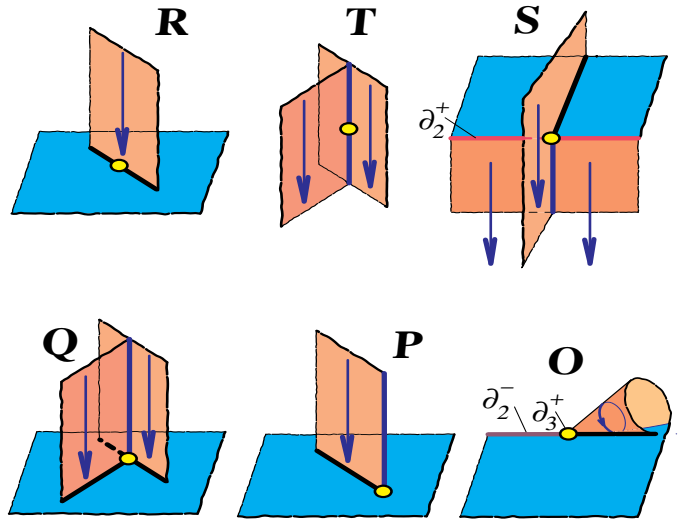}}
\bigskip
\caption{\small{Singularity types of  gradient spines.}}
\end{figure}

By its very construction, any generic gradient spine $K$ has a very particular  
stratified local geometry. In fact, as abstract cellular $2$-complexes, gradient spines
are  \emph{simple} spines in the sense of  [M]. However, in our context, (as well as in the case of branched spines) \emph{the smooth 
structure of} $\d_1^+X$  and waterfalls \emph{breaks the $3$-fold symmetry} of simple combinatorial spines: in other words, we distinguish between the ``$Y$-" and ``$T$-shapes".\smallskip

The stars of points in $K$ are shown in Fig. 13 and will be referred by their letter labels. The six types of links of points in $K$ are depicted in Fig. 19.
Each point from the set $\d_2^-X \cup_{\d_3^-X} \mathcal C(\d_3^-X)$ has  
a neiborhood  shaped as a half-disk and a link of type (2). The singular set $s(K)$ consists of points of types $R, T,$  and $S$, all having  links of type (5). The $T$-type is  
produced by the trajectories that belong to two distinct waterfalls or to  two branches of
the same waterfall. The $R$-points are  generic to loci where waterfalls hit the ground $\d_1^+X$.
The $S$-points are  generated when a waterfall transversally hits an arc from $\d_2^+X$.  
They are hybrids of $T$ and $R$ types.  Topologically $R, T,$  and $S$ types are indistinguishable. 
The $Q$-type is  generated where two waterfalls hit the ground. Stars of  $Q$-points 
are  shaped  as a union of a disk with a half-disk with a quoter-disk, and their links are of type (6). 
The $Q$-type singularities are isolated in $K$. 
Singularities of the $P$-type are located where a free trajectory  
through $\d_3^-X$ hits the ground. Hence, they are in $1$-to-$1$  
correspondence with points of $\d_3^-X$. A neighborhood of a  
$P$-singularity is a union of a disk with a quater-disk which share a  
common radius;  its link is of type (4). Finally, the singularities of the $O$-type are just points from  
$\d_3^+X$. They also have stars  shaped as cones over a 
circle with a radius (see Fig. 11).  Topologically  the $O$ and $P$-types are the same.
\bigskip

Next, we prove that 3-manifolds have gradient spines which are rather special \emph{"origami"} folded from 2-disks (see Fig. 14). This result is very similar to Theorem 1.2 from [I], where an origami is built from a normal pair. However, one technical difference is evident: in [I] the $v$-flow is transversal to the normal disk, while in our construction, the gradient flow is not necessarily transversal to $\d_1^+X$ along its boundary $\d_2X$.

\begin{thm}
The spine  $K =  \d_1^+ X \cup \mathcal C(\d_2^+ X)$ of a 3-manifold $X$, produced from Morse data as in Theorem 5.1 is an image under a cellular map  
$\Phi: S \to K$ of a cellular 2-complex $S$ homeomorphic to the surface $\d_1^+X$.  The map $\Phi$ is $1$-to-$1$ in the interiors of the 2-cells in which $S$ is subdivided, at most  $2$-to-$1$ on the $1$-skeleton without vertices, and at most $3$-to-$1$ on the set of vertices. The local geometry of $\Phi$ can be described by the four identification patters  in Fig. 15.

Moreover, any 3-manifold $X$ has a gradient spine $K$ which is an image of a 2-disk $D^2$ under a cellular map  $\Phi: D^2 \to K$ which is $1$-to-$1$ in the interiors of the 2-cells in $D^2$, at most  $2$-to-$1$ on the $1$-skeleton of $D^2$ without vertices, and at most $3$-to-$1$ on the set of vertices.
\end{thm}

\begin{figure}[ht]
\centerline{\includegraphics[height=2.5in,width=3.5in]{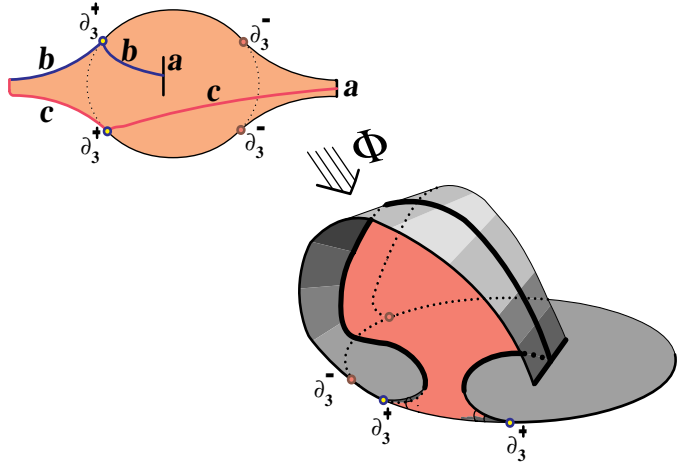} 
}
\bigskip
\caption{\small{An ``origami" map $\Phi: D^2 \rightarrow K$ with a collapsable $K$.}}
\end{figure}

\begin{figure}[ht]
\centerline{\includegraphics[height=2.0in,width=2.5in]{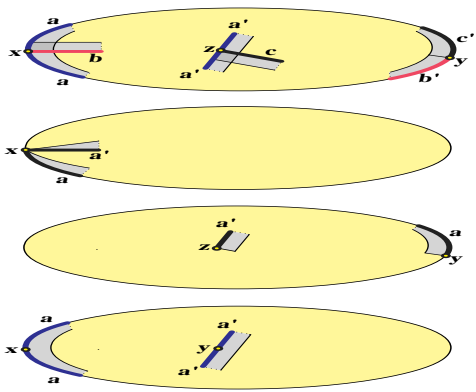}}
\bigskip
\caption{\small{Identification patterns for gradient spines.}}
\end{figure}

{\it Proof.}\;
One can resolve $K$  along its singular set  $s(K)$ so that the resulting space is  a  
celluar 2-complex $S := \d_1^+X \cup_{\d_2^+X} (\coprod_j W_j)$.   
It is obtained from  $\d_1^+X$ by attaching induvidual waterfalls $\{W_j\}$ along the loops  
and arcs forming $\d_2^+X$. As a result, $S$ is homeomorphic to   
the surface $\d_1^+X$. In particular, when $\d_1^+X$ is a 2-disk, so is the resolution  
$S$.

This resolution can be done by performing cuts of several types. First,  
at each singularity of any type, but the $T$-type, (see Fig. 13, 16) the  
cut separates the cascade from the ``ground" surface $\d_1^+X$. 
Then at each singularrity of the $T$-type (see Fig. 17),  
there is a preferred ``half-waterfall" which is separated from the adjacent waterfall 
by a cut along a trajectory that they both share.
\begin{figure}[ht]
\centerline{\includegraphics[height=1.8in,width=3in]{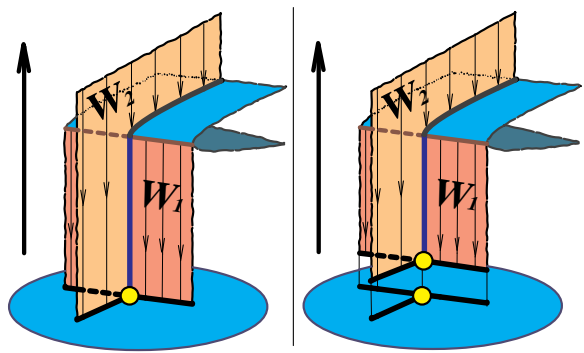}}
\bigskip
\caption{\small{A singularity of the  $Q$-type and its resolution.}}
\end{figure}
\begin{figure}[ht]
\centerline{\includegraphics[height=1.8in,width=3in]{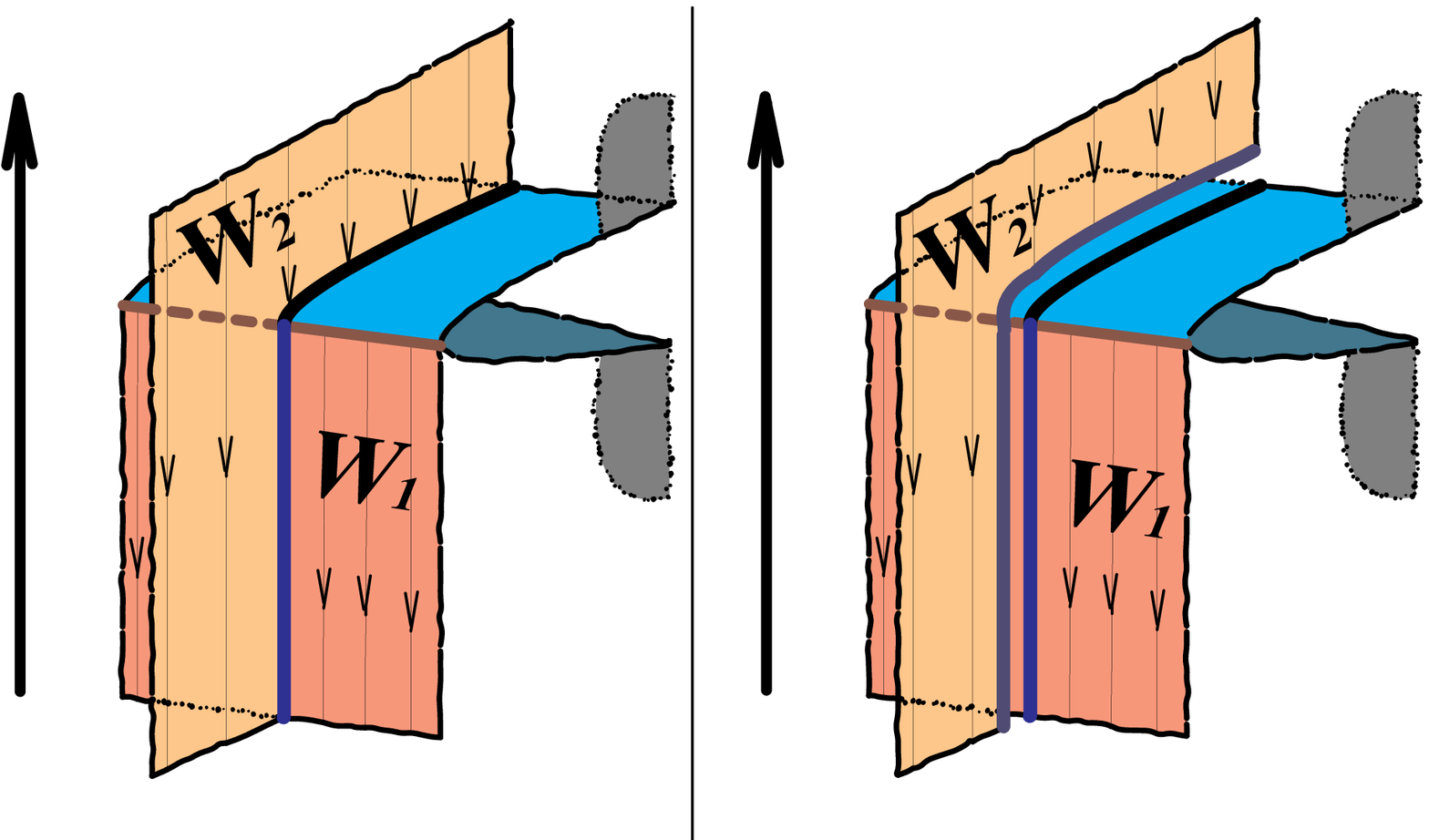}}
\bigskip
\caption{\small{A singularity of the  $T$-type and its resolution.}}
\end{figure}
\begin{figure}[ht]
\centerline{\includegraphics[height=1.8in,width=3in]{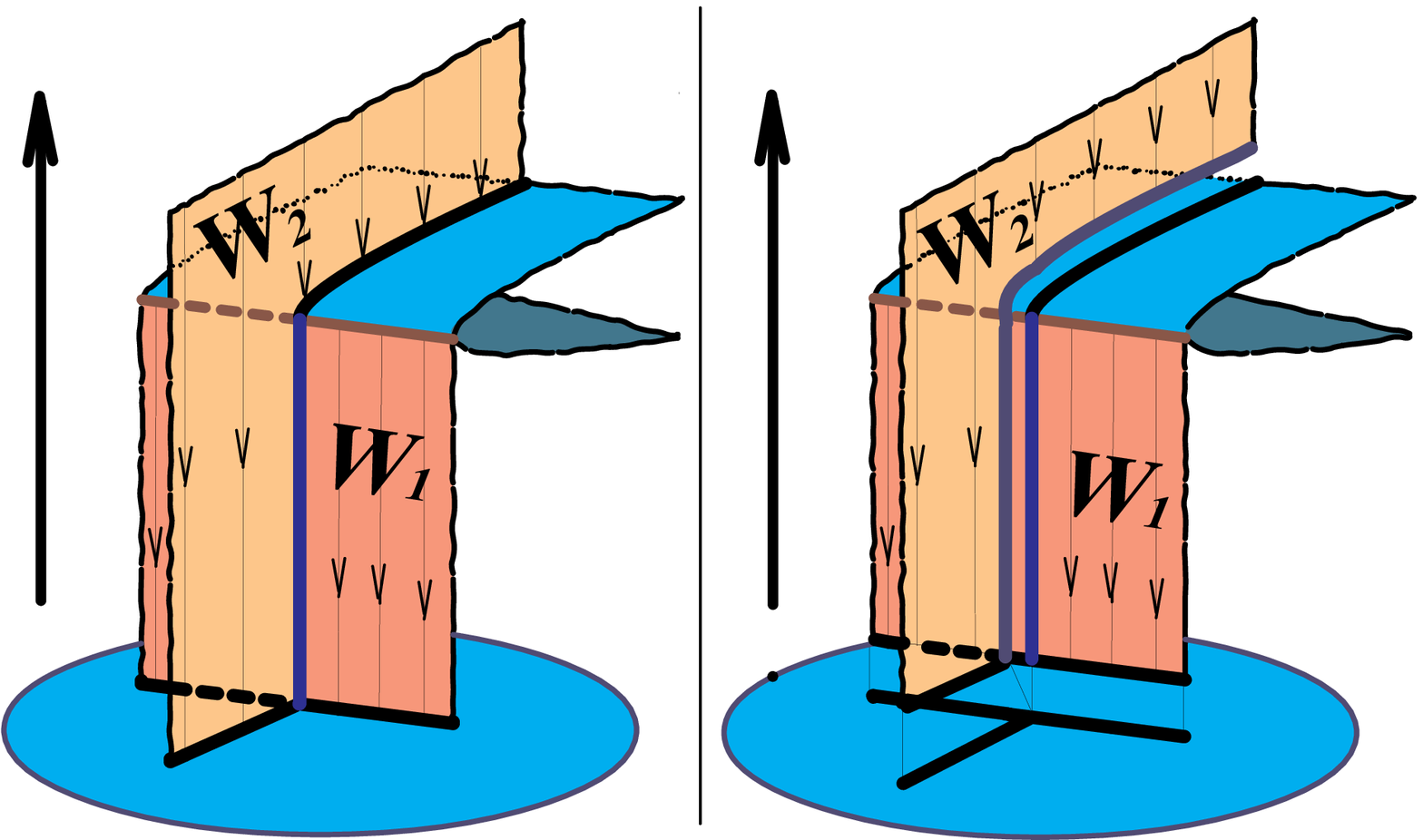}}
\bigskip
\caption{\small{A complete resolution of singularities.}}
\end{figure}

The shape of the waterfalls $W_j$'s depends on the type of the arc or  
loop in $\d_2^+X$  from which it falls. Recall that, after an  
appropriate change of the vector fields, the arcs come in three  
flavors: $A, B$ and $C$. The $A$-type waterfall has two ``free" edges which are not affected by the gluing $\Psi$.  The edges emanate from the two ends 
$x, y \in \d_3^-X$ of an $A$-arc. The $B$-type waterfall has no ``free" edges at all. It falls from a $B$-arc whose eds are in $\d_3^+X$. The $C$-type waterfall has one ``free" edge which emanates from a point in $\d_3^-X$---the end of the arc. 

By Corollary 3.1, for  appropriate Morse data, we can assume that $\d_1^+X$ is a disk.

Evidently, there is a mapping $\Phi: S \rightarrow K$ which glues  
the complex $K$ back. The map $\Phi$ is $1$-to-$1$ on the compliment to the $1$-skeleton of  $S$, it is at most $2$-to-$1$ on the interior of its edges and at most $3$-to-$1$  on its vertices (see Fig. 15---18). With the local topology of $K$ being restricted by the list in Fig. 13, one can verify that Fig. 15 lists all possible gluing patterns for $\Phi$ (see Appendix A in [K] for details).
\qed
\smallskip

Globally, the gluing patterns of $\Phi$ can be quite intricate as shown in Fig. 19 which depicts an example of a map  $\Phi: D^2 \rightarrow K$ as in Theorem 5.2 with a contractible, but not collapsable $K$.  The second  diagram in Fig. 19 is obtained from the first one by  elementary collapses  performed through the free trajectories. The third diagram suggests  contracting the shaded disk to a point, thus ignoring the subtleties of  the $\Phi$-gluing pattern inside the disk, but keeping the  presentation.
\begin{figure}[ht]
\centerline{\includegraphics[height=2.5in,width=3in]{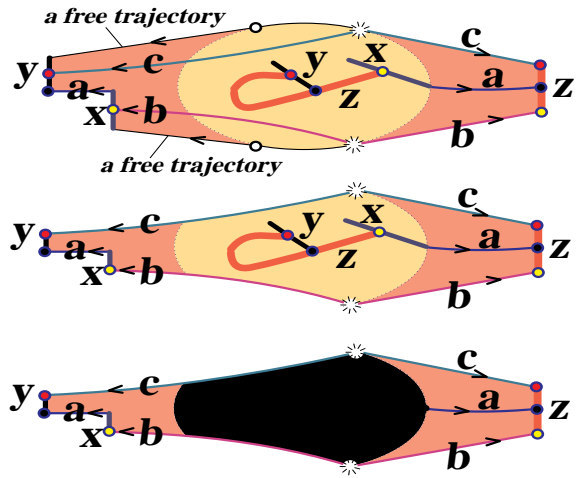}}
\bigskip
\caption{\small{The combinatorics of a map $\Phi: D^2 \rightarrow K$ which  
gives rise to a non-collapsable $K$ with $\pi_1(K) = 1$ presented as  
 $\{a, b, c; \; \hfill\break ba = c^{-1},\, ab^{-1} = 1,\, ac^{-1} = 1\}.$}}
\end{figure}
\bigskip

When $\d_1^+X$ is a 2-disk, Morse data (equivalently, gradient spines $K$) produce a \emph{presentation} of the  fundamental group $\pi_1(X) \approx \pi_1(X/\d_1^+X)$. 
Its generators are the oriented $v$-trajectories  that are shared by pairs of waterfalls, or by two branches of the same waterfall, together with the free trajectories that emanate from $\d_3^-X$ (each waterfall has at most two free trajectories, and two waterfalls can share a number of trajectories). Let $\{\omega_{i\beta}\}_\beta$ be the set of shared trajectories in $W_i$ and $\{\alpha_{i \delta}\}_\delta$ the set of free ones. Evidently, the homotopy class of each  loop $\gamma$ in $K/\d_1^+X$ is characterized by its traces in the waterfalls (see Fig. 10 and 19). In each waterfall $W_i$, the trace of $\gamma$ can be replaced by a word in $\{\omega_{i\beta}\}_\beta$ and $\{\alpha_{i \delta}\}_\delta$.  The relations are produced by looking independently at each waterfall marked with its shared trajectories.  If we cut $W_i$ along the $\{\omega_{i\beta}\}_\beta$'s, it will break into a number of polygons. Each polygon has certain shared and free trajectories in its boundary, and the rest of the boundary consists of arcs that belong to $\d_1^+X$. A free trajectory belongs to a single polygon, and each polygon has two free trajectories at most. Each polygon is oriented and contributes a single relation: moving along its oriented boundary produces a word in the alphabet that contains free and shared trajectories and their inverses (the arcs of the  boundary that were attached to $\d_1^+X$ are ignored). In fact, when a polygon contains a free trajectory, we can delete it from the list of generators and the polygon itself from the list of relations. Eventually, this will eliminate all  free trajectories with the exception of the pairs that belong to a single waterfall free of shared trajectories in its \emph{interior}. One of the free trajectories in each of such pairs can be dropped from the list of generators and its polygon from the list of relations. 

In general, the same recipe produces a presentation of $\pi_1(X/\d_1^+X)$. Hence,
\begin{thm} Let $X$ be a connected $3$-manifold.  
Then generic Morse data  $(f, v)$ such that $\d_1^+X$ is   
a $2$-disk  give rise to a finite presentation of $\pi_1(X)$.  Thus,  $(f, v)$ determine the $Q^{\ast \ast}$-class (as defined in {\rm [Me]})  of that presentation.  By {\rm  [R]},  this $Q^{\ast \ast}$-class is a topological invariant of  $X$. 
\end{thm}


\section{Abstract gradient spines}

We already noticed that the geometry of gradient spines $K$ provides us with a particular way of orienting their 2-cells. This orientation  is induced by a  preferred orientation of $\d_1X$\footnote{When $X$ is oriented, the preferred orientation of $\d_1X$ is naturally induced.} and is spread along the cascade by the $v$-flow. Furthermore, not only the surfaces $K^\circ$ acquire a \emph{preferred orientation}, but they also have \emph{preferred sides} in the ambient $X$, the sides picked by the inner normals to $\d_1^+X \subset X$.  The inner normals uniquely extend by continuity to each waterfall $W$ and pick its side in $X$.

These properties of gradient spines can be captured in the notion of an \emph{abstract gradient spine}. We start with a few preliminary definitions. 

\begin{defn}
A \emph{simple polyhedron}  is a compact 2-dimensional polyhedron such that each its point has a link homeomorphic  to one of the  shapes in Fig. 20\footnote{Actually, type (3) does not appear in generic gradient complexes, but is present in their 1-parameter deformations.} and the corresponding star from  Fig 13.
\begin{figure}[ht]
\centerline{\includegraphics[height=0.65in,width=3in]{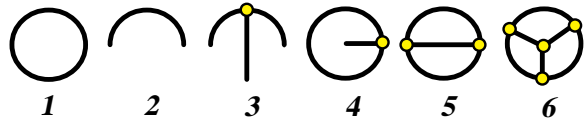}}
\bigskip
\caption{\small{Link types of points in simple polyhedrons.}}
\end{figure}
 \end{defn}
The six types have a partial order induced by the inclusions of closures of the appropriate strata: 
$$(1) > (2),\; (1) > (5), \; (2) > (4), \; (5) > (3).$$

Note, that generic gradient spines are simple polyhedra. However, typically they admit $2$-collapses; thus, the need for considering link types $(2), (3)$, and $(4)$.  \smallskip

The ``topological boundary" $s_\bullet(K)$ of $K$ is formed by points with links of type (2) degenerating into types (3) and (4). The graph $s_\dagger(K)$ is formed by the points with  links of type (5) degenerating into types (3), (4), and (6). Thus a generic point from $s_\bullet(K)$ belongs locally to the boundary a single surface in $K$, while a generic point from $s_\dagger(K)$ belongs to the boundaries of three surfaces. Let $s(K) = s_\bullet(K) \cup s_\dagger(K)$. The finite set of points of  types (3), (4), and (6) is denoted $ss(K)$. The  bivalent vertices from $ss_{\bullet\dagger}(K) := s_\bullet(K) \cap s_\dagger(K)$ are of type  (4), and the rest of the vertices from $ss(K)$ are the four-valent ones of type (6). The set  formed by the points of type (6) is also denoted $Q(K)$.
\bigskip

Recall that a simple spine $K$ is call \emph{special} in [M] or \emph{standard} in [BP] if the stratification $ss(K)   \subset s(K)  \subset K$ gives $K$ a structure of cellular 2-complex, i.e. if $K^\circ := K \setminus s(K)$ is a disjoint union of open 2-cells and  $s(K)^\circ := s(K) \setminus ss(K)$ is a disjoint union of open 1-cells. \bigskip

If $K$ is a simple special spine of $X$ whose points have local models of types (1), (5), and (6), then $X$ can be uniquely reconstructed just from a regular neighborhood  $N(s(K)) \subset K$ of the singular set $s(K)$ ([M], Theorem 1.1.17).  For the reconstruction to work, it is important  that $K \setminus s(K)$, the disjoint union of disks, does not support  non-trivial line bundles. For any gradient spine $K$, the  bundle, normal to $K \setminus s(K)$ in $X$, is also trivial.  This property of gradient spines $K$ will permit a unique reconstruction of $X$ from $K$ as well. 
\bigskip 

In order to distinguish intrinsically the $T$-shaped configurations from the $Y$-shaped ones, we use a particular \emph{system of markers} placed along the edges of the graph $s_\dagger(K)$,  $K$  being a simple polyhedron. The marker 
is a short segment emanating from a generic point $x \in s_\dagger(K)$. It is transversal to $s_\dagger(K)$ and is contained in 
one of the three surfaces (pages) that join at $x$. We call such segments $T$-\emph{markers}. A $T$-marker $m$, the vertical leg of letter $T$, tells us that the two pages that do not contain $m$ are thought ``to form $180^\circ$ angle" in the ambient $X$. In the category of \emph{branched} spines $K$, the local geometry of $K$ along $s(K)$ also picks one page out of three: recall that only two out of three tangent pages  form a cusp, and the preferred  page is the third one [BP]. It suffices to place a single $T$-marker at each edge or loop of the graph $s_\dagger(K)^\circ := s_\dagger(K) \setminus ss(K)$: by continuity, the marker spreads itself along the edge or loop until it reaches an isolated singularity  from $ss(K)$. There the marker's pattern  requires an additional explanation to be provided below.


The $N$-marker is attached to  generic points $x \in s_\dagger(K)$  and is contained in one of the two pages that  do not carry a $T$-marker. 
Informally, one can think of the pages with $T$-markers as waterfalls.  The $N$-markers reside both in the "ground" and in  the waterfalls,  and can be regarded as substitutes for the preferred normals to the oriented pages that contains  $T$-markers.\smallskip

For gradient spines $K$, $v$ provides us with a $T$-marker at points of $\d_1^+X \cap s(K)$. For points on a trajectory shared by two waterfalls, the $T$-marker resides in the page complementary to the two pages that form $180^\circ$ angle. For each waterfall $W$, the vectors tangent to $ \d_1^+X$ and pointing to the preferred side of $W$, produce the rest of $N$-markers. Effectively, the $T$ and $N$-markers generate a \emph{coloring} of $\mathcal N(s_\dagger(K)) \setminus \mathcal N(ss(K))$ with three colors, the colors of the tree pages which share an edge being distinct. Here $\mathcal N(\sim)$ denotes a regular neighborhood of an appropriate set in $K$. 
\begin{figure}[ht]
\centerline{\includegraphics[height=1.3in,width=1.3in]{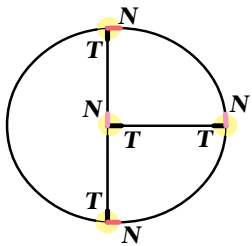}}
\bigskip
\caption{\small{The $TN$-markers on the link of a $Q$-singularity.}}
\end{figure}

A single pair of $TN$-markers approaches each singularity of types (3) or (4). In fact, the markers distinguish between the spines of type (4) in the vicinity of a point from $\d_3^+X$ (O-type in Fig. 13) and a point $y$ where a trajectory through  $\d_3^-X$ hits the ground $\d_1^+X$ (P-type in Fig. 13): in the first case, moving along the loop in (4) in the direction of $N$, we return from the direction marked by $T$; in the second case, moving along the loop in the direction of $N$, we return from the unmarked direction.

Four pairs of $TN$-markers approach each  singularity  $x$ of the type (6) from four different directions. The link $L_x$ of $x$ in $K$ is  shown in Fig. 21, and the markers reflect the ``plane - union half plane - union quoter plane geometry" (as seen from the North Pole). Among  the six edges of $L_x$, there is a single edge $\a$ with two $T$-markers and a single edge $\b$ with two $N$-markers, $\a$ and $\b$ sharing a  vertex. The rest of the markers are determined by  Fig. 21. 
Therefore, the $TN$-pattern in $L_x$ is completely determined by an ordered pair of edges $(\a, \b)$ that share a vertex. Hence there are $4\times C_3^2 = 12$ ways of attaching a pattern of $TN$ markers as in Fig. 21 to a complete graph on four vertices. On the other hand, there are $(3!)^4 = 1296$ ways to mark the four $Y$-shaped beams that join at the singularity $x$. So, the majority of the four beam patterns will not match the local geometry of a $Q$-singularity. 
\smallskip
\begin{defn}
An \emph{abstract $TN$-polyhedron} is a simple  polyhedron $K$  with $T$ and $N$-markers along the edges of the graph $s_\dagger(K)$. We insist that, at each vertex $x$  whose link is a complete graph in four vertices, these markers satisfy the Fig. 21 pattern.  
\end{defn}
Given an abstract $TN$-polyhedron $K$, using the $T$-markers, we can cut it open along $s_\dagger(K)$, so that locally each $T$-marked page is separated from the rest. For the resolution to work, it is important that the $TN$ markers at the isolated singularities will obey the combinatorial rules depicted in Fig. 21.
The result of this $T$-resolution is a compact surface (not necessarily connected) which we denote by $res_T(K)$. It is equipped with the canonical 2-to-1 map $\Phi: res_T(K) \to K$ whose generic fiber is of cardinality 1.

\begin{defn} An abstract $TN$-polyhedron  $K$ is said to be oriented if its resolution $res_T(K)$ is an oriented compact surface.
\end{defn}
\begin{defn} An abstract gradient spine is an abstract oriented $TN$-polyhedron.  
\end{defn}
Hence, any generic gradient spine is an abstract  gradient spine. 

\smallskip

\bigskip

%

%
\begin{thm} If two compact 3-manifolds $X_1$ and $X_2$ with boundaries have homeomorphic gradient spines $K_1$ and $K_2$, the homeomorphism $h: K_1 \to K_2$ being a diffeomorphism along the strata of the spines and respecting the  $TN$-markers, then the manifolds are diffeomorphic.
\end{thm}

{\it Proof.\;} The arguments are similar to the ones used in [M], Lemma 1.1.15 and Theorem 1.1.17. Throughout this proof,  an "embedding" of a stratified space $Z$ in a smooth manifold means an embedding which is an immersion of the smooth strata of $Z$ and which preserves their transversality. 

Let $D_+^3$ be the half-disk $\{x^2 + y^2 +z^2 \leq 1,\; z \geq 0\}$ with the equatorial disk $D^2 = \{z = 0\}$. Let $\mathcal Q \subset D_+^3$ be the union of $D^2$ with the half disk $D^2_+ = \{y = 0,\; z \geq 0\} \cap D_+^3$ with the quoter disk $D^2_{++} =\{x =0,\; y \geq 0, \; z \geq 0\} \cap D_+^3$. Up to diffeomorphisms of $D_+^3$ that preserve $D^2$, there is a unique way to embed the pattern $\Gamma_4$ in Fig. 21 (without the $TN$-markers) in $\d  D_+^3$ in such a way that the circle in $\Gamma_4$ is mapped onto $\d D^2$ and the rest of the graph into $\d  D_+^3 \setminus D^2$. 

Note that map $h$ must not only respect the stratifications in $K_1$ and $K_2$ by the link types in Fig. 20, but $h$ also discriminates between the  isolated singularities of types $O$, $P$, and $Q$ in Fig.13. Indeed, the patterns of $TN$-markers are different for the $O$-cusps and the $P$-intersections of free trajectories: in the case of cusps, moving from the singular point $y$ of the link (see Fig. 20, pattern (4)) in  the $N$-marked direction one returns to $y$ through the $T$-marked direction, while in the $P$-case, leaving $y$  in  the $N$-marked direction results in the return to $y$ through the unmarked direction.  Gradient spines of $O$ and $P$-type singularities admit preferred embeddings in $D_+^3$.  In the  case of a cusp,  $D^2 = \d_1^+X \cup \d_1^-X$ and the waterfall is a triangle whose interior resides in $Int(D_+^3)$ and whose two sides are attached to $D^2$; one of the two sides is realized as the arc $\d_2^+X$, along the other side, the triangle is transversal to $D^2$. Let $\mathcal O$ be the gradient spine of a cusp, that is,  sector $\d_1^+X$ union with this triangle.    In the case of a $P$-singularity, the germ of the spine can be identified with $\mathcal P := D^2 \cup D^2_{++} \subset D_+^3$.
Note that each of the preferred embeddinds $\mathcal Q, \mathcal O, \mathcal P \subset D_+^3$ admits two distinct $N$-markings  that pick the preferred side of the waterfall (they are mirror images of each other); the $T$-markings are uniquely determined by the geometry of the embeddings.   In what follows, we fix one of the two choices for the $N$-markers in $\mathcal Q, \mathcal O, \mathcal P$.
 
For each isolated singularity  $x \in ss(K_1)$ and $h(x) \in ss(K_2)$, consider sufficiently small regular neighborhoods $U_x \subset X_1$ and $V_{h(x)} \subset X_2$. Depending on the type of $x$, both pairs $(U_x,\, K_1 \cap U_x)$ and $(V_{h(x)},\, K_2 \cap V_{h(x)})$ are  diffeomorphic to one of the three models $(D_+^3, \mathcal Q)$, $(D_+^3, \mathcal O)$, and  $(D_+^3, \mathcal P)$ via the diffeomorphisms   which respect the markings. Thus, there exist a diffeomorphism $\a: U_x \to V_{h(x)}$ which maps $K_1 \cap U_x$ to $K_2 \cap V_{h(x)}$ and respects the markings.  Locally (at $x$) the composition $h^{-1}\circ \a : K_1 \cap U_x \to K_1 \cap U_x$ can be represented by the germ of a vector field $w$ tangent to $K_1$ (integrating $w$ over a unit of time produces $h^{-1}\circ \a$). This field $w$ extends to a field $\hat w$, defined in some regular neighborhood  of of $x$ in $U_x$. We use $\hat w$ to define a germ of a diffeomorphism $\b: X_1 \to X_1$ at $x$. Evidently, the germ of $\a \circ \b^{-1}$ is an extension  $\hat h$ of $h$ into a regular neighborhood of $x$ in $X_1$. Let $\hat U_x$ be a regular neighborhood of $x$ contained in the domain of $\hat h$ and $\hat V_{h(x)} := \hat h(\hat U_x)$. Therefore, we have constructed a  diffeomorphism $\hat h: K_1 \cup (\cup_x \hat U_x) \to K_2 \cup (\cup_x  \hat V_{h(x)})$ of stratified spaces which preserves their $TN$-markings. \smallskip

Let $D^0 := (0, 0, 0)$,\; $D^1 := D^2 \cap D^2_+$, and $D^1_+ := D^2_+ \cap D^2_{++}$.\smallskip

In a similar way, after "fattening" of $h$ in the vicinity of $ss(K_1)$,  we can extend $\hat h$ into a regular neighborhood $\mathcal N(s(K_1)$ of $s(K_1) \subset X_1$, while respecting the markings in the source and the target. To accomplish this we need the following  models:
\begin{itemize}
\item $\big\{(D^1\cup D^1_+)\times [0, 1]\big\} \cup D^2_+\times \{0\} \cup D^2_+\times \{1\}\; \subset D^2_+ \times [0, 1]$
\item $\big\{(D^1\cup D^1_+)\times S^1 \big\} \; \subset D^2_+ \times S^1$
\item $D^0 \times [0,1]\; \cup D^2_+ \times \{0\}\; \cup D^2_+ \times \{0\} \subset D^2_+ \times [0,1]$
\item $D^0 \times S^1 \subset D^2_+ \times S^1$
\end{itemize} 
The first couple models the vicinity of an arc or a loop from $s_\dagger(K_i), \;(i =1,2)$, the second one of an arc or a loop from $s_\bullet (K_i)$.\smallskip

The third, most problematic,  extension of $h$ occurs into a regular neighborhood $W$ of $K_1^\circ : = K_1 \setminus \mathcal N(s(K_1)$ in $X_1\setminus \mathcal N(s(K_1)$. It is possible because the surfaces from $K_i^\circ$ have preferred normals, which results in the normal bundles $\nu(K_i^\circ, X_i)$ being trivial. This crucial observation is valid due to the gradient nature of the two spines.\smallskip

Recall that any smooth regular neighborhood of a spine $K \subset X$ is diffeomorphic to $X$; thus, $X_1$ is diffeomorphic to $X_2$.\qed

\begin{lem} Let $G_4(c)$ denote the number of connected regular four-valent  graphs with $c$ vertices, taken up to a homeomorphism. Then the number of connected \emph{special} abstract gradient spines  whose points are of the types $(1), (5)$, and $(6)$ does not exceed $G_4(c)\cdot 12^c$.\footnote{By Theorem 7.3,  the same number $G_4(c)\cdot 12^c$ gives an upper bound on the number of irreducible and boundary irreducible with no essential annuli 3-manifolds.} In turn, $G_4(c)$ can be crudely estimated from above by the number of elements in the symmetric group $S_{4c}$ that move every symbol in $(1, 2, 3, \dots, 4c)$.
\end{lem}

{\it Proof.}\; A special spine $K$,  with all its points modeled after types (1), (5), and (6), is completely determined by a regular neighborhood $\mathcal U$ of its one-skeleton, a regular 4-valent graph $\Gamma$; to reconstruct $K$ from $\mathcal U$ we just attach a disk to every circular component of $\d \mathcal U$. 
There are $G_4(c)$ such graphs $\Gamma$. Each edge $\g \subset \Gamma$ is a core of a beam $B_\g \subset \mathcal U$ with a $Y$-shaped section. Let $V$ be the disjoint union of $c$ copies of the star $V_\star$ of a $Q$-singularity. With the $TN$ markers on both ends $a, b \in \d V$ of the beam in place, intrinsically, there is a unique way to attach the beam $B_\g$ to $V$. Therefore, \emph{any} $TN$ pattern as in Fig. 20 assigned to each copy of $V_\star$ in $V$ will determine the rules for attaching the beams, and thus the reconstruction of $K$. Since $V$ supports $12^c$ $TN$-patterns, the total number of special abstract gradient spines does not exceed $G_4(c) \cdot 12^c$. 

Fix a graph $\Gamma$ as above. Then, $\mathcal U$ is determined by assigning, for each edge $\g$ of $\Gamma$ a pairings between the two $Y$-shaped plugs  in $V$ that correspond to $\g$. For each $g$, there are six pairings, so that the total number of non-homeomorphic $U$'s is bounded from above by $6^{2c}$. This number $2^{2c}\cdot 3^{2c}$ should be compared with the estimate $12^c = 2^{2c}\cdot 3^{c}$ of the $\mathcal U$'s that are consistent with the abstract gradient spine structure. \qed
\bigskip

For a simple spine $K$, consider the  2-chains $C_2^\Box(K)$ that are combinations of the fundamental cycles of connected surfaces that form $K^\circ$, the coefficients in the combinations being $\pm 1$. Also consider the 1-chains $C_1^\Box(K)$ that are combinations of the fundamental cycles of arcs  that form $s(K)^\circ : =  s(K) \setminus ss(K)$, the coefficients in the combinations being $\pm 1$. Following [GR1], [GR2], and [BP], we introduce the following notion:

\begin{defn} An oriented branching on a simple complex $K$ is a 2-chain $\a \in  C_2^\Box(K)$ such that its boundary $\d\a \in C_1^\Box(K)$.
\end{defn}

In other words, an oriented branching is a  special choice of orientations for each of the components of $K^\circ$; note that, for an arbitrary $\a \in  C_2^\Box(K)$, some edges of $s(K)$ can contribute to $\d\a$ with multiplicity $\pm 3$. 
\bigskip 

Next, we introduce the notion of $\vec Y$-\emph{structure} for spines $K \subset X$. It resembles to the notion of branched spines (see Definition 6.6 and Lemma 6.3) and  plays a significant role in the sections to follow.  \smallskip

Consider a configuration $\mathcal Y$ in $\R^3$ of three distinct half-planes that share a line $l$. For any nonzero vector $w \in \R^3$ that is not parallel to $l$, consider a linear surjection $p_w: \R^3 \to \R^2$ with the kernel generated by $w$. There are two possibilities for the map $p_w: \mathcal Y \to \R^2$: 1) generic points in $\R^2$ have preimages of cardinalities one and two, and $p_w$ is onto; or 2) the cardinalities are zero and three ($p_v$ is not onto). We attach symbol $\vec Y$ to the first situation and symbol $\vec W$ to the second one.\footnote{The shapes of the letters are mimicking the desired properties of the half-plane configuration with respect to a horizontal $w$.}  The $\vec Y$-configurations are generated when $w$ and $-w$ point into distinct chambers in which the three half-planes divide $\R^3$.


At each point $x \in s_\dagger(K)^\circ$ of a simple spine $K \subset X$, the linearization of the three surfaces that join at $x$ generates a configuration $\mathcal Y_x$ in the tangent space $T_xX$ of $X$. At each point $x \in Q(K)$, the linearization of the four surfaces that meet at $x$ generates a 2-complex $\mathcal X_x$ which divides  $T_xX$ into four pyramids. We will prefer configurations $\mathcal X_x$ for which $w$ and $-w$ point into distinct pyramids. Such configurations are said to be of the $\vec X$-type. For them, the fibers of $p_w$ will be of cardinality 1, 2, and 3.

\begin{defn} We say that a spine $K \subset X$ is a $\vec Y$-\emph{spine} if there exists a vector field $w$ along $s(K)$ in $X$ which is transversal to each of the surfaces that form $K$ and join along $s(K)$. Moreover, for $x \in s_\dagger(K)^\circ $, the configuration $\mathcal Y_x \subset T_xX$ is of the $\vec Y$-type, and, for any $x \in Q(K)$, the configuration $\mathcal X_x$ is of the $\vec X$-type with respect to $w(x)$.
\end{defn}

For example, consider a union $Y$ of three radii in a disk $D^2$,  $Y$ being symmetric under the rotation $\phi$ on the angle $2\pi/3$. Let $X$ be the mapping torus of $\phi: D^2 \to D^2$ and $K$ the mapping torus of  $\phi: Y \to Y$. Evidently, $K$ is a spine of $X$, but not a $\vec Y$-spine.\smallskip

\begin{lem} Any  gradient spine $K = K(f, v)$  in $X$ is a $\vec Y$-spine.
\end{lem} 

{\it Proof.}\; 
For a $\d_2^+$-generic $v$, the  waterfalls and the ground $\d_1^+X$ are transversal along their intersections. 
At a generic  point  $x \in s_\dagger(K)$ two out of three half-planes form an angle $\pi$. Thus, for an open and dense cone of vectors $w(x) \in T_xX$ the configuration $\mathcal Y_x$ is of the $\vec Y$-type. Consider the two sheets $S_T$ and $S_N$ of $K$ at $x$ that are marked with the $T$ and $N$-markers and a  cone $C_x \subset T_xX$ --- the convex closure of the two half-spaces tangent to $S_T$ and $S_N$ at $x$. Since the unmarked  tangent half-space is never in the interior of  $C_x$,  the cone $C_x$ picks a unique chamber $C$ among the three cambers in which $K$  divides $X$ in the vicinity of $x \in s_\dagger(K)^\circ$. 

At each $Q$-singularity $x$,  $K$ divides the star of $x$ in $X$ into four pyramids.  
For a gradient spine $K$,  the  $NT$ markers pick one of these four  pyramids: its triangular base is built out of three edges that are marked with $TT$, $NN$ and $TN$-markers in Fig. 21. This choice is consistent with the choice of chambers $C$ of the four beams  that merge at $x$. 

We already noticed that for  vectors $w(x) \in C_x^\circ$  the configuration $\mathcal Y_x$ is of the $\vec Y$-type.  By partition of unity and convexity arguments, we conclude that there is a vector field $w$ along the graph $s_\dagger(K)$ such that $w$ belongs to the $NT$-preferred chamber along  
$s_\dagger(K)^\circ$ and to the preferred pyramid $P$ at the points of $Q(K) \subset ss(K)$. At the same time, $-w$ does not belong to $P$. Moreover, $w$ (which has a nonzero projection on the $N$-marked normal to the waterfalls) extends to a field which is transversal to $K$ along $s_\bullet(K)$ as well. \qed

\bigskip

\begin{lem} For an oriented  $X$, the notions of  an oriented branched spine $K \subset X$ and of a $\vec Y$-spine are equivalent, provided $K^\circ$ being orientable.
\end{lem}

{\it Proof.}\; Since $X$ is oriented, an orientation of each component $S$ of $K^\circ$ picks a particular normal $\nu_S$ to $S$ in $X$.  For any $x \in s_\dagger(K)$, an  oriented branching $\a$ on $K$, picks a preferred surface $S_x$ out of the three oriented surfaces that join at $x$. Here is the recipe for choosing $S_x$: if $\g$ is the edge of $s(K)$ through $x$, locally, there are exactly two surfaces, say $S_1$ and $S_2$, such that $\g$ contributes to $\d S_1$ and $\d S_2$ with the same sign. Then $\g$ contributes to $\d S_x$ with the opposite sign. The normal $\nu_{S_x}$ points into a single (among the three) chamber $C_x$ of $X \setminus K$\footnote{The chamber is not necessarily convex in $T_xX$.}. In a sense, $\nu_{S_x}$ is a piece-wise smooth surrogate of the desired field $w(x)$. Note that  $-\nu_{S_x}$ and  $\nu_{S_x}$ point into distinct chambers.

By Proposition 3.1.6, [BP], any oriented branching extends to $Q(K)$, that is, for any $x \in Q(X)$,  the preferred chambers of the four $Y$-beams  that approach $x$ "share" a vector $w_x \in T_xX$ transversal to $K$.  More accurately, infinitesimal parallel shifts of this $w_x$ in the directions of the four edges  $\g_1, \g_2, \g_3, \g_4 \subset s_\dagger(K)$ that share $x$ belong to the preferred chambers $C_{\g_1}, C_{\g_2}, C_{\g_3}, C_{\g_4}$, while the infinitesimal shifts of $-w_x$ do not belong to the preferred chambers. Next, we smoothly interpolate between the parallel shifts of $\{w_x\}_{x \in Q(X)}$ in the vicinity of $Q(X)$ and the fields $\nu_{S_x}$ along the edges of $s_\dagger(K)$.  This interpolation gives the desired field $w$ along $s(K)$. 

Conversely, if $K$ admits a $\vec Y$-structure and $K^\circ$ is orientable, then $K$ is an oriented branched spine. Indeed, the  field $w$, transversal to $K$ along $s(K)$, picks a particular orientation of each component $S$ of $K^\circ$ along its boundary $\d S$. Since $S$ is orientable, it gets a particular global orientation. Thus, $w$ generates an element $\a \in C_2^\Box(K)$ which, because of the very nature of the $\mathcal Y_x$-configurations, has the required property $\d \a \in C_1^\Box(K)$.
\qed


\section{Combinatorial and Gradient Complexities of 3-Manifolds}

Following Matveev [M], the \emph{complexity} $c(K)$ of a simple 2-polyhedron $K$ is defined to be the cardinality of the set $Q(K)$ formed by  points of type (6) in Definition 6.1. This definition of $c(K)$ can be applied to gradient spines as well. \smallskip

\begin{figure}[ht]
\centerline{\includegraphics[height=2in,width=3.0in]{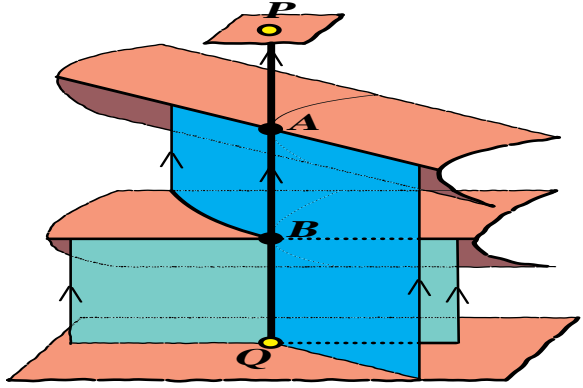}}
\bigskip
\caption{\small{Any $v$-trajectory $QP$, tangent to $\d_1X$ at $A$ and $B$, generates a singularity of type (6) at $Q$.}}
\end{figure}

Here and on, by a \emph{trajectory} of a vector field we mean an integral curve that does not admit a continuation.\smallskip

We notice that the $Q$-singularities of a gradient spine are in 1-to-1 correspondence with segments $[BQ]$ of the $v$-trajectories that are shared by either a pair of waterfalls or, locally, by two branches of the same  waterfall (see Fig. 22). The number of such shared  segments  in a  cascade can be given another, less technical, interpretation. We notice that, for a $\d_2^+$-generic field $v$, a shared segment $[BQ]$ corresponds to a unique pair of distinct points  $A, B \in \d_2^+X$ that are linked by a $v$-trajectory $[PQ]$. In turn, such trajectories are exactly the ones that link distinct points  $A, B \in \d_1X$ and that are \emph{tangent} to $\d_1X$ at $A$ and $B$. Indeed, $\d_2X$ is the locus where $v$ is tangent to $\d_1X$, and points of $\d_2^-X$ do not communicate through the bulk $X$. We call such trajectories $[PQ]$ \emph{double-tangent}.
\smallskip

Generic Morse data $(f, v)$ provides us with an \emph{oriented tangle} $\mathcal T(v) \subset X$ of segments $[AB]$ of double-tangent trajectories, the orientation of  $\mathcal T(v)$ being induced by $v$. 
When $\d_1^+X$ is a disk,  $\mathcal T(v)$ produces a coupling of points in its circular boundary $\d_2X$. 

{\it Question:} \quad 
How does the tangle $\mathcal T(v)$ and the coupling change as $v$ deforms  in the space of nonsingular (gradient-like) fields?\bigskip

Inspired by [M], we propose the following two definitions.

\begin{defn} The complexity $gc(f, v)$ of generic Morse data $(f, v)$ is defined to be the number of double-tangent $v$-trajectories $\g$.\footnote{By a general position argument, we can assume that $\g$ is tangent to $\d_1X$ only at two points.}
\end{defn}

For gradient spines,  \emph{polarities} $\oplus, \ominus$ can be given to the isolated singularities of $Q$-type.  
We have seen that  any $x \in Q(K)$ corresponds to a trajectory $\g$ tangent to $\d_1X$ at a pair of points  $A, B \in \d_2^+X$. Morse data $(f, v)$ help to brake symmetry between $A$ and $B$: indeed, $f(A) > f(B)$. Recall that the preferred orientation of $\d_1^+X$ induces an orientation of $\d_2X$. Consider a vector $v_A$ tangent to $\d_2X$ at $A$ and a vector $v_B$ tangent to $\d_2X$ at $B$, the directions of both vectors agreeing with the orientation of $\d_2X$. The $(-v)$-flow spreads  $v_A$ and $v_B$ and produces an ordered normal frame $(\tilde v_B, \tilde v_A)$ along the trajectory $\g = [P, Q]$. At $x = Q$, the orientation induced by $(\tilde v_B(x), \tilde v_A(x))$ can agree or disagree with the preferred orientation of $\d_1^+X$. In the first case, the polarity of $x$ and $\g$ is defined to be positive ($\oplus$), in the second case, it is negative 
($\ominus$).  We denote by $Q(K)^\oplus$ and  $Q(K)^\ominus$, respectively,  the sets of positively and negatively polarized points in $Q(K)$.  The same polarities $\oplus, \ominus$ can be assigned to the double-tangent trajectories $\g$. Note that reversing the orientation of $\d_1^+X$, reverses the orientation of $\d_2X$, and thus the frame $(v_B, v_A)$ is replaced by $(-v_B, -v_A)$. Therefore, the $(\oplus, \ominus)$ polarity is independent on the choice of  an orientation in $\d_1X$.

\begin{defn} The polarized gradient complexity $gc^\oplus_\ominus(X, f, v)$ of generic Morse data $(f, v)$ is defined to be the difference between the number of positive and negative double-tangent trajectories. 
\end{defn}

The polarized gradient complexity $gc^\oplus_\ominus(X, f, v)$ can be given another interpretation which has the flavor of a ``self-linking number'' for the 1-cycle $\d[K]$. 

By Lemma 6.2, the gradient spine $K = K(f, v)$ has a preferred vector field $w$ along $s(K)$ which gives $K$ its $\vec Y$-structure.  This field  is transversal to $K$ and points inside $X$ along $\d_1^+X$ and into the preferred side of each waterfall from the cascade.  

We view the 1-cycle  $\d[K]$ with support in $s(K)$ as an oriented graph. Employing a preferred field $w$, we push $\d[K]$ a bit in the direction  of the field. Denote $\d[K]_w$ the perturbed graph. By the construction of $w$, the 2-chain $[K]$ and the 1-cycle $\d[K]_w$ are in a general position in $X$, and their intersection points occur only in the vicinity of $Q(K)$, a single intersection per each point of $Q(K)$. Therefore, $gc(f, v) := c(K) = \#(\d[K]_w \cap [K])$.  We claim that an algebraic  count $\d[K]_w \circ [K]$ of intersection points in $\d[K]_w \cap [K]$ also makes sense\footnote{Recall that intersection theories of  spaces with singularities of codimension one fail miserably.}, provided that $X$ is oriented. Let us explain this observation.  Each  vertex of the oriented graph $\d[K]$ whose multiplicity $> 2$ is a $Q$-singularity. It has valency four, a pair of  incoming, and a pair of outgoing edges. Therefore, at each $x \in Q(K)$, there are exactly two oriented resolutions of the graph $\d[K]$ into a pair of arcs. One of the resolutions is the boundary of the resolved surface $res_T(K) \subset X$ (see Fig. 18, the right diagram). The other resolution of $\d[K]$ is denoted by $\overline{res}(\d[K])$. We denote  by  $\overline{res}(\d[K])_w$ its $w$-shift. Consider the algebraic intersection $\overline{res}(\d[K])_w \circ res_T(K)$ of the curve $\overline{res}(\d[K])_w$ with the surface $res_T(K)$. It is easy to see that the points of $\overline{res}(\d[K])_w \circ res_T(K)$ and of $\d[K]_w \cap [K]$ are in 1-to-1 correspondence. I requires more effort to check that the standard orientation assigned to each point of $\overline{res}(\d[K])_w \cap res_T(K)$ is positive  if and only if the corresponding singularity $x \in Q(K)$ has positive polarity $\oplus$.  Since the preferred sides of the surfaces from $K$ that join at $x$ depend only on Morse data (and not on the orientation of $\d_1X$) and so does the choice of $w$, one needs to consider only the interaction of a particular $w$-shift with the eight choices for the the orientations of the two waterfalls and the ground. We leave to the reader the rest of the verification. 

Define the self-linking number $lk_w(\d[K], \d[K])$ by the formula $\overline{res}(\d[K])_w \circ res_T(K)$. By standard homological reasoning, $lk_w(\d[K], \d[K])$ depends only on the graph $\d[K]$ and its preferred framing $w$ (which, in turn, depends only on the marked germ of $s(K)$ in $X$).

Combining Definitions 7.1 and 7.2 with the arguments above, we get

\begin{thm} For generic Morse data $(f, v)$ and its gradient spine $K= K(f, v)$ carrying the preferred framing $w$ along $s(K)$, 
\begin{eqnarray}
gc(f, v) = \#(\d[K]_w \cap [K])
\end{eqnarray}
\begin{eqnarray}
gc^\oplus_\ominus(f, v) = lk_w(\d[K], \d[K])
\end{eqnarray}
Hence, the number of double-tangent $v$-trajectories  is at least  $|lk_w(\d[K], \d[K])|$.
\end{thm}
\smallskip

We can refine  the gradient complexity and view it as an ordered \emph{pair} of nonnegative integers $(\#[Q(K)^\oplus],\, \#[Q(K)^\ominus])$, where $K = K(f, v)$.  
We will see that once a pair $(c_\oplus, c_\ominus)$ has been realized by some Morse data, then $(c_\oplus +1, c_\ominus +1)$ can be realized as well.
\smallskip

Contemplating about formula (7.2), we realize that it makes perfect sense for any simple $\vec Y$-spine $K \subset X$. This leads to
\begin{defn} Let $X$ be a compact 3-manifold  with  boundary. Its $\vec Y$-\emph{complexity}, $c_{\vec Y}(X)$, is the minimal  combinatorial complexity $c(K)$ of $\vec Y$-spines $K$ in $X$.\smallskip

Alternatively, $c_{\vec Y}(X)$ can be defined as $min_{\{K, w\}}\; \#(\d[K]_w \cap [K])$, where $(K, w)$ runs over the set of $\vec Y$- spines $K \subset X$ equipped with their preferred fields $w$.
\end{defn}
\begin{defn} Let $X$ be a compact 3-manifold  with  boundary.  Its 
\emph{gradient complexity}, $gc(X)$, is  the minimum, over all nonsingular $\d_2^+$-generic gradient-like fields $v$, of the number of double-tangent $v$-trajectories.\smallskip

Alternatively, we can define $gc(X)$ as $min_{\{K, w\}}\; \#(\d[K]_w \cap [K])$, where $K = K(f, v)$ runs over generic gradient spines equipped with their preferred fields $w$.
\end{defn}
Evidently, $gc(X) \geq c_{\vec Y}(X) \geq c(X)$. Our Theorem 8.1 implies that actually $gc(X) = c_{\vec Y}(X)$. In general, $gc(X) > c(X)$. For example, for the punctured lens space $L_{3, 1}$, one gets  $c(L_{3,1}) = 0$, while $gc(L_{3,1}) > 0$. \smallskip

Since, for any handlebody $X$ and appropriate Morse data, $\d_2^+X = \emptyset$, it follows that  the gradient complexity of handlebodies equals to zero.  \smallskip
 
In general, computing $c(X)$ is hard; it is much easier to estimate it from above or below. For instance, if $X$ admits a triangulation comprising  $n$ tetrahedrons, then $c(X) \leq n$ ([M], Proposition 2.1.6) . 
\smallskip

For geometrical pieces of the JSJ decomposition (see [J]), we can reduce generic $\vec Y$-spines to special ones without compromising their combinatorial complexity (cf.  Theorem 2.2.4 in [M]). 

\begin{thm} Let $X$ be an  irreducible and boundary irreducible\footnote{Recall the $X$ is \emph{irreducible} if any $S^2 \subset X$ bounds a 3-ball; $X$ is \emph{boundary irreducible}, if for any proper 2-disk $(D, \d D) \subset (X, \d X)$, the curve $\d D$ bounds a disk in $\d X$.} with no essential annuli\footnote{A proper annulus $(A, \d A)  \subset (X, \d X)$ is called \emph{nonessential} if either $A$ is parallel to $\d_1X$ or the core of $A$ is contractible in $X$. The rest of annuli are called \emph{essential}.}  3-fold with $\d_1X \neq \emptyset$, and let $K \subset X$ be 
a simple spine. Then $X$ has another simple spine $K'$ such that
\begin{itemize}
\item $K' \setminus s(K')$ is a collection of 2-disks
\item $s(K')$ is a connected  graph whose vertices are of multiplicity $1$ and $4$
\item $c(K') \leq c(K)$
\item when $K$ is a $\vec Y$-spine, so is $K'$
\item when $K$ is an abstract gradient spine, so is $K'$,
\end{itemize}
If $c(X) > 0$,\footnote{That is, if $X \neq D^3, L_{3, 1} \setminus D^3$.} then $K'$ is a special spine.
\end{thm}

{\it Proof.}\;  
Recall that $X$ is $PL$-homeomorphic to the mapping cylinder of a cellular map $q:  \d_1X \to K$. Denote by $r: X \to K$ the retraction induced by the $q$-mapping cylinder structure. In fact, one can choose $r$ so that, for any $x \in K^\circ = K \setminus s(K)$, $r^{-1}(x)$ is homeomorphic to a closed interval. 

Let $S$ be a typical component  of the set $K^\circ$. When $S$ is a closed surface, $K = S$.

Consider a simple closed path $\g \subset S$ which is not contractible in $S$ and such that the restriction of the normal line bundle $\nu(S, X)|_\g$ is trivial (for gradient spines, $\nu|_\g$ is automatically trivial). 
In fact, such a loop $\g$ always exists, unless $S = D^2,\, S^2,$ or $\R P^2$. 
The last two exceptions can occur only if $s(K) = \emptyset$ and $K = S^2,\, \R P^2$. They correspond to $X$ being a line bundle over $S^2$ or $\R P^2$. The space of a nontrivial bundle over $\R P^2$---the punctured $\R P^3$---is reducible, and the space of a trivial bundle has an essential annulus. The space of a trivial line bundle over $S^2$ is reducible as well. Therefore, unless $S = D^2$, the desired $\g$ always exists.

Denote by $A$ the annulus $r^{-1}(\g)$. Since no essential annuli are permitted, either 1) $A$ is parallel to the boundary $\d_1X$, or 2)  $\g$ is contractible in $X$. 

In the first case, $\g$ must divide $S$. Indeed, if it does not, we can find a loop $\delta \subset S$ which intersects with $A$ at a singleton; this will imply that $A$ is essential.
Moreover, by the definition of $A$ being parallel to $\d_1X$, there is a solid torus $T \subset X$ whose boundary is divided into $A$ and the complementary annulus $A' \subset \d_1X$. As we delete $T$ from $X$, we do not change the topological type of $X$ but  do change $K$ to a new $K' = K \setminus (K \cap T^\circ)$. Again, due to the construction of $A$, deleting $T$ preserves the $r$-induced product structure in $X \setminus K$, and thus $K'$ is a spine. Now $\g \subset s_\bullet(K')$.

\smallskip   

Next, consider the case when $\g$ is nullhomotopic in $X$. Then, by Dehn's Lemma, each of the two loops $\g_1$ and $\g_2$ comprising the boundary of $A$ and residing in $\d_1X$  bounds a disk in $X$. Since  $X$ is boundary irreducible, $\g_1$ bounds a disk $D_1\subset \d_1X$ and $\g_2$ bounds a disk $D_2\subset \d_1X$. Push $D_2$ slightly inside $X$ so that the loop $\d D_2$  slides along $A$ and denote by $D_2'$ the pushed disk. Consider the sphere $\Sigma$ formed by $D_1, D_2'$ and the portion of $A$ between them. Note that $\Sigma \cap K = \g$. Recalling that $X$ is irreducible, we conclude that $\Sigma$ bounds a ball $B \subset X$. Thus $\g$ must divide $S$ into $S'' := B \cap S$ and its complement $S' := (X \setminus B^\circ) \cap S$, and $K$ into $K'' := B \cap K$ and its complement $K' := (X \setminus B^\circ) \cap K$ . Denote by $A'$ the subannulus of $A$ bounded by $\g$ and $\d D_2'$. We notice that $K^\sharp :=  K' \cup A' \cup D_2'$ is a spine of $X \setminus B^\circ$. In fact, we have replaced $S$ by a new component $S^\sharp :=  S'  \cup A' \cup D_2'$ in which $\g$ is contractible.\smallskip

We already noticed that any closed loop $\g$ in $S$ with $\nu|_\g$ being trivial, can be a source for one of the two previous spine modifications. By modifications of the second type, we can insure that each $S$ is incompressible in $X$. Then, by modifications of the first type, we can make sure that all the nontrivial homotopy classes of loops in $S$ will be represented in $\d_\bullet S := s_\bullet(K)\cap S$ and thus will be disjoint from $s_\dagger(K)$.  Therefore, we can assume that no $S$ has a handle and  any nontrivial loop in $S$ is contained in $\d_\bullet S$, unless $S = S^2, \R P^2$---the case that has been ruled out.  Thus,  $S$ must be either 1) a disk, or 2) a disk with a number of holes.  In the second case,  all the boundaries of the holes must be in  $\d_\bullet S$.
\smallskip

Next, we claim that a non-separating simple path $\g$ with   $\d\g \subset \d_\bullet S$ is absent in $S$. Assume to the contrary that such a $\g$ exists.   
Since $r^{-1}(x)$ is a segment for all $x \in S$, $r^{-1}(\g)$ is a disk $D \subset X$ whose boundary $\d D \subset \d_1 X$. Since $X$ is irreducible and boundary irreducible, there exists a 3-ball $B \subset X$ whose boundary $\d B$ comprises $D$ union with another  disk $D' \subset  \d_1X$. The intersection $K \cap B$  bounds $\g$ (on one side), contrary to the hypothesis about $\g$. Therefore, no more than one hole in $S$ is possible: otherwise a non-separating path $\g$ as above exists. If $S$ is an annulus with one of its boundary loops in $s_\bullet(K)$, one can collapse $S$ on the other loop, thus simplifying $K$. Therefore, we managed to construct a spine $K$ with all the components $S$ homeomorphic to a disk. Some of these disks $S$ could have the property $\d_\bullet S \neq \emptyset$, in which case they can be collapsed, further simplifying the spine. 

The moment we arranged for $K^\circ$ to be a disjoint union of 2-disks, the graph $s(K)$ becomes  connected. Suppose to the contrary that $s(K) = s'(K) \coprod s''(K)$, where $s'(K), s''(K) \neq \emptyset$. Consider regular neighborhoods $\mathcal N'$ and $\mathcal N''$ of $s'(K)$ and $s''(K)$ in $K$, respectively. Then the boundaries  $\d \mathcal N'$ and $\d \mathcal N''$ each is a disjoint union of circles. In order to form $K$, we attach disks to $\d \mathcal N'$ and $\d \mathcal N''$; however this will lead to a disconnected $K$, clearly a contradiction. 
The vertices of  $s(K)$ are isolated singularities of $K$ of types (3), (4), and (6) from Fig. 20. Types (3) and (4) have a single edge of $s(K)$  of type (5) and at least one free edge of type (2)\footnote{These free edges are not in $s(K)$ but in the topological boundary of $K$.} that terminates there; type (6) is a four-valent vertex.  After all the collapses, the free edges will disappear and  $s(K)$ will become a regular four-valent graph.  

By [M], Theorem 2.2.4, the only spines $K$ with $K^\circ$ being a union of 2-disks and $Q(K) = \emptyset$ are the spines of $D^3$, $S^3$, $\R P^3$, and $L_{1, 3}$; however, only $D^3$ on this exceptional list satisfies our hypotheses. 
\smallskip

Let us examine how the spine modifications above affected  given
abstract gradient or  $\vec Y$-structures of the original spine. Let $w$ be a preferred vector field  along $s(K)$. Note that all we did amounts to deleting from $X$ a number of relative 3-balls $B$ and solid tori $T$. In all cases, but one, the effect on a spine $K$ was deleting its portion $K \cap B$ or $K \cap T$. Evidently, these operations neither increase the complexity of $K$, nor 
destroy the orientations of $S$ and the $TN$ markers (in the case when $K$ is an abstract gradient spine). In the case of $\vec Y$-spines, the cuts do not 
change the fact that  $w$ along $s(K)$ is transversal to each surface $S$ and that the configurations $\mathcal Y_x\subset T_xX$ are of the $\vec Y$-type with respect to $v$. The only less trivial case occurred when we replaced $K$ with $K^\sharp :=  K' \cup A' \cup D_2'$, but again, the procedure could only eliminate a portion of $s(K)$ which was \emph{disjoint} from the rest and thus did not affected the transversality of $w$ to $S$ or the $TN$ markers. Note that the orientation of $S'$ uniquely spreads  over the cup $A' \cup D_2'$.

The only simplifying move that could harm the abstract gradient structure of a given spine is the collapse of some of the disks or annuli $S$ in the very end of the game. So, if we want to keep the abstract gradient structure of $K$, we should stop there.
\qed
\smallskip

{\it An important warning:}\quad  Note that elementary expansions and collapses of abstract gradient (or oriented branched) spines $K$ are very different from the elementary expansions and collapses of $K$, viewed just as the underlying 2-dimensional complexes. The orientations can prevent us from executing some collapses which non-oriented complexes would support. Also, one need to play close attention to the choice of $NT$ markers, an integral part of the abstract gradient complex structure. 

For example, take a plane on which a circular fence is erected. The fence divides the plane into two domains, the disk and its exterior. If the orientations of the two domains disagree, we cannot collapse the fence.
\bigskip

Our Theorem 8.2 claims that  $gc(X) = c_{\vec Y}(X)$.  
In any case, the obvious inequality $gc(X) \geq c(X)$ can be combined 
with a number of results about $c(X)$ in order to get a lower bound on the number of double-tangent trajectories.

Matveev proved that, up to a homeomorphism,  there are only finitely many compact irreducible and 
boundary irreducible 3-manifolds $X$ that have no essential proper 
annuli and with a bounded 
combinatorial complexity $c(X)$ (see [M], Theorem 2.2.5). We can be a bit more specific:

\begin{thm} Let $X$ be a compact 3-manifold $X$ with a nonempty boundary. For any nonsingular 
$\d_2^+$-generic gradient-like field $v$, the number of  double-tangent $v$-trajectories  is greater or equal to $c(X)$. 

The number of topological types of irreducible and boundary irreducible $X$ with no essential annuli and gradient complexity  $gc(X) = c$ does not exceed 
$\Gamma_4(c)\cdot 12^c$, where $\Gamma_4(c)$ is the number of topological types of four-valent connected graphs with $c$ vertices at most.\footnote{The number  of \emph{labeled} regular four-valent graphs with $c$ vertices is less than $(4c -1)!!$ where $k!!$ denotes the product of all odd numbers that do not exceed $k$}.
In particular,  there are no more than $\Gamma_4(c)\cdot 12^c$ distinct orientable hyperbolic 3-folds $X$ with $gc(X) = c$. 
\end{thm}
{\it Proof.}\; Since we established that gradient spines are special kind of spines and in view of the arguments centered on Fig. 22, the first claim is clear.

In order to prove the second claim, consider a  gradient spine $K(f, v)\subset X$ with $gc(f, v) = gc(X) = c$. Then by Theorem 7.2, we can simplify $K(f, v)$ to a special abstract gradient spine $K$ with $c$ $Q$-singularities at most. Examining the constructions that lead to the proof of Theorem 7.2, we see that deleting 3-disks and solid tori adjacent to $\d_1X$, did not change the combinatorics of the special abstract gradient spine in the vicinity of the remaining $Q$-singularities (as depicted in Fig. 21).  After collapsing all its 2-cells whose boundary touches the topological boundary $s_\bullet(K)$ of $K$,  we could eliminate some of  the $Q$-singularities and transform the $Y$-beams that connect them into $I$-shaped ones.
We conclude that, in the end, $s_\dagger(K)$ must be a connected regular 4-valent graph with $c$ vertices at most.  By Lemma 6.1, the number of such abstract gradient spines has a upper boundary $\Gamma_4(c)\cdot 12^c$.  By Lemma 1.1.15 in [M], any special simple spine $K \subset X$ determines the topological type of its ambient $X$.  

\smallskip

\begin{cor} Let $X$  be an irreducible and  orientable  3-fold $X$  obtained from a closed manifold $Y$ by removing a 3-disk. If  $X$ admits Morse data $(f, v)$ with  no  double-tangent trajectories,  then $X$ is a disk $D^3$.
\end{cor}
{\it Proof.}\; If $gc(X) = 0$, then $c(X) = 0$ as well.
According to  Matveev's  classification list, the only closed irreducible manifolds $Y$ of complexity zero
are $S^3$, $\R P^3$, and the lens space $L_{3,1}$. However, $\R P^3 \setminus D^3$ and $L_{3,1} \setminus D^3$ are reducible. \qed
\smallskip

\bigskip

The $v$-flow through the ``bulk" $X$ and the $v_1$-flow trough its boundary $\d_1X$ are intimately linked. For instance, we get the following proposition:

\begin{cor} Let $(f, v)$ and $(f|_{\d_1X},  v_1)$ be generic Morse data on $X$ and $\d_1X$, respectively, and $v  \neq 0$. If there is no ascending $v_1$-trajectory that connects in $\d_1X$ a point of $\Sigma_1^-$ to a point of $\Sigma_1^+$, then  $v$ has no double-tangent trajectories; in other words, $gc(X) = 0$. 

On  the other hand, $c(X) \neq 0$  implies that, for each nonsingular $f$, there is  an ascending $v_1$-trajectory that connects in $\d_1X$ a point of $\Sigma_1^-$ to a point of $\Sigma_1^+$.
\end{cor}

{\it Proof.}\; By Theorem 4.1, the absence of an ascending $v_1$-trajectory $\gamma$ that connects in $\d_1X$ a point in $\Sigma_1^-$ to a point in $\Sigma_1^+$ implies covexity of the Morse data. On the other hand, if $\d_2^+X = \emptyset$, then no 
$v$-trajectory  links a pair of points $x, y \in \d_1X$ and is tangent to $\d_1X$ at $x$ and $y$. \qed
\bigskip

Given manifolds $X_1$ and $X_2$ with boundary,  two types of connected sum operations are available: $X_1 \# X_2$ and $X_1 \#_\d X_2$. In the first construction, a 1-handle is attached to the interiors of $X_1$ and $X_2$; in the second one, a 1-handle is attached so that the boundary $\d_1 X_1 \coprod \d_1X_2$ is subjected to 1-surgery as well. 

\begin{thm} The Morse complexity is a semi-additive invariant: 
$$gc(X_1 \# X_2) \leq  gc(X_1) + gc(X_2),$$ 
$$gc(X_1 \#_\d X_2) \leq  gc(X_1) + gc(X_2).$$
In particular, attaching a solid handle to $X$ does not increase its gradient complexity.
Also deleting a ball from the interior of $X$ does not change its gradient complexity. Therefore,  only 2-surgery on $X$ has the potential to increase its gradient complexity.
\end{thm}
{\it Proof.}\; The semi-additivity $gc(X' \# X'') \leq gc(X') + gc(X'')$ is easy to validate.  Let $(f', v')$ delivers $gc(X')$ and $(f'', v'')$ delivers $gc(X'')$, where $X'$, $X''$ are 3-manifolds with boundary. Assume that $f''$ attains its minimum at $a \in \d_1^+X''$ and $f'$ attains its maximum at 
$b \in \d_1^-X'$. By adding a positive constant to $f''$, we may assume that $min(f'') > max(f')$. By perturbing $v$ a bit, we can arrange that the cascade in $X''$ has an empty intersection with a small disk $D^2_a \subset \d_1^+X''$, centered at $a$, without changing the original topology of the gradient spine $K''$. Similarly, by perturbing $v'$ if necessary, we can pick a sufficiently small disk $D^2_b\subset \d_1^-X'$ so that the $-v'$ trajectories through $D^2_b$ do not intersect the cascade in $X'$. Then we  attach an one-handle $H \approx D^2 \times [0, 1]$ to $X' \cup X''$ at $a$ and $b$ and extend the Morse data from the top disk  $D^2_a$ and the bottom disk $D^2_b$ inside $H$. The neck  $\g : = \d D^2 \times 1/2$ belongs to the set $\d_2^+(X' \# X'')$ and the annulus $\d D^2 \times [1/2, 1]$ to the set $\d_1^+(X' \# X'')$. Thus, the cylindrical waterfall streaming from  $\g$ does not intersect with the cascade in $X'$ and hits $\d_1^+X'$ without producing new shared trajectories. By the choice of $D^2_a$, the cascade in $X''$ does not fall through $H$. 

Let $T$ denote a solid torus. Since $gc(T) = 0$, we have 
$gc(X \#_\d T) \leq gc(X)$. 
\smallskip

Deleting a ball $B$ from the interior of $X$ does not change its gradient complexity. Indeed, if $(f, v)$  delivers $gc(X)$, then we pick a sufficiently small ball $B$ whose center lies on a $v$-trajectory $\g$ that is not in the cascade. Then $\d B$ is concave with respect to $v$, and the trajectories tangent to $\d B$ are separated from the old cascade. As a result, they do not contribute to the set of double-tangent trajectories of $v$ in $X \setminus B$.  Hence 
$gc(X \setminus B) \leq gc(X)$. On the other hand, if $(f, v)$ delivers $gc(X \setminus B)$, then there is a trajectory $\g$ which connects a point $x \in \d B$ to a point in $y \in \d_1 X$ and is transversal at $x$ and $y$ to the boundary of $X \setminus B$. Drilling a narrow tunnel $W \subset X \setminus B$  centered on $\g$ and with a concave bottleneck  produces a manifold $\hat X$ homeomorphic to $X$ (we need to smoothen $\hat X$ at both ends of the tunnel). This can be done in such a way that $X \setminus B$ and $\hat X$ will share the same cascade. Hence, $gc(X \setminus B) \geq gc(\hat X) = gc(X)$.\smallskip

Therefore,  only 2-surgery on $X$ has the potential to increase its gradient complexity.
\qed
\bigskip

To establish the additivity of $gc(X)$ seems to be much harder. The additivity of  $c(X)$ is a nontrivial fact which relies on Haken's theory of normal surfaces ([M], Theorem 2.2.9). We notice that the defect $gc(X) - c(X)$ is also semi-additive under the connected sum operation.
\bigskip

Each time we have a lower bound on the combinatorial complexity $c(X)$, for any gradient-like nonsingular flow, the number of  double-trajectories must be at least $c(X)$ (Theorem 7.3). Fortunately, [M] contains a complete list of punctured closed irreducible 3-manifolds of combinatorial complexity $\leq 6$. To assemble this list is a labor-intense accomplishment.  
For example, using computations in $[M]$, pp. 77 and 407-408, for \emph{punctured} elliptic manifolds we get: 
$gc(S^3/P_{24}^\circ) \geq 4$,  $gc(S^3/P_{48}^\circ) \geq 5$, $gc(S^3/P_{120}^\circ) \geq 5$. In fact, for $2 \leq n \leq 6$, $gc(S^3/Q_{4n}^\circ) \geq n$. For punctured lense spaces $L_{p, q}^\circ$, one has: $gc(L_{4, 1}^\circ) \geq 1$, 
$gc(L_{5, 2}^\circ) \geq 1$; $gc(L_{5, 1}^\circ) \geq 2$, $gc(L_{7, 2}^\circ) \geq 2$, $gc(L_{8, 3}^\circ) \geq 2$; $gc(L_{6, 1}^\circ) \geq 3$, $gc(L_{9, 2}^\circ) \geq 3$, $gc(L_{10, 3}^\circ) \geq 3$, $gc(L_{11, 3}^\circ) \geq 3$, $gc(L_{12, 5}^\circ) \geq 3$, $gc(L_{13, 5}^\circ) \geq 3$, etc. \smallskip

Reinterpreting Theorem 2.6.2 in [M] (see also [MP]), one gets a lower homological bound on the number of double-tangent gradient-like trajectories in closed manifolds with holes.
	
\begin{cor} Let $X$ be a manifold obtained from a closed irreducible and orientable 3-fold $Y$, different from the lens space $L_{3,1}$, by removing a number of 3-disks.  Then any generic gradient-like nonsingular flow on $X$ has at least $$2\cdot log_5 |Tor(H_1(X; \Z)| + rk (H_1(X; \Z)) - 1$$ double-tangent trajectories. Here $|Tor(H_1(X; \Z)|$ denotes the order of the torsion subgroup $Tor(H_1(X; \Z)) \subset H_1(X; \Z)$.
\end{cor}

For hyperbolic $X$, both $c(X)$ and $gc(X)$ exhibit at least linear growth as  functions of the hyperbolic volume.

\begin{thm} Let $X$ be a compact manifold obtained form a closed hyperbolic 3-manifold $Y$ by removing a number of 3-balls. Let  $V(Y)$ denote the hyperbolic volume of $Y$, and $V_0$  the volume of a regular ideal tetrahedron in the hyperbolic space $\mathbf H^3$. Then, any generic gradient-like nonsingular flow $v$ on $X$ has at least $V(Y)/V_0$ double-tangent trajectories. .

\end{thm}
{\it Proof.}\; The statement follows from [M], Lemma 2.6.7 and Corollary 2.6.8, coupled with the inequality  $c(X) \leq gc(X)$. \qed
\bigskip

Let $G$ be a finitely presented group. The \emph{length of a relation} is the number of generators and their inverses that are present in the relation. The \emph{complexity of a presentation} is defined to be the sum of lengths of all its relations (for example, the complexity of the presentation in Fig. 19 is $3+2+2 = 7$).
\emph{Presentation complexity} $c(G)$ is the minimum complexity among all $G$-presentations.\smallskip

Let us return to the description of the presentation of $\pi_1(X/ \d_1^+X)$ given by a generic gradient cascade and described prior to Theorem 5.3. Its generators are shared segments of trajectories of the cascade together with some free trajectories of waterfalls that have two free trajectories and no shared segments at all. Each waterfall of this kind contributes a single "free" generator and no relations.  We notice that each shared segment in a cascade belongs to three polygons  (see Fig. 10). Therefore, it is present in the relations three times at most.

As a result, the complexity of the presentation of $\pi_1(X/ \d_1^+X)$ induced by a cascade does not exceed three times the number of shared segments of the $v$-trajectories.

This leads to the following analogue of Proposition 2.6.6 from [M] which claims that for any closed irreducible orientable 3-fold $X$, different from $S^3, \R P^3$, and $L_{1, 3}$, \, $c(X) \geq -1 + \frac{1}{3} c(\pi_1(X))$.
\begin{cor} For generic Morse data $(f, v)$, 
$gc(f, v) \geq \frac{1}{3}\, c(\pi_1(X/ \d_1^+X)).$ 

In particular, for any generic Morse data  with a disk-shaped $\d_1^+X$,  \hfill\break
$$gc(f, v) \geq  \frac{1}{3}\, c(\pi_1(X)).$$

Also, for any 3-fold $X$ whose boundary is a union of spheres, 
$$gc(X) \geq \frac{1}{3}\, c(\pi_1(X)).$$
\end{cor}

{\it Proof.}\; In view of the discussion above, we need to clarify only the last statement. Since the image $\pi_1(\d_1^+X) \to \pi_1(\d_1X) \to \pi_1(X)$ is trivial, attaching a cone with the base  $\d_1^+X$ to $X$ only adds new "free" generators to a representation of $\pi_1(X/ \d_1^+X)$ given by an optimal (that is, $gc(X)$-realizing)  cascade.  Hence $\pi_1(X)$ and $\pi_1(X/ \d_1^+X)$ share the same set of relations which implies that $gc(X) > c(\pi_1(X))/3$. \qed
\bigskip

Each self-indexing Morse function on a closed manifold $Y$ provides an upper bound for $gc(X)$, where $X$ is obtained from $Y$ by removing a number of 3-balls.

Let $Y$ be a closed 3-fold with a self-indexing Morse function $h: Y \to \R$ which has a single minimum. Then $h$ and its gradient-like field $v$ give rise to a presentation $\mathcal P_v$ of $\pi_1(Y)$: its generators are in 1-to-1 correspondence with the $h$-critical points $x$ of index one and its relations are in 1-to-1 correspondence the $h$-critical points $y$ of index two. Recall that, for all $x, y$, we have $h(x) = 1$ and $h(y) =2$. Denote by $D_x$ the ascending 2-manifold of $x$ and by $D_y$ the descending  2-manifold of $y$. Each relation $R_y =1$ is obtained by moving along the boundary $\d_y$ of the unstable disk $h^{-1}[1.5,\, 2] \cap D_y$ and recording oriented transversal intersections of the loop $\d_y$ with the disks $\{h^{-1}[1,\, 1.5] \cap D_x\}_x$ (in the surface $h^{-1}(1.5)$) in the order they appear along $\d_y$. This generates a word $R_y$ in the alphabet $\{x, x^{-1}\}$ and thus a presentation $\mathcal P_v$. \smallskip 

Consider the number $c_M(Y) = min_{\{h, v\}}\; length(\mathcal P_v)$, the minimum being taken  over all pairs $(h, v)$ as above. Evidently, $c_M(Y) \geq c(\pi_1(Y))$.

\begin{thm} Let $X$ be a  3-fold obtained from a closed manifold $Y$ by removing a number of 3-balls. Then $\frac{1}{3} \cdot c(\pi_1(X)) \leq  gc(X) \leq 4 \cdot c_M(Y)$.
\end{thm}

{\it Proof. \;}  In view of Corollary 7.4, we need to validate only the second inequality. Consider $(h, v)$ as above. Delete from $Y$ some balls centered on the critical points of $h$. Their size is picked so that, in the vicinity of critical points of indices 1 and 2, the variation of $h$ in each of the balls is less than $\e < 0.5$ . Also we assume that, in the balls, $h$ admits a canonical Morse form. Denote by $X$ the complement to the balls in $Y$ and restrict $(h, v)$ to $X$.  Then each critical point $x$ of index 1  contributes an ``equatorial" annulus $A_x$, and each critical point $y$ of index 2 contributes a pair of ``polar" disks  $B_y^\pm$ to the set $\d_1^+X$. The boundaries  of $A_x$ and $B_y^\pm$ form the set $\d_2^+X$. These observations are based on the quadratic nature of $h$ in the Morse coordinates. Moreover, $\d_1^+X$  consists of  $\{A_x\}_x$ and $\{B_y^\pm\}_y$ together with the sphere $S = h^{-1}(1- \e)$ centered on the point $h^{-1}(0)$ of the absolute minimum $0$. The cascade $\mathcal C(h, v)$ falls from the union of loops $(\cup_x \, \d(A_x))\, \cup \, (\cup_y\, \d(B_y^\pm))$ on the ``ground" $S$. Note that this particular choice of Morse data $(h, v)$ on $X$ has the property $\d_3X =\emptyset$ and thus  is 3-\emph{convex}! 

Denote by $W_x^\pm$ the two waterfalls from $\d A_x$. By choosing $\e$ small enough we insure that each of the two waterfalls $W_y^\pm$ from $\d B_y^\pm$ follow the unstable manifold $D_y$ very closely. In particular, we make sure that, for each $x$, the two loops $\g^\pm := cl\{W_y^\pm \cap [S \cup A_x \cup W_x^\pm ]\}$ 
are ``parallel" and close to the loop $\g:= cl\{D_y \cap [S \cup A_x \cup W_x^\pm] \}$, where $cl$ stands for the closure. Again, by choosing a small $\e$, each intersection point from  $D_y \cap h^{-1}(1 + \e) \cap D_x$  corresponds to a unique point in $\g \cap D_x$, the corresponding points acquiring similar orientations. Thus, each point from $D_y \cap h^{-1}(1 + \e) \cap D_x$ corresponds to  four points in $\g^\pm \cap S \cap W_x^\pm$. These 4-tupples are exactly the $Q$-singularities of the spine in $X$ generated by $(h, v)$. 

Therefore $gc(h, v) = 4\cdot \#\big\{
(\cup_x\, D_x ) \cap h^{-1}(1+\e) \cap (\cup_y\, D_y )\big\}$ ---the length of the representation $\mathcal P_v$. Employing Theorem 7.4, we conclude that the same  holds for any other manifold $\tilde X$ obtained from $Y$ by deleting any positive number of disks.

Note that locally the picture of the spine is symmetric with respect to the planes of 
$D_x$ (and $D_y$), so that the $Q$-singularities of the spine occur in pairs of \emph{opposite} polarities $\oplus, \ominus$. As a result,  $gc^\oplus_\ominus(h, v) = 0$.
\qed
\bigskip

As an example, consider a matrix $A \in SL_2(\Z)$ whose first row is $(a, b)$ and the second one is $(c, d)$. We employ $A$ to form a closed manifold $Y_A$ from two solid tori $\mathbf T_1$ and $\mathbf T_2$. Their boundaries $T_1= S^1 \times S^1$ and $T_2 = S^1 \times S^1$ are glued via $A$.  Here we assume that the second multipliers in the two products are, respectively, the meridians of $\mathbf T_1$ and $\mathbf T_2$. We denote by $X_A$ any manifold obtained from $Y_A$ by removing a number of balls.

\begin{cor} For any $A$ as above, $gc(X_A) \leq 4|b|$. 
\end{cor}

{\it Proof.\;}  The corollary is obtained by constructing a Morse function $h: Y_A \to \R$ such that: 1) on $T_1 = T_2$, $h$ is constant, 2) $h$ has a single minimum and a single critical point of index 1 in $\mathbf T_1$, 3) $h$ has a single maximum and a single critical point of index 2 in $\mathbf T_2$.  Then  apply  to $h$ the construction in  Theorem 7.6. \qed


\section{Which Spines Are Gradient?}

Next, we address  a natural question: Which spines $K$ are produced in the form of a cascade $\d_1^+X \cup \mathcal C(\d_2^+X)$ by appropriate generic Morse data 
$(f, v)$? In other words, \emph{Which spines are gradient?}
\smallskip

It turns out that any $\vec Y$-spine, and thus any \emph{branched} spine, can be approximated by a gradient spine without compromising its combinatorial complexity.

\begin{thm} Let  $K \subset X^\circ$ be a simple 
$\vec Y$-spine. Then, there exist  a nonsingular function 
$f: X \to \R$ and its gradient-like field $v$ so that the pair 
$(f, v)$ produces a gradient spine $\tilde K \subset X$ with the following properties:
\begin{itemize}
\item  $K$ is homeomorphic to 2-complex obtained from  $\tilde K$  by elementary collapses of  some of its $2$-cells 
\item $c(\tilde K) = c(K)$. 
\end{itemize}

Moreover, for any simple spine $K$  and generic Morse data $(\tilde f, \tilde v)$ in $X$, such that $\tilde v$ delivers $K$ as a $\vec Y$-spine, and  for any $\e > 0$, there exist new Morse data $(f, v)$ as above which satisfy two additional properties:
\begin{itemize}
\item $f $ and $\tilde f$ coincide when restricted to $K$
\item the $f$-controlled size of the 2-collapses does not exceed $\e$.\footnote {i.e. the $f$-images of the collapsing $2$-cells are $\epsilon$-small.}
\end{itemize}
\end{thm}
\begin{figure}[ht]
\centerline{\includegraphics[height=2in,width=3in]{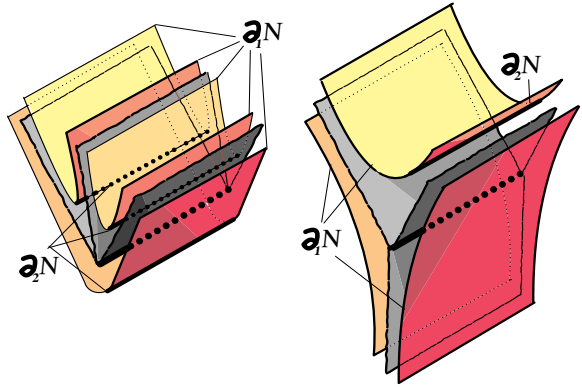}}
\bigskip
\caption{\small{Field $v$ is orthogonal to the plane of drawing $\Pi$.  Let  $p: \d_1N \to \Pi$ be the projection along $-v$. Note the  three bold curves in the left diagram and the single bold curve in the right one, all marked by $\d_2N$. These are the $p$-images of the folds of $p$. On the right,  $p(\d_2N)$ is a simple curve, so no double-tangent $v$-trajectories are present.}} 
\end{figure}
\begin{figure}[ht]
\centerline{\includegraphics[height=1.8in,width=2.6in]{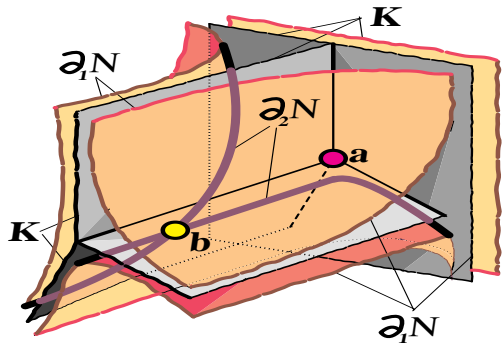}}
\bigskip
\caption{\small{Note the intersection point $b$ of the two bold curves---the $p$-images of the folds of the projection $p: \d_1N \to \Pi$ along $-v$. Point $b$ is the $p$-image of the unique double-tangent $v$-trajectory.}} 
\end{figure}
{\it Proof.}\;  Let $w$ be a field which delivers the $\vec Y$-structure to $K \subset X$.  We start a construction of the appropriate Morse data $(f, v)$ by  extending  $w$  from the graph $s(K)$ to a smooth field $v$ in a open neighborhood of $s(K)$.  Since $s(K)$ is one-dimensional and $w$ is transversal  to $K$ along $s(K)$, 
for a sufficiently small $\epsilon > 0$, the $v$-flow $\phi_t$ will have the property  
$\phi_t(s(K)) \cap \phi_{t'}(s(K)) = \emptyset$ for all $0 \leq t < t' \leq \epsilon$.  Hence, the set $L := \{\phi_t(x)\}_{x \in s(K),\; 0 \leq t \leq \epsilon}$ can be given a product structure $H: L \approx s(K)\times [0,\epsilon]$, $H$ being a diffeomorphism of stratified spaces. We define a function $\hat f: L \to \R$ as the composition of $H$ with the projection on $[0,\epsilon]$. Next, we  extend $\hat f$ to a smooth function $\tilde f: X \to \R$. We can assume that the singularities of $\tilde f$ are not located  on $K$ and thus can be removed by finger moves which originate at $\d_1 X$ and are confined to $X \setminus K$. The resulting nonsingular smooth function $f: X \to \R$ and its gradient field $v$,\, $v|_{s(K)} = w$,  give $K$ its $\vec Y$-spine structure.

Consider a small open regular neighborhood $\mathcal N$ of $K \subset X$. 
The  Regular Neighborhood Theorem (see [Hu], Theorem 2.11, 2.16) implies that $K$ is a spine for $\mathcal N$  and that there is a $PL$-homeomorphism $g: \mathcal N \to X$ which is an identity on $K$. In dimension three, we  can assume that $\d \mathcal N$ is a smooth surface and that $g$ is a diffeomorphism.

Using the smooth product structure in $\mathcal N \setminus K$, we can construct a smooth  function $F: \mathcal N \to [0,1]$ so that: 1) $F^{-1}(0) = K$, 2) $F^{-1}(1) = \d \mathcal N$, 3) $(0, 1]$ being the set of  regular values,
and 4) $K$ being the critical  set. The minus gradient-like flow of $F$ defines a collapse of $\mathcal N$ on $K$. 

Next, we pick $\epsilon > 0$ so small that $\mathcal N_\epsilon : =  F^{-1}([0, \epsilon])$ interacts with the $v$-flow as is depicted in Fig. 23, right diagram,  Fig. 24, and Fig. 27, diagrams 1 and 2. This depiction---a linearization at a point from $s_\dagger(K)$ of the surfaces that form $K$---is based on  $v$ being in general position with respect to $K^\circ$ and giving $K$ its $\vec Y$-structure. 

Let us take a closer look at the interaction of $v$ with the boundary of 
$\mathcal N_\epsilon$: 
\begin{enumerate}
\item in the vicinity of $s_\dagger(K) \setminus Q(K)$,
\item in the vicinity of $Q(K)$,  
\item along the loops of $v$-tangency which are located in $K^\circ$, 
\item in the vicinity of $s_\bullet(K) \setminus s_{\bullet \dagger}(K)$,
\item in the vicinity of $s_{\bullet \dagger}(K)$.
\end{enumerate}
In the first case (see the three-page book in Fig. 23, (2) and its 2D-section in Fig. 27, (2)), the gradient spine $K^\sharp_\epsilon = K(\mathcal N_\epsilon, f, v)$ generated by the Morse data $(f, v)$ on $\mathcal N_\epsilon$ is locally homeomorphic to the given spine $K$. This conclusion depends heavily on the $\vec Y$-nature of the field. Moreover, because of this nature, for a sufficiently small $\epsilon$, the image of the fold is a simple arc---no double-tangent trajectories are present in the vicinity of $s_\dagger(K) \setminus Q(K)$ (see Fig. 23, (2)). 
Compare diagrams 2 and 3 in Fig. 27: diagram 2 reflects the fact that $v$ is of the $\vec Y$-type,  while in diagram 3 the field is of the $\vec W$-type.

In the second case (see Fig. 24 and 27, diagram 1, where $\mathcal N_\epsilon$ is abbreviated to $N$), the surface $\d_1\mathcal N_\epsilon$ and the field $v$ at a 
$Q$-singularity are transversal in two (out of four) chambers-pyramids in which $X$ is divided by $K$. Each of the other two chambers  provides a fold of $\d_2 X$, shown  in Fig. 24 as a bold arc. The  intersection $b$ of the two arcs is the image of  a single double-tangent $v$-trajectory for the $v$-generated gradient spine $K^\sharp_\epsilon$ of $\mathcal N_\epsilon$. 
\begin{figure}[ht]
\centerline{\includegraphics[height=1.5in,width=3in]{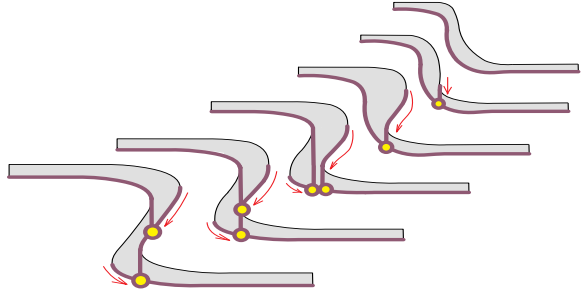}}
\bigskip
\caption{\small{Slicing the gradient spine $K^\sharp_\epsilon$ of $N_\epsilon$ in the vicinity of a  cusp---a point where $v$ is tangent to the arc $L$. The arrows indicate collapses of $K^\sharp_\epsilon$ onto a $2$-complex homeomorphic to the given $K$.}} 
\end{figure}

\begin{figure}[ht]
\centerline{\includegraphics[height=1.8in,width=2.5in]{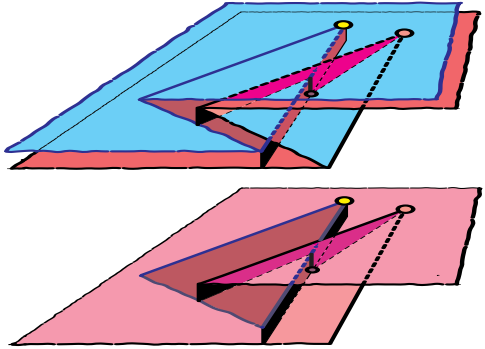}}
\bigskip
\caption{\small{A regular neighborhood of a cusp from $\d_3(K^\circ, v)$ and its gradient spine with a single double-tangent trajectory.}} 
\end{figure}

\begin{figure}[ht]
\centerline{\includegraphics[height=2.5in,width=3in]{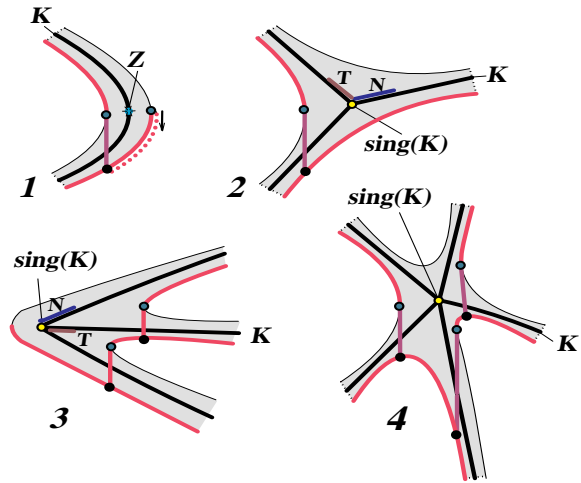}}
\bigskip
\caption{\small{Sections of a spine $K$ and its regular neighborhood by a plane transversal to $s(K)$ 
and containing (vertical) $v$. Diagrams 3 and 4 illustrate complications arising when $v$ is of $\vec W$-type and when $K \subset X$ is not simple.}} 
\end{figure}
\begin{figure}[ht]
\centerline{\includegraphics[height=1.2in,width=3.5in]{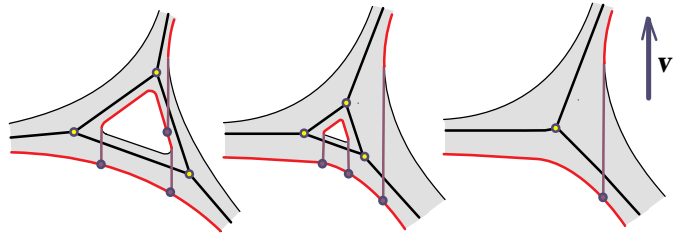}}
\bigskip
\caption{\small{Sections of the gradient spine $K^\sharp_\e \subset N_\e$ by a family of ``parallel" surfaces in vicinity of a $Q$-singularity. The surfaces are invariant under the $v$-flow.}} 
\end{figure}

The third case (see Fig. 27, (1), and Fig. 25), the behavior of $v$ with respect to $K^\circ$ is similar to the behavior of $v$ with respect to $\d_1X$ along $\d_2 X$. Consider a  smooth arc $L \subset K^\circ$ were the field $v$ is tangent to the surface $K^\circ$  and transversal to $L$.  In the vicinity of  a  point $x \in L$, $\mathcal N_\epsilon$ and the flow are represented, up to a diffeomorphism, by the product of first diagram  in Fig. 27 with a segment. Note the portion of the cascade $K^\sharp_\epsilon$ generated in $\mathcal N_\epsilon$ by $v$ and marked with a dotted line. 
Collapsing this portion produces a complex locally homeomorphic to the original $K$. Points $x \in L$ where $v$ is tangent to $L$ require special attension. There $K^\circ$ and $v$ have local geometry similar to the geometry of $\d_1X$ in the vicinity of  the cusp point from $\d_3X$ (see Fig. 25 and 26). 
Fig. 25 demonstrates how $K^\sharp_\epsilon$ collapses onto the given $K$.  We notice that each cusp from $\d_3(K^\circ, v)$ contributes a single $Q$-singularity to $K^\sharp_\epsilon$ (i.e., one double-tangent trajectory to $\mathcal N_\epsilon$), a singularity which has no counterpart in the original $K$! Fortunately, according  to Therem 9.6, we can modify our Morse data $(f, v)$ away from a neighborhood of $s(K)$ (where $v$ is transversal to  $K \subset X$ and gives it its $\vec Y$-structure) so that the modified data will have no cusps in $K^\circ$. Therefore, we get $gc(\mathcal N_\epsilon, f, v) = c(K^\sharp_\epsilon) = c(K)$. \smallskip

Simpler cases (4) and (5) (cf. Fig. 25) are treated similarly. None of them contributes $Q$-singularities to the gradient spine of $\mathcal N_\epsilon$. 
\smallskip

Next, we use the diffeomorphism $g^{-1}_\epsilon: X \to \mathcal N_\epsilon$ to transplant the previously constructed Morse data $(f, v)$ from $\mathcal N_\epsilon$ to $X$. Evidently, the transplanted data and their gradient spine still have all the desired properties.
\smallskip

We leave to the reader to verify that the size of elementary collapses, as measured in terms of the $f$-variation, can be made arbitrary small. The argument uses the uniform continuity of the smooth functions $F$ and $f$ together with the fact that $v$ is in general position with respect to $K$. These imply that, for a sufficiently small $\epsilon' > 0$, there is $\epsilon > 0$ so that in $\mathcal N_\epsilon : = F^{-1}([0, \epsilon])$ the variation of $f$ along each 
$v$-trajectory does not exceed $\epsilon'$ (see Fig. 27, (1)). Moreover, the $f$-controlled size of elementary collapses does not exceed the $f$-controlled size of the trajectories inside $\mathcal N_\epsilon$. This property is preserved under the diffeomorphism $g^{-1}_\epsilon$. \qed
\smallskip

\smallskip


Theorem 8.1 has an important implication:
\begin{thm}   For any compact 3-fold $X$ with a nonempty boundary,  
$c_{\vec Y}(X) = gc(X)$. 
Moreover, for any $X$ as in Theorem 7.2, $c(X) \leq gc(X) \leq 6 \cdot c(X)$.
\end{thm}
{\it Proof.}\; Let $K \subset X$ be a simple spine with $c(X)$ $Q$-singularities. If  such a spine  admits a $\vec Y$-structure, by Theorem 8.1, $c(X) = c_{\vec Y}(X) = gc(X)$, and we are done. 

By [M], Theorem 2.2.4 (cf. Theorem 7.2), any irreducible and boundary irreducible $X$ with no essential annuli has a \emph{special} spine $K$ of complexity $c(X)$. According to [BP], Theorem 3.4.9,   for any special spine $K$, there exists a sequence of not more than $5\cdot c(K)$ oriented Matveev-Piergallini moves (see [M], Figure 1.12) which convert $K$ into a \emph{branched} spine. Furthermore, each move increases the combinatorial complexity of the modified spine  by one. Thus, these moves result in a branched spine of complexity $\leq 6 \cdot c(X)$. Next, by Lemma 6.3,  any branched spine is a $\vec Y$-spine. By Theorem 8.1, it admits an approximation by a gradient spine of  complexity $\leq 6 \cdot c(X)$.  We conjecture that $gc(X) \leq 6 \cdot c(X)$ is valid for any $X$.
\qed
\smallskip



\section{How Deformations of Morse Data Affect the Spine}

Now we investigate the effect of deforming the Morse function $f$ and its gradient-like field $v$ on the gradient spines they generate. We start  with deformations of nonsingular functions that do not introduce singularities in the process. \smallskip

New  \emph{orientations} can be given to  isolated singularities from $\d_3X$. These orientations depend on the Morse data and the preferred orientation of $\d_1X$\footnote{which, when $X$ is orientable, is induced by an orientation of $X$}, i.e. on the structure of an abstract gradient spine inherited by the cascades.
Any point of $\d_3X$ comes equipped with two orientations  marked by ``$\oplus$" and ``$\ominus$": if at $x \in \d_3X$ the the orientation of the arc from $\d_2X$ defined by the tangent vector $v(x)$ agrees with the orientation of that arc induced by the preferred orientation of $\d_1^+X$, then we assign $\oplus$ to $x$; otherwise, we assign $\ominus$. 
The orientations $\oplus, \ominus$ divide $\d_3X$ into two sets $\d_3^\oplus X$ and $\d_3^\ominus X$. As a result, the set  $\d_3X$ is subdivided into four disjoint subsets:  $\d_3^{+\oplus}X$, $\d_3^{-\oplus}X$,  $\d_3^{+\ominus}X$, $\d_3^{-\ominus}X$. Put $\d_3^AX := \d_3^{+\oplus}X \cup \d_3^{-\ominus}X$ and $\d_3^BX := \d_3^{-\oplus}X \cup \d_3^{+\ominus}X$.
\smallskip
\bigskip

As we deform Morse data and watch the transformations of the corresponding spines, we will  pay a close attention to the evolution of signs attached to the $Q$-singularities  and to the points from $\d_3X$.   
\bigskip

\begin{figure}[ht]
\centerline{\includegraphics[height=1.6in,width=3in]{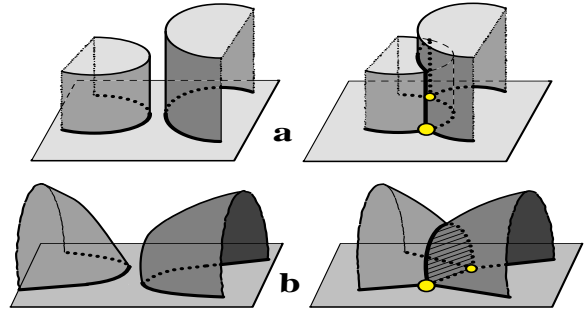}}
\bigskip
\caption{\emph{The $\a$-move shown as a change in the shape of cascades $(a)$ and as a change in the shape of $2$-complexes $(b)$.}}
\end{figure}
\begin{figure}[ht]
\centerline{\includegraphics[height=2in,width=4in]{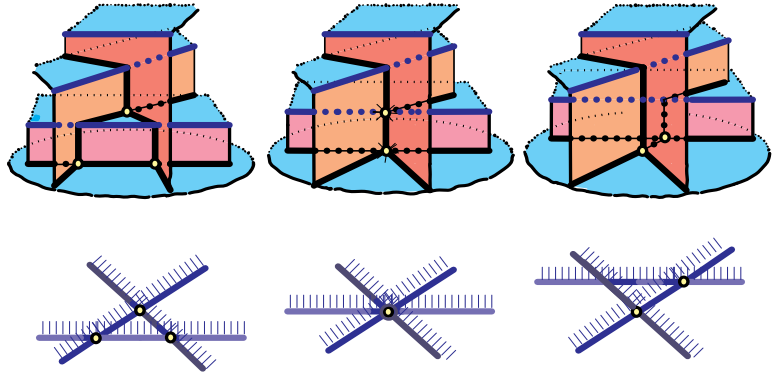}}
\bigskip
\caption{\emph{The $\b^{-1}$-move shown as a change in the shape of cascades and as a $v$-projection of $\d_2^+X$ on $\d_1^+X$.}}
\end{figure}

Two types of elementary transformations of gradient spines are instrumental. The first one is depicted in Fig. 29 as the passage from left to right diagrams and back.  We  call it an $\a$-\emph{move}. An $\a$-move is similar to the second Reidemeister move, where the role of the link diagram is played by the folds $\d_2^+X$. The two $Q$-singularities generated  in an $\a$-move have opposite signs (that is, $\oplus$ and $\ominus$). If we forget the markers (which break the symmetry between the left and right surfaces in Fig. 29, (b)),  $\a$-moves make sense for any unmarked 2-complex $K$. 

The $\b$-move is an analogue of the third Reidemeister move. In a $\b$-move, three branches of $\d_2^+X$ form a triangular configuration, as viewed from the $v$ direction. In the deformation process, the configuration degenerates into one with ``triple intersection", i.e. into a cascade that has a trajectory tangent to $\d_1X$ at \emph{three} distinct points (see Fig. 30). The $\b$-moves come in different flavors depending on the orientations and coorientations of the three folds from $\d_2^+X$ and their ordering by the function $f$ (Fig. 30 shows only one of the possible flavors).   

\begin{thm} Let $\{f_t: X \to \R\}$,\, $t \in [0,1]$,  be a continuous family of smooth non-singular functions and $\{v_t\}$ a corresponding family of gradient-like vector fields. Let $K(v_t)$ denote the spine $\d_1^+X \cup \mathcal C(\d_2^+X)$ generated by $v_t$. Assume that $v_0$ and $v_1$ are $\d_2^+$-generic. Then the $2$-complexes $K(v_0)$ and $K(v_1)$ are linked by a finite sequence of elementary expansions and  collapses of  two-cells combined with a sequence of $\a$- and $\b$-moves and their inverses.  
\smallskip

In the deformation process, the cusps from $\d_3X$ could cancel in pairs as shown in Fig. 32, diagrams 1---6, and their mirror images. The change in topology of $\d_1^+X$ is accompanied by cancelations of pairs from $\d_3X$  with the opposite second  polarity $(\oplus, \ominus)$. The cancellations of pairs with the opposite first polarity $(+, -)$ do not change the topology of $\d_1^+X$. The pairs sharing the same first and second polarities cannot be canceled.  
\end{thm}

\begin{figure}[ht]
\centerline{\includegraphics[height=2.2in,width=3.2in]{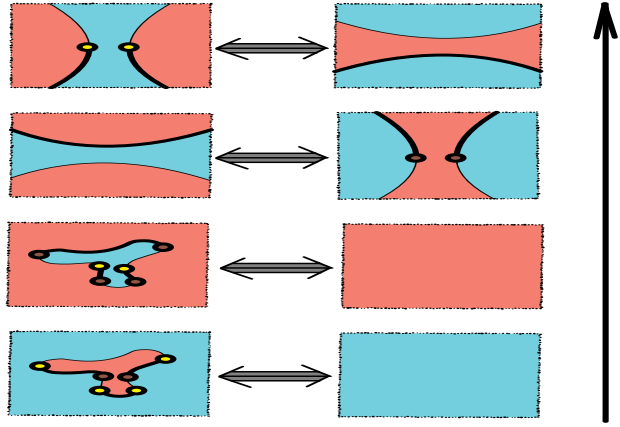}}
\bigskip
\caption{\small{Elementary surgery on $\d_1^+(X)$ and its effect on 
$\d_2^+(X)$ (bold arcs) and $\d_3^+(X)$ (bold dots).}} 
\end{figure}

\begin{figure}[ht]
\centerline{\includegraphics[height=4in,width=3in]{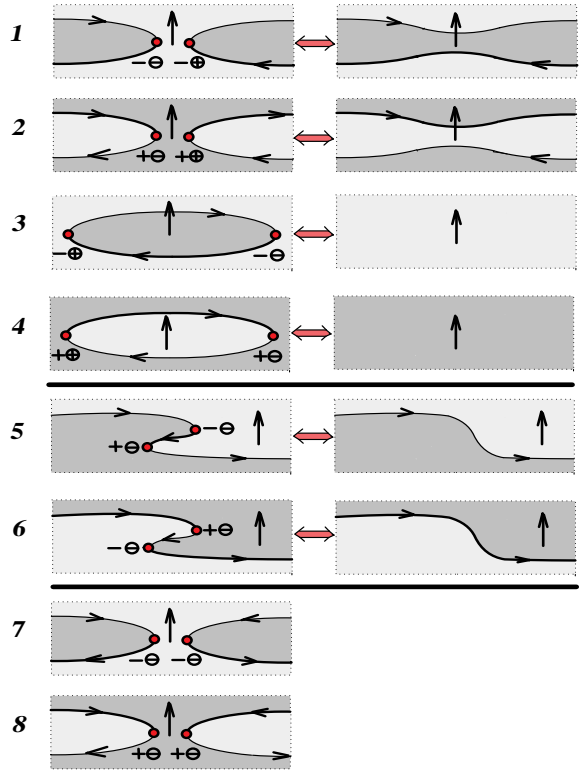}}
\bigskip
\caption{\small{In diagrams 1---4 the topology of $\d_1^+X$ changes, in diagrams 5---6 it does not; diagrams 7---8 show ``impossible cancellations".}}
\end{figure}

We can assume that there are only finitely many $t \in (0, 1)$ for which $v_t$ is not $\d_2^+$-generic. As we deform the Morse data $(f_t, v_t)$,  the topology of the sets $\d_1^+X := \d_1^+(X, v_t)$ is changing by surgery which can be decomposed into a sequence of elementary surgeries depicted in Fig. 31.\footnote{Actually, the third and fourth moves can be decomposed into similar moves applied to a number of convex (round) holes with a single arc for $\d_2^+(X, v_t)$ followed by 1-surgeries  as in the first and second moves.}  All these events take place in the vicinity of a point  $x_\star \in \d_1X$ (where the two singularities from $\d_3X$ merge at a moment $t_\star$)  and propagate inside of $X$ via the waterfalls. In Fig. 31, the  vector field $v_1(x_\star, t_\star)$ is directed upward. 

Diagrams 1 and 2 from Fig. 32 depict 1-surgery, and diagrams 3 and 4 depict 0-surgery and 2-surgery on $\d_1^+X$. In each diagram, the curves from $\d_2X$ are oriented; so the points from $\d_3X$ acquire four flavors: $(+, \oplus), (+, \ominus),  (-, \oplus), (-, \ominus)$. Reversing the orientation of $\d_2X$ flips the polarities  $\oplus \Leftrightarrow \ominus$. Such a flip leads to another four diagrams (not shown in the Fig. 32) which are the mirror images of the diagrams 1---4. They complete the elementary surgery list. We notice that in the process, the first polarities $(+, -)$ of the canceling singularities from $\d_3X$ are the same, and the second polarities 
$(\oplus, \ominus)$ are opposite.  The orientation of $\d_2X$ prevents pairs of the types $(+, \oplus), (+, \oplus)$ or of the types  $(-, \ominus), (-, \ominus)$ from cancelation (see diagrams 7 and 8). 
Again, we think about moves depicted in Fig. 32 as localized events, occurring  at a point of the boundary $\d_1X$,  and whose effect on the gradient spine propagates inside $X$.

Note that each pair of canceling singularities from $\d_3X$ has one point in $\d_3^AX$ and the other in $\d_3^BX$. Hence, the difference $\#(\d_3^AX) - \#(\d_3^BX)$ stays \emph{invariant} under the transformations from Fig. 32. Actually, by Theorem 9.5, this difference is always zero.

Diagrams 5 and 6 from Fig. 32 show another generic mechanism by which the singularities from $\d_3X$ cancel, a mechanism which has no effect on the topology of $\d_1^+X$, but which modifies $\d_2^+X$. In fact, this type of cancellation is generated by the universal dove tail family (2.3). Again, reversing the orientations of the arcs from $\d_2X$ will produce two  more diagrams not shown in  Fig. 32.  This time, the cancelation occurs among the points of opposite first and same second polarities.  \smallskip

Let us first examine the effect of 1-surgery (the first and the second moves) on the shape of the cascades. Fig. 33 shows what happens with cascades as two cusps  merge (as in Fig. 32, diagrams 1 and 2). Viewed in terms of the local projections $p_x : \d_1X \to f^{-1}(f(x))$, this transformation is generated via a generic one-parameter family of mirror-symmetric merging cusps.
\begin{figure}[ht]
\centerline{\includegraphics[height=2.5in,width=3in]{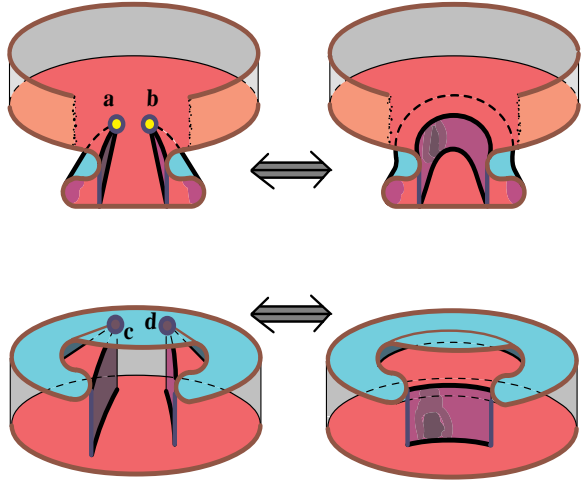}}
\bigskip
\caption{\small{Change of cascades by elementary collapses and expansions of $2$-cells that 
corresponds to $1$-surgery on $\d_1^+(X)$.}} 
\end{figure}
The upper diagram depicts the case when the ``groovings" are located in $\d_1^+(X)$ (Fig. 33, diagram 2), and the lower diagram when they are in $\d_1^-(X)$ (Fig. 33, diagram 1).\footnote{In the first case, the $f$-controlled size of collapses and expansions is small, while in the  second case it can be rather big---a localized change in convexity of the boundary $\d_1X$ causes a distributed change in the shape of cascades.} In the upper  diagram the cusp points $a$ and $b$ of opposite second polarities belong to the set 
$\d_3^+(X)$. Consider the arc-shaped band marked with a dotted line in the upper-right diagram. The band belongs to $\d_1^+(X)$\footnote{most of the band is in the rear and thus invisible.}. In order to get the upper-left diagram, we collapse a 2-cell (a middle rectangle) in that band onto a segment. In the lower diagram, the  cusp points $c$ and $d$  of opposite second polarities belong to the set 
$\d_3^-(X)$ and are assumed to be very close to each other. Therefore, the $(-v)$-trajectories through $c$ and $d$ will hit the same plato (shown  as the lower disk)---the cascades are assumed to be generic. Locally, the spine in the left-lower diagram can be obtained from this plato by two elementary expansions of 2-cells, while the spine in the right-lower diagram by a single elementary expansion.

The spirit of considerations centered on 0-surgery and 2-surgery of $\d_1^+(X)$ is similar. It is illustrated in Fig. 8 which portrays a small indent in a round 3-ball $D^3$ and its effect on the spine. (The indent generates a hole in the set $\d_1^+(D^3)$ of the original round ball.)

Fig. 11 shows the cascades that corresponds to  diagram 6 in Fig. 32. Again, the cascades before and after the dove tail cancelation are linked by elementary $2$-collapses or expansions.

In short, we have seen that changing topology of  the pair $\d_1^+(X) \supset \d_2^+(X)$ via surgery induces  two-dimensional expansions and collapses  of the corresponding gradient spines. \smallskip

Next, let us examine how the spines are affected by deformations of the Morse data $(f, v)$ that do not change the topology of the pair $\d_1^+(X) \supset \d_2^+(X)$. Evidently, under such deformations, the corresponding spine can change only via  changing interactions of waterfalls among themselves and with the ``ground" $\d_1^+(X)$. 

\begin{figure}[ht]
\centerline{\includegraphics[height=2.2in,width=4.2in]{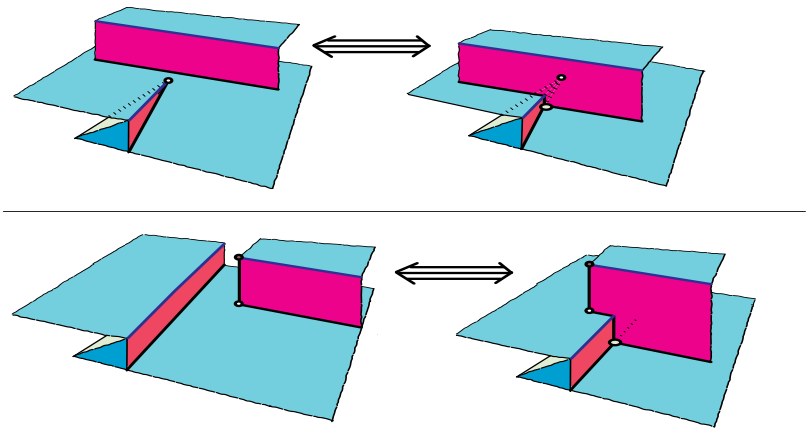}}
\bigskip
\caption{\small{Generating a new $Q$-singularity by deforming the cascade, while keeping the topology of  
$(\d_1^+X, \d_2^+X)$ fixed.}}
\end{figure}

Generically, these changes can be decomposed into sequences of three basic moves, the second and third of which resemble to the second and third Reidemeister moves, respectively.  Let us describe them. 

1) Consider two waterfalls $W_1$  and $W_2$ that fall from two arcs/loops $C_1, C_2 \subset \d_2^+(X)$.  First, $W_1$  and $W_2$ are disjoint in the vicinity of a given trajectory. Then, as the the gradient-like field $v$ changes,  $C_2$ can pierce $W_1$ transversally at a single point $x$ (see Fig. 34) (this is change of the spine  is impossible when $\d_3X = \emptyset$).

2) Alternatively,  $C_2$ can touch $W_1$ and then penetrate it  at a pair of nearby points $y$ and $z$ (see Fig. 29, a). In the vicinity of $x$  on the plato $P_2$ that contains $x$, $W_1 \cap P_2$ is an arc $D$ which is transversal to $C_2$. Consider a small neighborhood  of $D \cup E_x$ in $W_1$, where $E_x$ denotes the  $(-v)$-trajectory trough $x$.  This neighborhood, shaped as a ``fat $\Gamma$" (see Fig. 17), is collapsible. The resulting 2-complex is homeomorphic to the spine formed by the original data that produced  disjoint waterfalls $W_1$ and $W_2$.  In the case of bifurcation shown in Fig. 29 (a) and (b), the topology of the cascade is changing via the $\a$-moves. Recall that the $(\oplus, \ominus)$-polarities of the two $Q$-singularities generated by an $\a$-move are opposite.  

3) Finally, the $\b$-move is generated where the three branches of $\d_2^+X$ form a "triangular" configuration that deforms into a new "triangular" configuration via Morse data that has a trajectory tangent to $\d_1X$ at three distinct points (see Fig. 30). As a result of the $\b$-move, the complexity of the gradient spine jumps  by one; moreover, the polarized complexity also changes by one.
 \qed \smallskip

It seems that it is hard to control the birth and death of $Q$-singularities that accompany  deformations of Morse data. For example, canceling two cusps via a dove tail deformation reduces the number of $Q$-singularities (equivalently, double-tangent trajectories) by one. The situation becomes more manageable in the category of Morse data with no cusps (see Corollary 9.3).
\bigskip

Next, consider the space $\mathcal F(X)$ of smooth nonsingular functions on a compact oriented manifold with boundary. It coincides with the space of all submersions $f : X \to \R$.  Let $\mathcal V(X)$ be the space of smooth nonsingular vector fields on $X$. Philips' remarkable Theorem B  [Ph] claims that, for a fixed metric, the gradient map 
$\nabla: \mathcal F(X) \to \mathcal V(X)$ is a weak homotopy equivalence. When $X$ is an oriented connected 3-fold with boundary, the tangent bundle of $X$ is trivial, and  Philips' Theorem reduces to  the following known proposition:

\begin{thm} Let $X$ be a connected oriented  riemannian 3-fold with boundary. Fix a trivialization of the tangent bundle $\tau X$ of $X$. Then the trivialization-induced normalized gradient map 
$\nabla/ \|\nabla\|: \mathcal F(X) \to Map(X, S^2)$, where  $Map(X, S^2)$ stands for the space of $C^\infty$-maps from $X$ to $S^2$, is is a weak homotopy equivalence.
\end{thm} 

\begin{cor} Let $X$ be as in Theorem 6.2.  Then the trivialization-dependent map 
$\nabla/ \|\nabla\|$ induces a homotopy groups isomorphism
$$\pi_n(\mathcal F(X)) \approx \pi_n(Map(X, S^2))$$
In particular, as sets,  $\pi_0(\mathcal F(X)) \stackrel{h}{\approx} H^2(X; \Z)$ and $\pi_1(\mathcal F(X), f)$ can be identified with  the set $H^1(X; \Z) \oplus H^2(X; \Z)$.
\end{cor}

{\it Proof.}\;
In order to prove the corollary, we need to explain just the validity of the last two isomorphisms. 
With  a trivialization $\beta$ of $\tau X$ being fixed,  any nonsingular $f$ gives rise to a map $\nabla f / \| \nabla f\|: X \to S^2$ whose homotopy class is an element of $\pi^2(X)$. Since $X$ is 3-co-connected\footnote{that is, $H^i(X; G) = 0$ for all $i \geq 3$ and any coefficient group $G$}, by Hopf's Theorem (see Theorem 11.5, [H1]), the natural map $h : \pi^2(X) \to H^2(X; \Z)$ is an isomorphism.  We denote by $h(f)$ the element 
$h(\nabla f / \| \nabla f\|) \in H^2(X; \Z)$\footnote{The choice of $h(f)$ is equivalent to the choice of a $Spin^c$-structure on $X$.}. Thus, $h(f)$ detects the element  
$[f] \in \pi_0(\mathcal F(X))$. 

Any loop in $\g \subset \mathcal F(X)$ can be viewed as function $F: X \times S^1 \to \R$ which is nonsingular when restricted to each fiber $X_\theta := X \times \theta, \; \theta \in S^1$. Hence, $\g$ produces a map 
$G: X \times S^1 \to S^2$, and the homotopy class $[\g] \in \pi_1(\mathcal F(X), f)$ corresponds to an element 
$H(\g) \in  \pi^2(X \times S^1)$ which restricts to $h(f) \in H^2(X; \Z)$. Thus, $\pi_1(\mathcal F(X), f)$ can be identified with the elements of  the set  $\pi^2(X\times S^1)$ that map to $h(f)$ under the natural map $\pi^2(X\times S^1) \to H^2(X\times S^1; \Z) \to H^2(X; \Z)$. Obstructions to linking any pair $F_0, F_1$ of such maps by a  homotopy  lie in  $\oplus_j H^j(X\times S^1; \pi_j(S^2))$. Since $X$ is 3-coconnected, $\pi_2(S^2) \approx \Z \approx \pi_3(S^2)$, and $F_0|_{X\times 0} = f =  F_1|_{X\times 0}$, via the K\"{u}nneth formula, these obstructions lie in $H^1(X; \Z) \oplus H^2(X; \Z)$. With a little more work, one can show that any element of $H^1(X; \Z) \oplus H^2(X; \Z)$ is realizable as an obstruction between $X \times S^1\stackrel{p}{\to} X \stackrel{f}{\to} \R$  and some function $F$ as above. \qed
\bigskip 

Another invariant $e(f) \in H^2(X;\Z)$ of nonsingular Morse functions $f: X \to \R$ on an orientable 3-fold $X$ is available. Its definition is independent on the choice of a trivialization of $\tau X$. Consider an oriented 2-dimensional vector bundle $\eta$ on $X$ formed by the planes tangent to the constant level surfaces of $f$.  The orientation of $\eta$ is induced by the orientation of $X$ and by $\nabla f$. Note that $\eta\oplus \R$ is isomorphic to the tangent bundle of $X$ and therefore is trivial. Let 
$e(f)\in H^2(X;\Z)$ be the Euler class of $\eta$. The element $e(f)$ is invariant under homotopies of $f$ through nonsingular functions.

Since the  bundle  $\eta \oplus \R$ is trivial, for each choice of the trivialization $\beta$, the isomorphism class of $\eta$ is described by the homotopy class of the appropriate map $E(f): X \to SO(3)/SO(2) \approx S^2$. In turn, the class of $E$ is detected by the Euler class of the bundle $\eta$.  The relation between $E(f)$ and $\nabla f/ \|\nabla f\|$ is well-known. It is described by the lemma below whose proof is left to  the reader.
\begin{lem} The homotopy class of the map $E(f): X \to S^2$ is twice  the homotopy class of the map $h(f) : X \to S^2$. Thus, $e(f) = 2h(f)$.
\end{lem}

\smallskip

Combining Theorem 9.1 with Corollary 9.1 we get one of our main results:

\begin{thm} Let $(f, v)$ and $(\tilde f, \tilde v)$ be a pair of  generic  Morse data, where the fields $v, \tilde v$ are nonsingular, and such that $h(f) = h(\tilde f)$\footnote{Both fields define equivalent $Spin^c$-structures on $X$.}. Then there exists a sequence of elementary $2$-expansions, $2$-collapses,  $\a,\, \a^{-1}, \b,\, \b^{-1}$-moves which transforms the gradient spine $K(v)$ into the gradient spine $K(\tilde v)$. Therefore, if $H^2(X; \Z) = 0$, then any two gradient spines of $X$ are linked by a sequence of elementary $2$-expansions, $2$-collapses, intermingled with $\a$ and $\b$-moves and their inverses.
\end{thm}

It remains to sort out what happens to a gradient spine $K(v)$ when the value of the invariant 
$h(f) \in H^2(X; \Z)$ jumps. In fact, due to Theorem 9.3 , it suffices to analyze how the spine changes as a result of  critical  points of $\{f_t\}, 0 < t < 1,$ ``traveling through $X$" along arcs that represent a generator of $H_1(X, \d X; \Z)$ (see Fig. 35). Here we assume that $f_0$ and $f_1$ are nonsingular and that the traveling critical points are of the Morse type. \smallskip

\begin{figure}[ht]
\centerline{\includegraphics[height=3in,width=3in]{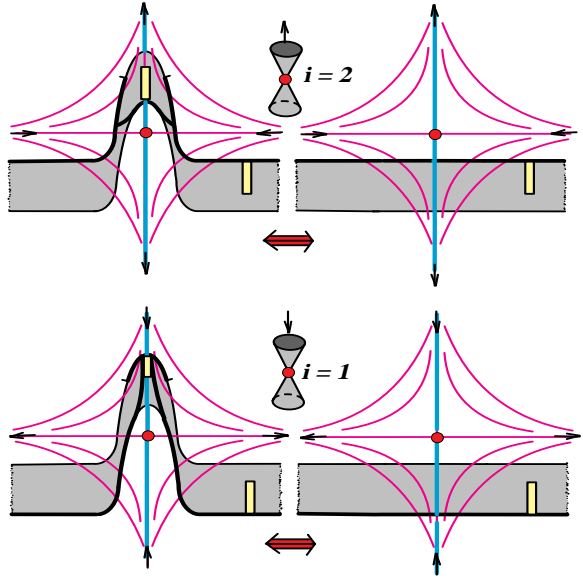}}
\bigskip
\caption{\small{Changing $h(f)$ by the dual of the 1-cycle $(-1)^i[J]$ via an isotopy. The spine changes by mushroom flips. Note the change in the $T$-markers (shown as small vertical rectangles).}} 
\end{figure}

The class $h(f) \in H^2(X; Z)$ is Poincar\`{e}-dual to the oriented 1-dimensional locus $J \subset X$ where $v$ has fixed, up to proportionality, coordinates in the basis that trivializes $\tau X$.  $J$ is as a union of oriented loops and arcs in $X$ with end in $\d X$. \smallskip

Suppose we have a homotopy $\{f_t\}_ {t\in [0,1]}$, of Morse functions on $X$ so that the singular set $\Sigma\{f_t\} = \{(x,t)\in X\times [0,1]\ |\ x\text{ is a (Morse) critical point of }f_t\}$ consists of a collection 
$\tilde J$ of arcs  in $X\times (0,1)$ with endpoint in $\d X\times (0,1)$. Let $J$ be the image of $\tilde J$ under the projection onto $X$. By transversality, we assume that $J\subset X$ is a union of disjointly embedded arcs $\{J_\alpha\}$ with end points in $\d X$. These arcs are oriented according to the direction in which the $t$ parameter is increasing.

\begin{lem}\label{lem e(f) changes by 2TP on J} 
The cohomology class $h(f_1)-h(f_0)$ is dual to the relative 1-cycle 
$\sum_\alpha (-1)^{i_\alpha}[J_\alpha]$ where $i_\alpha$ is the Morse index of the critical point that traces the arc $J_\alpha$. When all indices are even, $h(f_1)-h(f_0)= g^\ast[S^2]$, where $g:X\to S^2$ is the map defined by the Thom-Pontjagin construction on this union of arcs $J (=g^{-1}(\ast))$.
\end{lem}

\begin{proof} 
Fig. 35 reflects the spirit of the argument. It depicts the case of a critical point an index two tracing an oriented arc $J$ in $X$. In fact, the figure shows an isotopy  of the manifold $X$ against a background of a ``stationary" Morse function $f$ and a trivialization $\beta$, both defined on a larger manifold $\hat X \supset X$ ($\hat X $  is obtained from  $X$ by two elementary expansions using three dimensional cells). 

The isotopy of $X$ is supported in  a cylinder  
$C:= D^2 \times \hat J \subset \hat X$. When we isotop $X$ inside of $\hat X$, the original trivialization  changes by a homotopy which is constant outside $C$. As long as $v \neq 0$, this deformation does not change the  class of the map $\nabla f/ \|\nabla f\|$ in $\pi^2(X)$.  

Examining the locus where $v$ is vertical and points up (see the upper diagram in Fig. 35), we conclude that the Poincar\'{e} dual of the variation of  $h(f)$ can be represented by the oriented arc $J$ ---the bold arrow in the figure. The case of index one critical point is similar: the gradient $\nabla f$ will flip its direction (in comparison to the one shown in Fig. 35) causing the change in the orientation of $J$. In the case of the Morse index $i$, the variation of the Poincar\'{e} dual of $h(f)$ is given by the formula $(-1)^i [J]$. Note that when the critical point $x$ of $f$ is inside $X$, the map $\nabla f/ \|\nabla f\|$ is only well-defined in $X \setminus x$. \smallskip

Now we are in position to prove  Lemma 9.2. Put $Y = X\times [0,1]$.
Consider a  bundle $\tau$ tangent to the fibers of the obvious projection $Y \to X$ and its trivialization $\beta$. The gradient-like fields $v_t$ define a section $w$  of $\tau$ that vanishes on $\tilde J =  \cup_t \; \Sigma_t \times t$, and thus a map $\Phi: Y \setminus \tilde J \to S^2$ is well-defined. Denote my $U$ a small regular neighborhood of $\tilde J$ such that all the fibers $U_x$ of the projection $U \subset  Y  \to X$  are homeomorphic to  two-disks (the disks get  truncated as they approach $\d_1X$) in which the Morse function $f_t$ acquires its ``almost canonical" form $$a_1(t) x_1(t)^2 + a_2(t) x_2(t)^2 +a_3(t) x_3(t)^2,$$ with $|a_1(t)| < |a_2(t)| < |a_3(t)|$. Along $\tilde J$, the Morse coordinates $(x_1(t), x_2(t), x_3(t))$  could disagree with the trivialization $\b$.  However, this disagreement happens along a bunch of arcs  $\tilde J_\a \subset X\times [0, 1]$ that have \emph{disjoint} projections $J_\a$ in $X$. Because each arc is contractible, we can homotop $\b$ in the vicinity of each $J_\a$ so that the new trivialization will be adjusted to the Morse coordinates. (Recall, that a homotopy of $\b$ does not change the invariants $h(f_t)$.) Now, by general position, we can assume that 1) $\Phi$ is transversal to a base point $\ast \in S^2$, 2) $\Phi^{-1}(\ast) \cap (X\times 0) = (\nabla f_0/|\nabla f_0|)^{-1}(\ast)$ and 
$\Phi^{-1}(\ast) \cap (X\times 1) = (\nabla f_1/|\nabla f_1|)^{-1}(\ast)$, 3) $\Phi^{-1}(\ast) \cap U$ is given by $x_2(t) = 0 = x_3(t),\, x_1(t) > 0$. Let $Z$ be the surface $\Phi^{-1}(\ast) \cap (Y \setminus U)$ equipped with an orientation induced by the orientations of $S^2$ and $Y$.
The boundary of the 2-chain $Z$ is the 1-cycle which satisfies the equation 
$$\d Z = (\nabla f_1/|\nabla f_1|)^{-1}(\ast)  - (\nabla f_0/|\nabla f_0|)^{-1}(\ast) +  \sum_\alpha (-1)^{i_\a}  {\tilde J'}_\a$$where ${\tilde J'}_\a \subset \d U$ is parallel to $\tilde J_\a$. Hence, $h(f_1)-h(f_0)$ is dual to the relative 1-cycle $\sum_\alpha (-1)^{i_\alpha}[J_\alpha]$. 
\end{proof}
\bigskip

We will see that the effect of a Morse singularity passing through $X$ on the gradient spine $K = K(v)$ can be described as the following procedure: cut a two-dimensional disk $D$ out of $K^\circ$, stretch $D$ uniformly into a bigger disk $\hat D \supset D$, \emph{flip the orientation} of $\hat D$ and paste $\hat D$ back to $K \setminus D$ so that the boundary of $D\subset K$ is identified with the circle that bounds $D$ inside of $\hat D$. The result of this pasting, $K(\tilde v)$ looks like a mushroom with the ``head" $\hat D$. It is equipped with  the new $TN$-markers along $\d D$ that reflect the new orientations. We call such spine changes   \emph{mushroom  flips}. They are manifestations of jumps of the $Spin^c$-structure on $X$.

\begin{thm} Let $(f, v)$ and $(\tilde f, \tilde v)$ be two  generic pairs of Morse data on a compact oriented 3-fold $X$ with boundary, the fields $v, \tilde v$ being nonsingular. Then there exists a sequence of  $2$-expansions, $2$-collapses,  $\a,\, \a^{-1}, \b,\, \b^{-1}$-moves, and mushroom flips  which transform the gradient spine $K(v)$ into the gradient spine $K(\tilde v)$. 
\end{thm}

{\it Proof.}\; Combining Theorem 9.3 with Lemma 9.2 reduces the problem  to understanding the changes in the shape of gradient spine that are affected by of  critical points of an appropriate index  traversing $X$ along  oriented arcs $\{J_\a\}$ representing a given generator $\g$ in $H^1(X, \d X; \Z)$.  The arks representing $\g$ can be chosen in general position with respect to a given gradient spine $K$.  That is, they are transversal to $K^\circ$ and have an empty intersection with  $s(K)$. Furthermore, if an oriented arc $J_\a$ hits the cascade $C \subset  K$ transversally at a point $x$, then we can replace it with two oriented arcs $J'_\a$ and $J''_\a$ such that  
$J'_\a \cup J''_\a$ is homologous to $J_a$ and $(J'_\a \cup J''_\a) \cap C = (J_a \cap C) \setminus x$. To construct the new arcs, use the down trajectory $\g_x$ through $x$ and a band $B$ with the core $\g_x$, $B$ being  transversal to the waterfall $W$ containing $x$. 

By somewhat similar construction, we can find a representative of $\g$ so that that each arc $J_\a$ hits $\d_1^+X$ at a single point $x$. For example, if we have an arc $J_\a$ with two ends, $x, y \in \d_1^+X$, then we pick a path $\delta$ connecting a point $z \in Int(J_\a)$ to a point $w \in \d_1^-X$, form a narrow band $B$ to with the core $\delta$, and use two arcs in $\d B$ to replace $J_a$ with two arcs $J'_\a$ and $J''_\a$, each having the desired property. Therefore, we can assume  that each $J_\a$ either  is oriented in accordance with the vector $v(x)$, or opposite to it. In the first case, we send a critical point of an even index to trace $J_\a$, in the second case, we  send a critical point of an odd index. 

Thus we need to describe only the case where an arc $J_\a$  hits $\d_1^+X \subset K$ transversally away from the cascade $C \subset K$. Each intersection of this type will be responsible for one mushroom flip as shown in Fig. 35. 
\qed

\bigskip


\begin{thm} Let $X$ be an an oriented compact 3-fold. Any generic Morse data $(f, v)$ can be deformed into new data $(\tilde f, \tilde v)$, so that the new $\d_3X(\tilde f, \tilde v) = \emptyset$, and $gc(\tilde f, \tilde v) = gc(f, v)$. In particular, there are no topological obstructions to the $3$-convexity.

\end{thm}

Theorem 9.5  employes Theorem 9.6 below which is similar in spirit to some results  from [E1], [E2] concerned with folding maps of surfaces. 
    
\begin{thm}
Let $S \subset X$ be a connected compact oriented and two-sided surface regularly embedded in the ambient  $3$-fold $X$, and let $(f, v)$ be non-singular Morse data such that $v$ is transversal to $S$ along its boundary $\d S$. Then there exists a deformation of $(f, v)$ which is fixed in the vicinity of $\d S$ and such that the new generic data do not have cusps in $S$.
\end{thm}
Since the surface $S$ has a preferred side in $X$, it can be divided with the help of  $v$ into two domains $S^+$  and $S^-$ which share a common boundary $L$.  In $S^+$ the field points into the preferred side of $S$ and is tangent to $S$ along the locus $L$. Since $v$ is transversal to $S$ along $\d S$, for a generic $v$, $L$ is a collection of loops. The preferred orientation of $S^+$  induces a particular orientation on $L$. The 1-submanifold  $L$ is divided by the cusp locus $C$ into portions $L^+$ and $L^-$. Along $L^+$, $v$ points inside $S^+$.  As before, the points from $C$ acquire four flavors: $(+, \oplus ), (-, \ominus ), (-, \oplus ),$ and $(+, \ominus )$. The first $\{+, -\}$ polarity reflects the fact that $v$ points inside or outside of  $L^+$. The second polarity $\{\oplus, \ominus\}$ tells us whether the field  agrees or disagrees with the orientation of $L$. We denote by $C^A$ the points of the first two flavors $(+, \oplus ), (-, \ominus )$ and by $C^B$ of the last two flavors $(-, \oplus ), (+, \ominus )$.\smallskip

We divide the proof of Theorem 9.6 in three lemmas.

\begin{lem}\label{lemma A} Under the hypotheses of Theorem 9.6, $\#(C^A) = \#(C^B)$. In particular, for  generic Morse data $(f, v)$, 
\begin{eqnarray}
\#(\d_3^{+\oplus}X) - \#(\d_3^{+\ominus}X) + \#(\d_3^{-\ominus}X) - \#(\d_3^{-\oplus}X) = 0.
\end{eqnarray}

\end{lem}
{\it Proof.}\; We already remarked that $L$ must be a disjoint union of simple loops. Each loop $L_i$ from $L$ either has no cusps, or  the arcs from $L^+$ and $L^-$ alternate along $L_i$. For an arc from $L^+$, we examine the four possible flavors attached to its end points $a$ and $b$ and see that one  the two polarities of $a$ and $b$ must be different. Therefore, if $a \in C^A$, then $b \in C^B$. As a result,  $\#(C^A) = \#(C^B)$. In the case of $S = \d_1X$, this leads to (9.1).
\qed
\begin{lem}\label{lemma B} Let $a, b \in C$ have flavors either $1)$ $(+, \oplus)$ and $(+, \ominus)$ or  $2)$ $(-, \oplus)$ and $(-, \ominus)$, respectively. In the first case, assume that  $a$ and $b$  can be connected by simple path $\gamma \subset S^+$,  in the second case, assume that $\g \subset S^-$.  Then $a$ and $b$ can be cancelled via a  deformation of $(f, v)$ as in Fig. 32, diagrams 1---2. The deformation is an identity away from  a regular neighborhood of $\g$.
\end{lem}
{\it Proof.}\; The two cusps are mirror images of each other in the sense that there are ambient coordinates $x, y, z$ in the vicinity of $\g \subset X$ so that the Morse function is $f(x, y, z)$ has a form $z+g(x,y)$ and the surface $S$ is given by $y=z^3+(x-a)z$ at the cusp $(a, 0, 0)$ and by $y=z^3-(x-b)z$ at the cusp $(b, 0, 0)$. In  local coordinates, the path $\g$  can be given by $y=z=0$ with $x$ ranging from $a$ to $b$. The surface is transverse to  $\nabla f$ along $Int(\g)$,  and the space of all germs of two-sided oriented surfaces with this property has the homotopy type of $S^0$ (there is no topological obstruction to making the surface standard along $\g$). So, we can cancel the two cusps by embedding a standard model and the deformation as in Fig. 33.\qed
\begin{lem}\label{lemma C} If two consecutive cusps $a$ and $b$ along a loop $L$ are as in Lemma \ref{lemma B},
then $a, b$ can be canceled so that the new Morse data has a tangency locus $L' \subset S$ which is the result of 0-surgery on $L$. The deformation of the Morse data has a support in an arbitrary small neighborhood of the arc 
$[a, b] \subset L$.
\end{lem}

{\it Proof.\;}
There is a path $\gamma$ along the surface that links $a$ and $b$ inside of $S^+$ in the case of $(+,\ominus)$ and $(+,\oplus)$ or $S^-$ in the other case. This path can be found in an arbitrarily small neighborhood of the arc $[a,b] \subset L$. By Lemma \ref{lemma B}, the cancellation is possible.\qed
\begin{lem} Let $S \subset X$ be as in Theorem 9.6. Then in the vicinity of any point $x \in L^+$ of the tangency locus $L \subset S$ there is a deformation of the Morse data so that two new consecutive cusps of types $(+,\oplus)$ and $(-,\oplus)$ or of types $(+,\ominus)$ and $(-,\ominus)$ are introduced in $L$ via the dove tail deformation. The deformation is an identity away from $x$.
\end{lem}
{\it Proof.}\; See Fig. 32, diagrams 5---6, and Fig. 11 depicting the dove tail surface. \qed
\smallskip

{\it Proof of Theorems 9.5 and 9.6.}\; 
Take any two consecutive cusps $a,b$ along $L$ so that the connecting arc $[a,b]$ lies in $L^+$. Say $a$ is of type $(+,\oplus)$. Then $b$ is either of type $(+,\ominus)$ or $(-,\oplus)$. In the first case, according to Lemma \ref{lemma B}, $a$ and $b$ can be cancelled. In the second case, by Lemma \ref{lemma C}, we can introduce two cusps $c, d\in [a,b]$ of type $(+,\ominus)$ and $(-,\ominus)$, respectively. Then, by Lemma \ref{lemma B}, $a$ and $c$ can be cancelled, as well as $d$ and $b$. The same argument works in the case when $a$ is of any other type than $(+,\oplus)$. Continuing in this way, all cusps can be eventually eliminated, which completes the proof of Theorem 9.6. \smallskip

In order to prove Theorem 9.5, we need to examine carefully the previous argument. Consider all the arcs (but not loops) $[a, b] \subset \d_2^-X$ with the different first polarities of $a$ and $b$ (then, by an argument in Lemma 9.3, the second polarities of $a, b$ agree). For every such arc $[a, b]$,  we introduce a pair of cusps $c, d \in [a, b]$ as above. Because $[a, b] \subset \d_2^-X$, the new waterfall of $[c, d]\subset \d_2^+X$ is "protected" by $\d_1^-X$ and isolated from the rest of waterfalls; as a result, the original gradient complexity is not affected by the introduction of $c, d$ (contrast this with Fig. 11 where complexity increases by 1). Introducing $c, d$ also has no affect on the degree $deg(h) = \#(\d_3^+X) - \#(\d_3^-X)$ of the map $h: \d_2X \to S^1$ from Lemma 2.1. Thus, we can assume that every arc from $\d_2^-X$ has cusps of opposite second polarity.

Now,  for each arc $[a, b] \subset \d_2^-X$, we  pick a path $\g \subset \d_1^\pm X$ which connects two canceling cusps $a, b$ (as in Lemma 9.5) and resides in the vicinity of  $[a, b]$.  The new portion of the waterfall that is generated after the cusps' cancelation is also localized in the  vicinity of $[a, b]$. Hence it is "protected" by $\d_1^-X$ and separated from the old waterfalls that existed before the cancellation (see Fig. 33, upper diagram). Therefore, canceling $a, b$ via such $\g$, again,  does not  change the gradient complexity of the original Morse data. Since the first polarities of $a$ and $b$ agree, each cancelation does  change the degree  of the map $h$  by one. So, we will need at least $deg(h)$ cancellations to get to the Morse data with $\d_3X = \emptyset$.
\qed

\begin{cor}
In terms of the polarized  cusps, the degree of the map $h: \d_2X \to S^1$ can be expressed as follows:
\begin{eqnarray}
 deg(h: \d_2X \to S^1) = [\#(\d_3^{+\oplus}X) - \#(\d_3^{-\oplus}X)]  \\ \nonumber = [\#(\d_3^{+\ominus}X) - \#(\d_3^{-\ominus}X)]  = 
 \chi(X) - 2\chi(\d_1^+X).
\end{eqnarray}
Hence, for the Morse data with fixed values of $\chi(\d_1^+X)$\footnote{For example, for Morse data with $\d_1^+X$ being a disk} 
the number $\#(\d_3^{+\oplus}X) - \#(\d_3^{-\oplus}X) = \#(\d_3^{+\ominus}X) - \#(\d_3^{-\ominus}X)$ is a topological invariant.
\end{cor}

{\it Proof.}\;  By Lemma  9.3, $\#(\d_3^{+\oplus}X) - \#(\d_3^{+\ominus}X) + \#(\d_3^{-\ominus}X) - \#(\d_3^{-\oplus}X) = 0$.  Combined with Lemma 2.1, this leads to the formula for the degree of  $h: \d_2X \to S^1$ claimed in the corollary.  \qed
\bigskip

Now, consider only  nonsingular Morse data $(f, v)$ such that $v$ is transversal to the \emph{submanifold} $\d_2 X \subset \d_1 X$.  Denote by $\mathcal W(X)$ their space. In particular, for elements $(f, v) \in \mathcal W(X)$, 
$\d_3X = \emptyset$.

Similar spaces of smooth maps with folds only from a manifold $M^n$ to a manifold $N^n$  have been studied in great generality in [E1], [E2]. In our context, we are lacking a nice target space $N^2$. Its role is played by the space $X/\sim_{v}$ of $v$-trajectories, a space  which has a structure of a cellular 2-complex and is singular in general.\smallskip

Recall that, according to Theorem 9.5, 
\begin{eqnarray}
gc(X) = min_{\{(f, v) \in \mathcal W(X)\}} \; gc(f, v).
\end{eqnarray}
\smallskip

Let $\mathcal W{\Huge{\ast}}(X) \subset \mathcal W(X)$ be an open and dense subspace of Morse data $(f, v)$ for which no $v$-trajectories, tangent to $\d_1X$ at three distinct points, exist.  Note that no $\b$-move is possible within $\mathcal W{\Huge{\ast}}(X)$ (see Fig. 30). The codimension one walls  $\mathcal W(X) \setminus \mathcal W{\Huge{\ast}}(X)$ can be cooriented by the the following rule: a path $\g \subset \mathcal W(X)$, which represents a $\b$-move, defines a positive coorientation when (as a result of the $\b$-move)  the difference between the numbers of $\oplus$ and $\ominus$ double-tangent trajectories \footnote{equivalently, $\#[Q(K)^\oplus] - \#[Q(K)^\ominus]$} jumps by $+1$. This coorientation is similar in spirit but different from the one used by V. Arnold in his studies of the spaces of immersions of plane curves [A].

We have shown that any pair $(f, v)$ can be deformed into a pair with no cusps. For Morse data with $\d_3X = \emptyset$, both $\d_2^+X$ and $\d_2^-X$ are collections of simple oriented loops in $\d_1X$, the orientation being induced by the orientation of $\d_1^+X$. 

Within the space $\mathcal W(X)$, no surgery on $\d_2X$, induced by deformations of Morse data $(f, v)$, is possible  (see Fig. 31 and 32). Indeed, the transversality of $v$ to $\d_2X$ is imposed by the nature of  $\mathcal W(X)$ and prevents loops from $\d_2X$ from touching each other, or being born/annahilated. Thus, each component of $\mathcal W(X)$ has its own oriented and polarized loop pattern $\theta(v) = \d_2^+X \coprod \d_2^-X \subset \d_1X$. \smallskip

{\it Questions:}\quad For a given $X$, what is the minimal number of positive/negative loops for nonsingular  Morse data with $\d_3X = \emptyset$? Evidently, $\chi(\d_1^+X) = -\chi(X)$ imposes constraints on the number of loops in $\d_2X$. Which oriented loop patterns are realizable on $X$?\smallskip

Within the space $\mathcal W{\Huge{\ast}}(X)$ no $\b$-moves are permitted. Therefore,
there is a well-defined map from $\pi_0(\mathcal W{\Huge{\ast}}(X))$ to skew-symmetric integral-valued 
bilinear forms $\Psi$.  In a sense, the forms are induced by the intersections of  1-cycles forming $\d_2^+X$ in the $v$-orbit space $X/\sim_{v}$. Since $X/\sim_{v}$ has singularities of codimension one, in general, this intersection has no homological interpretation. However, within the constraints of a given camber of $\mathcal W{\Huge{\ast}}(X)$, it is well-defined.
Consider a free $\Z$-module $M$ generated by the oriented loops $\g_i$ from $\d_2^+(X)$. Define $\Psi(\g_i, \g_j)$ to be the sum of  $\pm 1$ which are contributed by the $Q$-singularities  $x_{ij} \in \d_1^+X$ that correspond to the double-tangent trajectories linking  $\g_i$  to $\g_j$. The sign contributed by $x_{ij}$ is the 
$\oplus / \ominus$ polarity that has been associated with $x_{ij}$. Evidently, the form $\Psi$ is preserved under the $\a$-moves---the only admissible transformations within a given camber of $\mathcal W{\Huge{\ast}}(X)$ (see Fig. 29, 30). 

If $\d_3X = \emptyset$, then, for each component $\d_1X_j$ of the boundary, the degree $\#(\d_3^+X_j) - \#(\d_3^-X_j)$  of the map $h_j :  \d_2X_j \to S^1$ is zero.  By Lemma 4.3, we get $\chi(\d_1^+X) = \chi(\d_1^-X)$. This property is shared by all Morse data from $\mathcal W(X)$ (cf. [E1]).
\smallskip

The considerations above produce

\begin{thm} 
The  oriented and polarized loop patterns $\theta(v) \subset \d_1X$  are locally constant  on the space $\mathcal W(X)$ of Morse data with the property $\d_3 X = \emptyset$.

The  skew-symmetric form $\Psi$  and  the linking number $$lk_v(\d[K], \d[K]) = \#(Q(K)^\oplus) - \#(Q(K)^\ominus),$$  $K = K(f, v)$, are locally constant  on the subspace 
$\mathcal W{\Huge{\ast}}(X) \subset \mathcal W(X)$ of Morse data $(f, v)$ with no triple-tangent trajectories. 
\end{thm}

\begin{cor}
If $(f, v)$ and $(\tilde f, \tilde v)$ belong to different chambers of $\mathcal W{\Huge{\ast}}(X)$, then any generic path $\g$ that links in $\mathcal W(X)$ the point $(f, v)$ with the point $(\tilde f, \tilde v)$ must have at least $|gc_\ominus^\oplus(f, v) - gc_\ominus^\oplus(\tilde f, \tilde v)|$ intersections with the walls $\mathcal W(X) \setminus \mathcal W{\Huge{\ast}}(X)$ of various chambers of $\mathcal W{\Huge{\ast}}(X)$, that is, the deformation family $\g$  must have at least $|gc_\ominus^\oplus(f, v) - gc_\ominus^\oplus(\tilde f, \tilde v)|$ members with triple-tangent trajectories.
\end{cor}

Most likely, these and many other results of the paper admit multidimensional generalizations.

\end{document}